\documentclass[12pt]{amsart}
\usepackage{amssymb,amsmath,mathrsfs}
\usepackage[latin1]{inputenc}
\usepackage{graphicx}
\usepackage[all]{xy}
\DeclareMathOperator{\op}{op}
\DeclareMathOperator{\Op}{op}

\DeclareMathOperator{\Gr}{{\mathcal G}}

\DeclareMathOperator{\id}{id}
\DeclareMathOperator{\ind}{ind}

\DeclareMathOperator{\im}{im}
\DeclareMathOperator{\ch}{ch}
\DeclareMathOperator{\coker}{coker}

\DeclareMathOperator{\supp}{supp}
\DeclareMathOperator{\Td}{Td}
\DeclareMathOperator{\tr}{tr}
\DeclareMathOperator{\Tr}{Tr}
\DeclareMathOperator{\sgn}{sgn}

\begin{document}
%BLACKBOARD LETTERS
\newcommand{\hhom}{{\mathbb H}}
\newcommand{\R}{{\mathbb R}}
\newcommand{\N}{{\mathbb N}}
\newcommand{\Z}{{\mathbb Z}}
\newcommand{\C}{{\mathbb C}}
\newcommand{\T}{{\mathbb T}}
\newcommand{\rn}{{\mathbb R}^n}

%CALLIGRAPHIC LETTERS
\newcommand{\cA}{{\mathcal A}}
\newcommand{\cB}{{\mathcal B}}
\newcommand{\cC}{{\mathcal C}}
\newcommand{\cD}{{\mathcal D}}
\newcommand{\cE}{{\mathcal E}}
\newcommand{\cF}{{\mathcal F}}
\newcommand{\cG}{{\mathcal G}}
\newcommand{\cH}{{\mathcal H}}
\newcommand{\cI}{{\mathcal I}}
\newcommand{\cJ}{{\mathcal J}}
\newcommand{\cK}{{\mathcal K}}
\newcommand{\cN}{{\mathcal N}}
\newcommand{\cL}{{\mathcal L}}
\newcommand{\cP}{{\mathcal P}}
\newcommand{\cQ}{{\mathcal Q}}
\newcommand{\cR}{{\mathcal R}}
\newcommand{\cS}{{\mathcal S}}
\newcommand{\cT}{{\mathscr T}}
\newcommand{\cW}{{\mathcal W}}

%NICER CALLIGRAPHIC LETTERS
\newcommand{\ccA}{\mathscr{A}}
\newcommand{\ccB}{\mathscr{B}}
\newcommand{\ccC}{\mathscr{C}}
\newcommand{\ccD}{\mathscr{D}}
\newcommand{\ccE}{\mathscr{E}}
\newcommand{\ccF}{\mathscr{F}}
\newcommand{\ccG}{\mathscr{G}}
\newcommand{\ccH}{\mathscr{H}}
\newcommand{\ccI}{\mathscr{I}}
\newcommand{\ccJ}{\mathscr{J}}
\newcommand{\ccK}{\mathscr{K}}
\newcommand{\ccN}{\mathscr{N}}
\newcommand{\ccL}{\mathscr{L}}
\newcommand{\ccP}{\mathscr{P}}
\newcommand{\ccQ}{\mathscr{Q}}
\newcommand{\ccS}{\mathscr{S}}
\newcommand{\ccT}{\mathscr{T}}
\newcommand{\ccW}{\mathscr{W}}

%GREEK LETTERS
\newcommand{\ga}{\alpha}
\newcommand{\gb}{\beta}
\renewcommand{\gg}{\gamma}
\newcommand{\gG}{\Gamma}
\newcommand{\gd}{\delta}
\newcommand{\eps}{\varepsilon}
\newcommand{\gve}{\varepsilon}
\newcommand{\gk}{\kappa}
\newcommand{\gl}{\lambda}
\newcommand{\gL}{\Lambda}
\newcommand{\tgo}{\widetilde{\go}}
\newcommand{\go}{\omega}
\newcommand{\gO}{\Omega}
\newcommand{\gvp}{\varphi}
\newcommand{\gt}{\theta}
\newcommand{\gT}{\Theta}
\renewcommand{\th}{\vartheta}
\newcommand{\gs}{\sigma}

%Gothic
\newcommand{\fA}{{\mathfrak A}}
\newcommand{\fT}{{\mathfrak T}}

%LINES
\newcommand{\ol}{\overline}
\newcommand{\ul}{\underline}

%DIVERSE
\newcommand{\Dbar}{D\hspace{-1.5ex}/\hspace{.4ex}}
\newcommand{\dbar}{d\hspace{-1.0ex}{}^-\hspace{-.4ex}}
\newcommand{\bracket}[1]{[#1]}
\newcommand{\ti}{\tilde}
\newcommand{\Pf}{{\em Proof.}~}
\newcommand{\Proof}{\noindent{\em Proof}}
\newcommand{\pitensor}{\hat\otimes_\pi}
\newcommand{\eproof}{{~\hfill$ \triangleleft$}}
\renewcommand{\i}{\infty}
\newcommand{\rand}[1]{\marginpar{\small  #1}}
\newcommand{\forget}[1]{}
\newcommand{\Mat}[1]{\begin{pmatrix}#1\end{pmatrix}}
\newcommand{\TEXT}[1]{\marginpar{\hspace{1.3cm} \small  #1}}	
\newcommand{\cut}{C^\infty_{tc}(T^-X)}
\newcommand{\Ctc}{C^\infty_{c}}
\newtheorem{leer}{\hspace*{-.3em}}[section]
\newenvironment{rem}[2]%
{\begin{leer} \label{#1} {\bf Remark. }  %\marginpar{\hspace{1.3cm} \tiny #1} 
 {\rm #2 } \end{leer}}{}
\newenvironment{lemma}[2]%
{\begin{leer} \label{#1} {\bf Lemma. }  %\marginpar{\hspace{1.3cm}  \tiny #1}  
{\sl #2} \end{leer}}{}
\newenvironment{thm}[2]%
{\begin{leer}\label{#1} {\bf Theorem. }  %\marginpar{\hspace{1.3cm} \tiny #1} 
{\sl #2} \end{leer}}{}
\newenvironment{dfn}[2]%
{\begin{leer}  \label{#1} {\bf Definition. } %\marginpar{\hspace{1.3cm} \tiny #1}    
{\rm #2 } \end{leer}}{}
\newenvironment{cor}[2]% 
{\begin{leer} \label{#1} {\bf Corollary. }%\marginpar{\hspace{1.3cm} \tiny #1}   
{\sl #2 } \end{leer}}{}
\newenvironment{prop}[2]%
{\begin{leer} \label{#1} {\bf Proposition. } %\marginpar{\hspace{1.3cm} \tiny #1} 
{\sl #2} \end{leer}}{}
\newenvironment{extra}[3]%
{\begin{leer} \label{#1} {\bf #2. } %\marginpar{\hspace{1.3cm} \tiny #1}  
{\rm #3 } \end{leer}}{}
\renewcommand{\theequation}{\thesection.\arabic{equation}}
\renewcommand{\labelenumi}{{\rm (\roman{enumi})}}

%LISTS
\newcounter{num}
\newcommand{\bli}[1]{\begin{list}{{\rm(#1{num})}%
\hfill}{\usecounter{num}\labelwidth1cm
\leftmargin1cm\labelsep0cm\rightmargin1pt\parsep0.5ex plus0.2ex minus0.1ex
\itemsep0ex plus0.2ex\itemindent0cm}}
\newcommand{\eli}{\end{list}}

%EXTRAS
\def\Im{{\rm Im}\,}
\def\lra{\longrightarrow}
\def\Re{{\rm Re}\,}
\def\rpb{\overline\R_+}
\def\sumj#1{\sum_{j=0}^{#1}}
\def\sumk#1{\sum_{k=0}^{#1}}
\def\vect#1#2#3{\begin{array}{c}#1\\#2\\#3\end{array}}
\def\vec#1#2{\begin{array}{c}#1\\#2\end{array}}
\def\skp#1{\langle#1\rangle}

%%%%%%%%%%%%%%%%%%%%%%%%%%%%%%%%%%%%%%%%%%%%%%%%%%%%%%%%%%%%%%%%%%%%%%%
\title[Index Theory for Boundary Value Problems]%
{Index Theory for Boundary Value Problems via Continuous Fields 
of C*-algebras}
%%%%%%%%%%%%%%%%%%%%%%%%%%%%%%%%%%%%%%%%%%%%%%%%%%%%%%%%%%%%%%%%%%%%%%%
\author{Johannes Aastrup}
\address{SFB 478 ``Geometrische Strukturen'',
Hittorfstraße 27, 
48149 Münster,
Germany}
\email{johannes.aastrup@uni-muenster.de}
%%%%%%%%%%%%%%%%%%%%%%%%%%%%%%%%%%%%%%%%%%%%%%%%%%%%%%%%%%%%%%%%%%%%%%%
\author{Ryszard Nest}
\address{Department of Mathematics, 
Copenhagen University, 
Universitetsparken 5,
2100 Copenhagen, Denmark}
\email{rnest@math.ku.dk}
%%%%%%%%%%%%%%%%%%%%%%%%%%%%%%%%%%%%%%%%%%%%%%%%%%%%%%%%%%%%%%%%%%%%%%%
\author{Elmar Schrohe}
\address{Institut für Analysis,
Leibniz Universität Hannover, Welfengarten 1, 30167 Hannover,
Germany} 
\email{schrohe@math.uni-hannover.de} 
%%%%%%%%%%%%%%%%%%%%%%%%%%%%%%%%%%%%%%%%%%%%%%%%%%%%%%%%%%%%%%%%%%%%%%%
\begin{abstract}
We prove an index theorem for boundary value problems in Boutet de
Monvel's calculus on a compact manifold $X$ with boundary. 
The basic tool is the tangent semigroupoid $\ccT^-X$ generalizing
the tangent groupoid defined by Connes in the boundaryless case, 
and an associated
continuous field $C_r^*(\ccT^-X)$ of $C^*$-algebras over $[0,1]$.
Its fiber in $\hbar=0$, $C^*_r(T^-X)$, can be identified with the symbol
algebra for Boutet de Monvel's calculus; 
for $\hbar\not=0$ the fibers are isomorphic
to the algebra $\cK$ of compact operators. 
We therefore obtain a natural
map $K_0(C^*_r(T^-X))=K_0(\cC_0(T^*X))\to K_0(\cK)=\Z$.
Using deformation theory we show that this is the analytic index map.
On the other hand, using ideas from noncommutative geometry, 
we construct the topological index map and prove that it coincides 
with the analytic index map.
\end{abstract}
%%%%%%%%%%%%%%%%%%%%%%%%%%%%%%%%%%%%%%%%%%%%%%%%%%%%%%%%%%%%%%%%%%%%%%%
\maketitle
%%%%%%%%%%%%%%%%%%%%%%%%%%%%%%%%%%%%%%%%%%%%%%%%%%%%%%%%%%%%%%%%%%%%%%%
\section*{Introduction}
%%%%%%%%%%%%%%%%%%%%%%%%%%%%%%%%%%%%%%%%%%%%%%%%%%%%%%%%%%%%%%%%%%%%%%%
Let $X$ be a smooth compact manifold with boundary 
$\partial X$.
An operator in Boutet de Monvel's calculus on  $X$   is a matrix 
\begin{equation}\label{g.2.1}
A = \left( \begin{array}{cc} P^+ + G & K \\ T & S \end{array} \right) :
\begin{array}{ccc} \cC^\i (X,E_1) && 
\cC^\i (X,E_2) \\ \oplus & \to & \oplus
\\ \cC^\i (\partial X, F_1) && \cC^\i (\partial X, F_2) \end{array} .
\end{equation}
of operators acting on smooth sections of (hermitian) vector bundles 
$E_1$ and $E_2$ over $X$ and 
$F_1$ and $F_2$ over $\partial X$.

Here $P$ is a pseudodifferential operator on the double $\widetilde X$
of $X$, and $P^+$ is the so-called truncated operator given by
$P^+=r^+Pe^+$, where $e^+:\cC^\infty(X,E_1)\to L^2(\widetilde X,E_1)$ 
denotes extension by zero and $r^+$ is the operator of restriction of 
distributions on $\widetilde X$ to the open interior $X^\circ$ of $X$. 

For a general pseudodifferential operator $P$, the truncation  
$P^+$ will map $\cC^\infty(X,E_1)$ to $\cC^\infty(X^\circ, E_2)$,
but the result may not be smooth up to the boundary. 
The operator $P$ is supposed to satisfy the  
transmission condition to ensure mapping property \eqref{g.2.1}.
%; for details see Section  \ref{BdM}, below.

The entry $G$ is a so-called singular Green operator. 
Roughly speaking it acts like an operator-valued pseudodifferential 
operator along the boundary with values in smoothing operators in the 
normal direction. 
In the  interior, $G$ is  regularizing.
Singular Green operators come up naturally: 
If $P$ and $Q$ are pseudodifferential operators on $X$, then
the composition of the associated truncated operators differs
from the truncation of the composition 
by the so-called leftover term 
\begin{eqnarray}
L(P,Q)= (PQ)^+-P^+Q^+,
\end{eqnarray}
which is a singular Green operator. 
%Again, details will be given in Section \ref{BdM}.

The operators $T$ and $K$ are trace and potential (or Poisson) operators,
respectively, and
\forget{Classical examples of trace operators are given by the evaluation 
of the normal derivatives of $u\in \cC^\infty(X)$ at $\partial X$.
Another example for a trace operator is an integral operator
with a smooth kernel on $\partial X\times X$.
A typical example for a potential operator is the classical 
Poisson operator for the Dirichlet problem. 
Finally }
$S$ is a pseudodifferential operator on $\partial X$.
%We shall skip details on $T,K$, and $S$ as 
We skip details here, 
since for the purpose of index theory it is sufficient to consider the case
where there are  no bundles over the boundary and 
the operator $A$ is of the form 
$$A=P^++G$$
with both $P$ and $G$ classical.
Moreover we can confine ourselves to  the case where $A$ is of order and class
zero.

In analogy to the classical Lopatinskij-Shapiro condition,
the Fredholm property is governed by the invertibility of two 
symbols.
The first, the interior symbol of $A$, 
simply is the principal symbol of $P$.
The second, the {boundary symbol}, 
is an operator-valued symbol on $S^*\partial X$ which we will explain, below.

The above calculus was introduced by L. Boutet de Monvel in 1971, 
\cite{BdM71}. He showed that the operators of the form \eqref{g.2.1} indeed
form an algebra under composition, assuming the vector bundles match.
Moreover, this calculus contains the parametrices to elliptic elements 
and even their inverses whenever they exist.

He also proved an index theorem: To every elliptic operator $A$ one can 
associate a class $[A]$ in $K_c(T^*X^\circ)$, 
the compactly supported $K$-theory of the 
cotangent bundle over the interior $X^\circ$ of $X$. 
The index  of $A$ then is given by the topological index of $[A]$.
It is far from obvious how to assign $[A]$ to $A$. Boutet 
de Monvel gave a brilliant construction, combining pseudodifferential
analysis and classical topological $K$-theory in an ingenious way.
Still, the proof is hard for readers not familiar with
the details of the calculus, and efforts have been made to make it 
accessible to a wider audience. 

Following work by Melo, Nest, and Schrohe \cite{MNS}, it has been shown
in \cite{MSS} by Melo, Schick and Schrohe that the mapping $A\mapsto [A]$
can be obtained more easily with the help of $C^*$-algebra $K$-theory,
which had not yet been developed in 1971. This proof relies on only 
a rudimentary knowledge of the calculus and the ideal structure of the 
algebra of operators of order and class zero, but it lacks the geometric
intuition of Boutet de Monvel's initial idea.  

In this paper we will make the link to geometry. 
We will show how the index of a Fredholm operator of 
order and class $0$ can be determined
from its two symbols with the help of deformation theory and a 
continuous field of $C^*$-algebras over $[0,1]$.

This field, $C^*_r(\ccT^- X)$, is the reduced $C^*$-algebra of the 
tangent semigroupoid $\ccT^- X$ associated to $X$; 
it was introduced in \cite{ANS}.
The construction of $\ccT^- X$ is similar to 
Connes' construction of the tangent groupoid for a closed manifold,
cf.\ \cite[Section II.5]{Connes}.
In the case at hand, the `half-tangent space' 
$T^-X$ is glued to  
$X\times X\times {]0,1]}$. 
Here, fixing a connection,  
$T^\pm X$ consists of all tangent vectors $(x,v)$ such that 
$\exp_x(\pm tv)\in X$ for small $t\ge0$, and 
$X\times X\times {]0,1]}$ is endowed with 
the pair groupoid structure; the gluing is performed with the 
help of the exponential map: 
$(x,v,\hbar)\mapsto (x,\exp_x(-\hbar v),\hbar)$.

We showed in \cite{ANS} that -- just as in the boundaryless case -- 
the fibers of $C_r^*(\ccT^- X)$ are isomorphic to the algebra $\cK$ of compact 
operators for $\hbar>0$ so that their $K$-theory is given by $\Z$.
For $\hbar=0$, the situation is different. 
The `symbol algebra' $C^*_r(\ccT^-X)(0)=C^*_r(T^-X)$ is generated
by two representations of $\cC^\infty_c(T^-X)$. The first is the 
representation on $L^2(X^\circ)$ via convolution in the fibers.
The second takes into account the boundary: We associate to 
$f\in \cC_c^\infty (T^-X)$ the operator $\pi^\partial_0(f)$ on 
$L^2(T^+X|_{\partial X})$ given by half-convolution:
$$\pi_0^\partial (f)\xi(x,v)=\int_{T^+X}f(x,v-w)\xi(w)\,dw.$$
The $K$-theory of  $C^*_r(T^-X)$ 
turned out to be given by $K_0(\cC_0(T^*X))=K_c(T^*X^\circ)$.

%Our argument follows a line of thought developed by
%A.\ Connes \cite[Section 2.6]{Connes} and Elliott-Natsume-Nest \cite{ENN}.

In order to make use of this, 
we first compose the operator with an order reducing operator of 
positive order $m>\dim X$. 
This gives us an operator $A$ of order $m$ and class $0$ in Boutet de 
Monvel's calculus (in fact, we might have started with $A$ %an operator of
this order and class). 

We denote by $p^{{m}}$ its principal pseudodifferential symbol, 
a homogeneous function on the nonzero vectors in the cotangent bundle, and 
by $c^{{m}}$ its homogeneous principal boundary symbol,
a homogeneous operator-valued function defined outside the 
zero section in the cotangent bundle over the boundary.

We then set out to compute the index of $A$. To this end we first
smooth out both symbols near the zero section in the corresponding
cotangent bundle, obtaining a smooth function $p$
on $T^*X$ and a smooth boundary symbol operator $c$ on 
$T^*\partial X$. 
 
Following an idea of Elliott-Natsume-Nest \cite{ENN} 
we then consider a semiclassical deformation 
$A_\hbar,\ 0<\hbar\le 1$ of $A$ with 
$A=A_1$ and study the associated graph projection 
\begin{eqnarray}\label{0.a}
\cG(A_\hbar) = 
\Mat{(1+A_\hbar^*A_\hbar)^{-1}& 
(1+A_\hbar^*A_\hbar)^{-1}A_\hbar^*\\ 
A_\hbar(1+A_\hbar^*A_\hbar)^{-1}& 
A_\hbar(1+A_\hbar^*A_\hbar)^{-1}A_\hbar^*}.  
\end{eqnarray}

The positivity of the order of $A$ is crucial here; it ensures that
all entries of this matrix are compact operators, except for the
one in the lower right corner which differs from a compact operator
by the identity.
The graph projection of $A_\hbar$ therefore is a projection 
in $\cK^\sim$, the unitization of the compact operators. 
It is closely related to the 
index of $A$. In fact, denoting by $e$ the 
projection $\Mat{0&0\\0&1}$ (we will use the same notation
for this projection in various algebras), 
we will show in Theorem \ref{2.5} that
\begin{eqnarray}\label{0.b}
[\Gr(A)]-[e]=
[\pi_{\ker A}]-[\pi_{\coker A}]
\end{eqnarray}
is the difference of the classes associated to the 
projections onto the kernel and the co\-kernel, respectively, 
thus the index of $A$.

Using formula \eqref{0.a} 
we can also define the graph projections for $p$ and $c$. 

The crucial step is to show that  
$$\hbar\mapsto \left\{\begin{array}{ll}
\Gr(A_\hbar),&0<\hbar\le 1\\
\Gr(p)\oplus \Gr(c),&\hbar=0
\end{array}\right.$$
defines a continuous section of the unitization of 
$C_r^*(\ccT^- X)$, cf.\ Proposition \ref{t.7}, 
Theorem \ref{t.8}.

A technical problem arises from the fact that 
Boutet de Monvel's calculus does not in general contain the 
adjoints of operators of positive order. The analysis 
of the graph projection therefore takes us out of the 
calculus. We overcome this difficulty by working with 
operator-valued symbol classes. Boutet de 
Monvel's calculus fits well into this concept, 
cf.~Schrohe and Schulze \cite{SS1}, Schrohe \cite{SIB}; 
moreover, it also allows to treat the adjoints.

The continuity of the section gives us a natural map
associating to an elliptic operator $A$ a class
in $K_0(C^*_r(T^-X))=K_0(\cC_0(T^*X^\circ))=K_c(T^*X^\circ)$
by evaluating the section in 
$\hbar=0$, more precisely by taking 
$$\left([\Gr(p)\oplus \Gr(c)]\right)-
\left( [e\oplus e]\right). $$
 
%More is true: Since $C_r^*(\ccT^- X)(0)=C^*_r(T^-X)$ and
%$C_r^*(\ccT^- X)(1)\cong \cK$, 
In addition, evaluation in $\hbar=0$ and in $\hbar=1$ 
defines a map in $K$-theory
$$\ind_a:K_0(\cC_0(T^*X^\circ))%= K_0(C^*_r(T^-X))
\to K_0(\cK);$$
it associates to  
$\left([\Gr(p)\oplus \Gr(c)]\right)-
\left( [e\oplus e]\right) $
the class  $[\Gr(A)]-[e]$, thus the index of $A$. 

In this way we obtain the analytic index map. In 
a second step we then construct the topological index
in order to obtain an index formula in cohomological terms.

A cohomological index formula had been established before by
Fedosov \cite[Chapter 2, Theorem 2.4]{Fedosov96}: He showed that
$$\ind A = \int_{T^*X}\ch([p])\Td(X)+ 
\int_{T^*\partial X}\ch'([c])\Td(X),$$
with the Chern character of  the $K$-theory class of $p$ and a variant
$\ch'$ of the Chern character of the $K$-class  of $c$. 
As usual, $\Td$ denotes the  Todd class. %ES not genus 
%In fact, Fedosov states this formula

In this article, we conclude the discussion in the spirit of noncommutative
geometry. We extend the fundamental class 
$$\int_{T^*X^\circ}:
H^*_c(T^*X^\circ)=HP^*(\cC^\infty_c(T^*X^\circ)\to \C$$
to a fundamental class
$$F:HP^*(\cC^\infty_{tc} (T^-X))\to \C.$$
Here 
$$\cC^\infty_{tc} (T^-X)=
\cC_c^\infty(TX)\oplus 
\cC^\infty_c(T\partial X\times\R_+\times\R_+)$$ 
is the `smooth' algebra associated to
the symbol algebra $C^*_r(T^-X)$
which was introduced in \cite[Definition 2.16]{ANS}. 

The construction leaves us a certain degree of freedom. In fact,
we obtain extensions $F_\go$ for every choice $\go$ of a 
closed form on $T^*X$ of even degree, which is the pull-back 
of a closed form on $X$. 

Using the equality of the analytical and the topological index
in the boundaryless case, established by Connes, we then obtain the 
index formula 
$$\ind A = F_{\Td (X)} (\ch( [\Gr(p)\oplus \Gr(c)]-
 [e\oplus e]))
$$
with the Chern-Connes character $\ch$.

%%%%%%%%%%%%%%%%%%%%%%%%%%%%%%%%%%%%%%%%%%%%%%%%%%%%%%%%%%%%%%%%%%%%%%%
\section{Operators in Boutet de Monvel's Calculus}
\label{BdM}
\setcounter{equation}{0}
%%%%%%%%%%%%%%%%%%%%%%%%%%%%%%%%%%%%%%%%%%%%%%%%%%%%%%%%%%%%%%%%%%%%%%%

\extra{B.12}{Manifolds with boundary}{%
In the sequel let $X$ be a compact manifold with boundary and
$\widetilde X$ its double.
\forget{We choose a covering of $X$ by open sets 
$X_j$ and associated coordinate maps $\chi_j:X_j\to\R^n$. 
We can assume that a coordinate neighborhood either is 
contained in a collar neighborhood of the boundary (`boundary chart') 
or else
does not intersect a neighborhood of the boundary (`interior chart').
It is also no restriction to suppose that 
the boundary charts all lie in a collar
neighborhood of $\partial X$ and that coordinate changes for boundary 
charts involve only the tangential variables.
}%forget
\forget{It will be sometimes convenient to know that, 
by going over to a refinement, 
we may also assume that for each pair $(j,k)$, both $X_j$ 
and $X_k$ lie in the 
support of a single coordinate map. 
We can assume that the boundary charts all lie in a collar
neighborhood of $\partial X$ and that coordinate changes for boundary 
charts involve only the tangential variables.}
We use the standard Sobolev space 
$H^s(\widetilde X)$ and the associated spaces  
$$H^s(X)=\{u|_{X^\circ}:u\in H^s(\widetilde X)\},\quad 
H^s_0(X)=\{u\in H^s(\widetilde X): \supp u\subseteq X\}.$$ 
The space $\cC^\infty(X)$ is dense in $H^s(X)$ for all $s$, while 
$\cC^\infty_c(X^\circ)$ is dense in $H^s_0(X)$ for all $s$.
The $L^2$-inner product allows us to identify $H^s(X)$ with the 
dual of $H^{-s}_0(X)$. 

In general, the spaces $H^s(X)$ and $H^s_0(X)$ are quite different.
For $-1/2<s<1/2$, however, they can be identified. 
In particular, we have for $s=0$ the
natural identification of $L^2(X)$ with the subset of all 
functions in $L^2(\widetilde X)$ which vanish on 
$\widetilde X\setminus X$.
}
%%%%%%%%%%%%%%%%%%%%%%%%%%%%%%%%%%%%%%%%%%%%%%%%%%%%%%%%%%%%%%%%%%%%%%%

%%%%%%%%%%%%%%%%%%%%%%%%%%%%%%%%%%%%%%%%%%%%%%%%%%%%%%%%%%%%%%%%%%%%%%%

\extra{B.15}{Operators, symbols, ellipticity}{%
Detailed descriptions of Boutet de Monvel's calculus were given  
by Grubb \cite{Grubb} and Rempel and Schulze 
\cite{RS}.
In order to overcome technical difficulties 
%it will be convenient in the later sections to 
we will also rely on the representation of these operators
by operator-valued symbols as presented in \cite{SIB}.	

In order to keep the exposition short, we will focus on the
algebra formed by the elements in the upper left corner:

$$\cA=\{A:\cC^\infty(X,E_1)\to \cC^\infty(X,E_2): A=P^++G \},$$
where $P$ is a pseudodifferential operator with the transmission
property and  $G$ a singular Green operator (sGo).
We assume all operators to be classical.

The operator $A$ is said to have order $\mu$ 
and class $d\in \N_0$, if  $P$ is of order $\mu$ and $G$ is of 
order $\mu$ and class $d$.
We speak of smoothing or regularizing operators, if the order 
is $-\infty$.

A sGo $G$ of order $m$ and class $d$ can be written 
$$G=\sum_{j=0}^d G_j\partial_\nu^j$$
where $G_j, j=0,\ldots d$, are sGo's of order $\mu-j$ and
class $0$ and $\partial_\nu$ is a differential operator
which coincides with the normal derivative in a neighborhood
of the boundary and vanishes farther away from the boundary.
This representation differs from the standard one in that it 
avoids the trace operators. The equivalence becomes clear
from \ref{1.5}, below. 

Let $\gvp\in \cC_c^\infty(X^\circ)$,
and denote by $M_\gvp$ multiplication by $\gvp$. Then
$G M_\gvp$ is regularizing of class $0$ and
$M_\gvp G$ is regularizing of class $d$.

We associate two symbols to $A$. 
The first is the pseudodifferential principal symbol, $\gs^\mu_\psi(A)$.
It is defined as the principal symbol  of $P$, restricted to 
$T^*X\setminus 0$:
$$\gs_\psi^\mu(A) =\left.\gs^\mu(P)\right|_{T^*X\setminus 0}.$$ 
This makes sense as $G$ is
smoothing in $X^\circ$ and therefore cannot contribute 
to the symbol.

The second is the principal boundary symbol, $\gs^\mu_\partial(A)$, 
defined on $T^*\partial X\setminus 0$. 
In local coordinates near $\partial X$, it is given by
$$\gs_\partial^\mu(A)(x',\xi') 
= p^\mu(x',0, \xi', D_n)^++g^\mu(x',\xi', D_n): 
\cS(\ol\R_+, \tilde E_1)\to \cS(\ol\R_+, \tilde E_2).$$
Here $p^\mu$ and $g^\mu$  are the homogeneous principal symbols of 
$P$ and $G$, respectively. 
For fixed $x',\xi'$, the operator $p^\mu(x',0,\xi',D_n)$ 
is the Fourier 
multiplier with symbol $p^\mu(x',0,\xi',\xi_n)$, while
$g^\mu(x',\xi',D_n)$ 
is an integral operator with smooth integral kernel.
By $\tilde E_j$ we have denoted the fiber in $(x',\xi')$ 
of the pullback of $E_j$ to $T^*X$.

The principal boundary symbol is homogeneous 
on $T^*\partial X\setminus 0$ in a sense we shall explain later,
see \eqref{twistedhom}, so that it can be viewed as a function 
on $S^*\partial X$.

An operator $A$ of order $\mu$ and class $d\le\max\{\mu,0\}$
is said to be elliptic, if 
\bli\roman\item $\gs^\mu_\psi(A)(x,\xi): \pi^*E_1\to \pi^*E_2$ 
is invertible for all $(x,\xi)\in T^*X\setminus 0$, and 
\item $\gs^\mu_\partial(A)(x',\xi'): 
H^\mu(\R_+,\widetilde E_1)\to L^2(\R_+,\widetilde E_2)$
is invertible for all $(x',\xi')\in T^*\partial X\setminus 0$. 
\eli 
Here $\pi:T^*X\to X$ is the base point projection.

Apart from these symbols we have, of course, in any coordinate 
neighborhood, the full symbols of $P$ and  $G$ in the corresponding 
classes. 
%the singular Green
%symbol kernels of the operators associated with $G$.
} 

%%%%%%%%%%%%%%%%%%%%%%%%%%%%%%%%%%%%%%%%%%%%%%%%%%%%%%%%%%%%%%%%%%%%%%%
%ES
\thm{B.18}{Let $A$ be an operator of order $\mu$ and class $d$ in 
Boutet de Monvel's calculus. Then $A$ induces a bounded linear map 
$$A:H^s(X,E_1) \to H^{s-\mu}(X,E_2)$$
for each $s>d-1/2.$
}
In general we cannot extend $A$ to $H^s(X,E_1)$ 
for $s\le d-1/2$. The reason is that neither extension by zero 
makes sense on these spaces nor do integral operators with smooth 
(up to the boundary) integral kernels act continuously on them. 
%when $s\le-1/2$. 

%%%%%%%%%%%%%%%%%%%%%%%%%%%%%%%%%%%%%%%%%%%%%%%%%%%%%%%%%%%%%%%%%%%%%%%

\thm{B.20}{ 
Let $A_1:\cC^\infty(X,E_1)\to \cC^\infty(X,E_2)$ and 
$A_2=\cC^\infty(X,E_2)\to \cC^\infty (X,E_3)$ 
be operators of orders $\mu_1$ and $\mu_2$ and classes
$d_1$ and $d_2$, respectively, 
in Boutet de Monvel's calculus on $X$. 
The composition
$A_2A_1$ is an operator of order $\mu_1+\mu_2$ and class
$\max\{d_1, \mu_1+d_2\}$.
Its principal symbols are given by 
\begin{eqnarray*}
\gs^{\mu_1+\mu_2}_\psi(A_2A_1)
&=&\gs^{\mu_2}_\psi(A_2)\  \gs^{\mu_1}_\psi(A_1);\\
\gs^{\mu_1+\mu_2}_\partial(A_2A_1)
&=&\gs^{\mu_2}_\partial(A_2)\ \gs^{\mu_1}_\partial(A_1).
\end{eqnarray*} 
%\eli
}

%%%%%%%%%%%%%%%%%%%%%%%%%%%%%%%%%%%%%%%%%%%%%%%%%%%%%%%%%%%%%%%%%%%%%%%

\thm{B.22}{Let $A$ be an operator of order $\mu$ and class 
$d\le\max\{\mu,0\}$ in 
Boutet de Monvel's calculus and $s>d-1/2$. Then
$$A:H^s(X,E_1) \to H^{s-\mu}(X,E_2)$$
is Fredholm if and only if $A$ is elliptic.
In this case we find an operator $B$ of order $-\mu$ and class
$d'\le\max\{-\mu,0\}$ such that 
\begin{eqnarray*}
AB=I+R_1\quad\text{and}\quad BA=I+R_2
\end{eqnarray*}
with regularizing operators $R_1$ and $R_2$ of class $d'$
and $d$, respectively. For the symbols we have 
\begin{eqnarray*}
\gs_\psi^{-\mu}(B)=\gs_\psi^{\mu}(A)^{-1}\quad\text{and}\quad 
\gs_\partial^{-\mu}(B)=\gs_\partial^{\mu}(A)^{-1}.
\end{eqnarray*}  
 }
%%%%%%%%%%%%%%%%%%%%%%%%%%%%%%%%%%%%%%%%%%%%%%%%%%%%%%%%%%%%%%%%%%%%%%%
\extra{B.21}{Adjoints}{%
Boutet de Monvel's calculus is not closed under taking adjoints. 
For the sake of completeness let us introduce a few basic concepts. 
To an operator $A$ of order $\mu$ and class $d\le\max\{0,\mu\}$ 
we can associate a minimal adjoint
$A^*_{\min}$ defined on $\cC^\infty_c (X^\circ)$ (we omit the bundles from
the notation) and taking values
in $\cC^\infty(X)'$ by the relation 
$$\skp{Au,v}=\skp{u,A_{\min}^*v},
\quad u\in \cC^\infty(X), v\in \cC^\infty_c(X^\circ).$$
For $s>d-1/2$, the adjoint $A^*$  of the bounded operator 
$$A: H^s(X)\to H^{s-\mu}(X)$$
is then given by extending 
$A_{\min}^*$  by continuity to an operator 
$A^*:H^{\mu-s}_0(X)\to H^{-s}_0(X).$

If the class is zero, we can explicitly determine the adjoint:
We write $A=P^++G$ with a pseudodifferential operator
$P$ on $\widetilde X$ and a 
singular Green operator $G$ of order $\mu$ and 
class $0.$ 
Next we recall from \cite[Lemma 1.3.1]{Grubb} 
that a pseudodifferential operator
$P$ of order $\mu\in\Z$ with the transmission property can be written
\begin{eqnarray}\label{B.4.1}
P=S+Q
\end{eqnarray}
where $S$ is a differential operator of order $\le \mu$ and $Q$ is a 
pseudodifferential operator of order $\mu$ 
which maps $e^+\cC^\infty(X)$, the extensions
(by zero) of smooth functions on $X$ to $\widetilde X$,  to 
$L^2(X)$ and satisfies -- with the formal adjoint $Q_f^*$ --
\begin{eqnarray}\label{B.4.2}
\skp{Q^+u,v}=\skp{u,Q_f^{*,+}v}, \quad u,v\in \cC^\infty_c(X).
\end{eqnarray}
Note that for $\mu\le 0$, 
equality \eqref{B.4.2} will hold for $P$. 
Still it it useful to know that we can choose $S$ in such a way that
the local symbol $q$ of $Q$ is of order  (at most) $-1$ with respect 
to $\xi_n$; we say that $q$ is of normal order  $-1$. 

For $u,v\in \cC^\infty(X)$ it is known 
\cite[Section 1.6]{Grubb} that  
$$\skp{S^+u,v}_{L^2(X)}=\skp{u,S_f^{*,+}v}_{L^2(X)}
+\skp{\fA \rho^+ u,\rho^+ v}_{L^2(\partial X)}$$
with the formal adjoint $S^*_f$ of $S$, 
the Green matrix $\fA$, and the vectors 
$\rho^+ u$ and $\rho^+v$ of boundary values
for $u$ and $v$. 
When $v$ lies in $\cC^\infty_c(X^\circ)$, 
the boundary terms on the right hand side vanish, so that 
$$\skp{S^+u,v}_{L^2(X)}=\skp{u,S_f^{*,+}v}_{L^2(X)},
\quad u\in \cC^\infty(X),
v\in \cC^\infty_c(X^\circ).$$
On $\cC^\infty_c(X^\circ)$, the operation $e^+$ of extending by 
zero is trivial; moreover $S^*_f$ is a differential operator,
so that $S_f^*v\in\cC^\infty_c(X^\circ)$. There are no singular
terms arising at the boundary, and $S^{*,+}_fv= S^*_fv$
as a functional on $\cC^\infty(X)$.  
For $Q$ the corresponding
identity \eqref{B.4.2} is valid by construction. 
Again, $Q^{*,+}_fv=Q^*_fv$ as a functional on $\cC^\infty(X)$.
Hence 
\begin{eqnarray}\label{P*}
\skp{P^+u,v}_{L^2(X)}=\skp{u,P_f^{*}v}_{L^2(X)},
\quad u\in \cC^\infty(X),
v\in \cC^\infty_c(X^\circ).
\end{eqnarray}

As the singular Green part $G$ is assumed to be of class
zero, it also has a formal adjoint $G^*_f$,
cf.\ \cite[(1.2.47)]{Grubb}, and
\begin{eqnarray}\label{G*}
\skp{Gu,v}_{L^2(X)}=\skp{u,G_f^*v}_{L^2(X)}.
\end{eqnarray}

Hence, as a functional on $\cC^\infty(X)$,  
$$A^*_{\min}v = P_f^{*}v+G^*_fv, \quad v\in \cC^\infty_c(X^\circ).$$

If  $A=P^++G$ is of order  
$\mu\le0$ and class $0$, 
then $A: L^2(X)\to L^2(X)$ is continuous, and 
$\cC^\infty_c(X^\circ)$ is dense in $L^2(X)$. 
We conclude from  \eqref{P*} and 
\eqref{G*} that the $L^2$-adjoint of 
$A$ is 
$$A^*= P^{*,+}_f + G_f^*$$
with the formal adjoints $P^*_f$ and $G_f^*$ of $P$ and
$G$, extended to $L^2$.
Thus $A^*$ again is an operator in Boutet de Monvel's calculus.
Its principal symbols are given by 
\begin{eqnarray*}
\gs^{\mu_1}_\psi(A^*)&=&\gs^{\mu_1}_\psi(A)^*;\\
\gs^{\mu_1}_\partial(A^*)&=&\gs^{\mu_1}_\partial(A)^*.
\end{eqnarray*}
}

\cor{B.21b}{%ES
Suppose $A:H^\mu(X,E_1)\to L^2(X,E_2)$ is 
invertible of order $\mu\ge0$ and class $d\le \mu$. 
Then $A^{-1}$ is an operator of order $-\mu$ and 
class $0$, 
cf.\ \cite[Theorem 4.5]{S2}.
It has an adjoint $(A^{-1})^*$ in Boutet de Monvel's
calculus,  which  extends to a bounded operator
from $H^{-\mu}_0(X,E_2)$ to $L^2(X,E_1)$, considering $e^+$ as a 
trivial operation on $H^{-\mu}_0$.
On the other hand, the minimal adjoint $A^*_{\min}$ of $A$ 
extends to an invertible operator $A^*:L^2(X,E_2)\to H^{-\mu}_0(X,E_1)$.
It is easily checked that $(A^{-1})^*=(A^*)^{-1}.$
}

\extra{or.1}{Order reducing operators}{There exists a family 
$\gL^m_-$, $m\in\Z$, of classical scalar pseu\-do\-differential 
operators on $X$ satisfying the transmission condition with the 
following properties:
\bli\roman\item $(\gL^{m}_-)^+:H^s(X)\to H^{s-m}(X)$ is an
isomorphism for all $s>-1/2$. 
\item In fact, the result of the application of $r^+\gL^m_-$ 
to $u$ in $H^s(\widetilde X)$ only depends on the restriction
of $u$ to $X^\circ$. The map in  (i) extends
to an isomorphism  $r^+\gL^m_-e_s:H^s(X)\to H^{s-m}(X)$ 
for an arbitrary choice of an extension operator $e_s:
H^s(X)\to H^s(\widetilde X)$.  
\item $(\gL^{m}_-)^+(\gL^\mu_-)^+=(\gL^{m+\mu}_-)^+$, $m,\mu\in\Z$.
\item The (extension of the) formal pseudodifferential adjoint 
defines an isomorphism 
$\gL^m_+=((\gL^{m}_-)^{+})^*:H^s_0(X)\to H^{s-m}_0(X)$.
%ES Johannes
\item The operators $\gL^m_-,\gL^m_+$ can be extended to operators with the same properties but acting in a vector bundle $E$. %ES 
\eli 
}
%%%%%%%%%%%%%%%%%%%%%%%%%%%%%%%%%%%%%%%%%%%%%%%%%%%%%%%%%%%%%%%%%%%%%%%

%%%%%%%%%%%%%%%%%%%%%%%%%%%%%%%%%%%%%%%%%%%%%%%%%%%%%%%%%%%%%%%%%%%%%%%
\section{Operator-valued Symbols}
\label{edge}
\setcounter{equation}{0}
%%%%%%%%%%%%%%%%%%%%%%%%%%%%%%%%%%%%%%%%%%%%%%%%%%%%%%%%%%%%%%%%%%%%%%%
It will be helpful to consider the operators in Boutet de Monvel's
calculus as operator-valued pseudodifferential operators.
We recall the basic concepts from \cite{SIB}.

We first fix a function $\R^n\ni\xi\mapsto \bracket{\xi}\in\R_{\ge0}$
which is positive for $\xi\not=0$ and coincides with $|\xi|$ for
$|\xi|\ge1$.  
\extra{1.2}{Group actions}{%
A strongly continuous group action on a Banach space $E$ is a family
$\gk = \{ \gk_\gl: \gl \in \R_+ \}$ of isomorphisms in $\cL(E)$ 
such that
$\gk_\gl \gk_\mu = \gk_{\gl\mu}$ and the mapping 
$\gl \mapsto \gk_\gl e$
is continuous for every $e \in E$.
Note that there is an $M>0$ such that 
\begin{eqnarray}\label{estgk}
\left\|\gk_\gl\right\|\le \left(\max\{\gl, \gl^{-1}\}\right)^M.
\end{eqnarray}  
For the usual Sobolev spaces on
$\R$ and $\R_+$ 
we shall use the group action defined on functions $u$ by
\begin{equation}\label{121}
(\gk_\gl u)(x) = \gl^{1/2} u(\gl x).
\end{equation}
It will be useful to consider also weighted Sobolev spaces: 
For $s=(s_1,s_2)\in\R^2$ we define 
\begin{eqnarray*}
H^{s}(\R)=H^{(s_1,s_2)}(\R)=\{\bracket{x}^{-s_2} u:u\in H^{s_1}(\R)\}
\end{eqnarray*} 
with the usual unweighted space on the right hand side. 
Similarly we define $H^{s}(\R_+)$. 
We then have 
$\cS(\R_+)=\text{projlim}_{s_1,s_2\to\infty}H^{s_1,s_2}(\R_+)$ and
$\cS'(\R_+)=\text{indlim}_{s_1,s_2\to\infty}H^{-s_1,-s_2}(\R_+).$
On $E=\C^l, l \in \N$, we use the trivial group
action $\gk_\gl \equiv {\rm id}$. Sums of spaces of the above kind will
be endowed with the sum of the group actions.
}

%%%%%%%%%%%%%%%%%%%%%%%%%%%%%%%%%%%%%%%%%%%%%%%%%%%%%%%%%%%%%%%%%%%%%%%

\extra{1.4}{Operator-valued symbols and amplitudes}{
Let $E,F$ be Banach spaces with strongly continuous group actions $\gk$
and $\ti\gk$, respectively. 
Let $a \in \cC^\i (\R^q \times \R^q \times
\R^q, \cL(E,F))$ and $\mu \in \R$. 
We shall write 
$a \in S^\mu (\R^q\times \R^q \times \R^q; E,F)$ 
and call $a$ an amplitude of order $\mu$ provided that, 
for all multi-indices
$\ga,\gb,\gg$, there is a constant $C = C(\ga,\gb,\gg)$ with
$$
\| \ti\gk_{\bracket{\eta}^{-1}} D^\ga_\eta D^\gb_y D^\gg_{\ti y}
a(y,\ti y,\eta) \gk_{\bracket{\eta}} \|_{\cL(E,F)} \leq C \,
\bracket{\eta}^{\mu-|\ga|} .
$$
If $a$ is independent of $y$ or $\tilde y$ 
we shall write $a \in S^\mu (\R^q
\times \R^q; E,F)$. 

For $E=F=\C$ we recover the usual pseudodifferential symbol classes.

The concept extends to the cases where $E$ is an inductive or 
$F$ a projective limit.

In order to avoid lengthy formulas we shall abbreviate this by saying
that $a$ is a symbol of order $\mu$ with values in $\cL(E,F)$.  

\forget{We call an operator-valued symbol $a$ of order $\mu$ 
classical, if it has an asymptotic expansion 
$a\sim\sum a_j$ (i.e., each $a_j$ is a symbol of order $\mu-j$
and $a-\sum_{j<N}a_j$ is of order $\mu-N$), where $a_j$ 
is positively homogeneous for large $\eta$ in the sense that
\begin{eqnarray}\label{twistedhom}
a_j(y,\lambda\eta) =
\lambda^{\mu-j}\tilde{\kappa}_\lambda\,a_j(y,\eta)\,\kappa_{\lambda^{-1}},
\quad x\in\R^q,|\eta|\ge R, \gl\ge1.
\end{eqnarray}
}
}

%%%%%%%%%%%%%%%%%%%%%%%%%%%%%%%%%%%%%%%%%%%%%%%%%%%%%%%%%%%%%%%%%%%%%%%

\extra{1.4a}{Example: Potential, trace and singular Green boundary 
symbol operators}{%
The boundary symbol operators associated to potential, trace or 
singular Green symbols in the usual presentation of Boutet de Monvel's
calculus have simple descriptions in the framework of operator-valued
symbols. 
\bli\alph\item The elements in 
$S^\mu(\R^{n-1}\times\R^{n-1};\C,\cS(\R_+))$ are precisely the 
boundary symbol operators associated with 
potential symbols of order
\footnote{In fact, the notion of order 
differs slightly in \cite{Grubb}, \cite{SIB} and \cite{SS1}; this will, however, 
not play a role in the sequel. } 
$\mu$ on $\ol\R^n_+$.
\item The elements of $S^\mu(\R^{n-1}\times\R^{n-1};\cS'(\R_+),\C)$
are precisely the boundary symbol operators associated with 
trace symbols of order $\mu$ and class $0$.
Those of the form 
$$\sum_{j=0}^dt_j(x',\xi')D_n^j$$
with $t_j$ in $S^{\mu-j}(\R^{n-1}\times\R^{n-1};\cS'(\R_+),\C)$
and the derivative $D_n$ on $\R$
are the  boundary symbol operators associated with 
trace symbols of order $\mu$ and class $d$ on $\ol\R^n_+$.

\item The elements of 
$S^\mu(\R^{n-1}\times\R^{n-1};\cS'(\R_+),\cS(\R_+))$
are precisely the boundary symbol operators 
associated with singular Green symbols of order $\mu$ 
and class $0$ on $\ol\R^n_+$. Those  of the form 
$$\sum_{j=0}^d g_j(x',\xi')D_n^j$$
with $g_j$ in $S^{\mu-j}(\R^{n-1}\times\R^{n-1};\cS'(\R_+),\cS(\R_+))$
are the singular Green boundary symbol operators 
of order $\mu$ and class $d$.
\eli
We therefore speak of these operator-valued symbols as potential,
trace and singular Green boundary symbol operators or, for short, 
symbols, of the corresponding orders and classes.}

%%%%%%%%%%%%%%%%%%%%%%%%%%%%%%%%%%%%%%%%%%%%%%%%%%%%%%%%%%%%%%%%%%%%%%%

\extra{1.5}{Example: Trace operators}{%
Let $\gg_j$  be defined on $\cS(\R_+)$  by
$
\gg_ju = \lim_{t \to 0^+} D^j_{n} u(t).
$
It extends to an
element of $\cL(H^{\gs} (\R_+), \C)$ for  $\gs=(\gs_1,\gs_2)$, 
$\gs_1 > j +1/2$, by the trace theorem for Sobolev spaces. 

Viewed  as an operator-valued symbol
independent of the variables $y$ and $\eta$, $\gg_j$ 
then is a symbol of order $j+1/2$ with values in 
$\cL(H^{\gs} (\R_+), \C)$:
Recalling that the group action on the Sobolev space is given by
\eqref{121} while on $\C$ it is given by the identity, we only 
have to check that
$$
\| \gg_j \gk_{\bracket\eta} \|_{\cL(H^\gs(\R_+),\C)} 
= O(\bracket{\eta}^{j+1/2}) .
$$
This is immediate, since
$$
\partial^j_t \left(\bracket{\eta}^{1/2} 
u(\bracket{\eta} t)\right)
= \bracket{\eta}^{j+1/2} (\partial^j_t u) (\bracket{\eta} t) .
$$
Let us next show that $\gg_j$ is a trace symbol of order $j+1/2$ 
and class $j+1$ in the sense of \ref{1.4a}. It clearly suffices to 
do this for $j=0$.

Choose $\gvp\in \cS(\R_+)$ with $\gvp(0)= 1$. The identity
\begin{eqnarray}\label{gg0}
\ u(0) =- \int_0^\infty \bracket{\xi'}^{}\gvp'(\bracket{\xi'}s)u(s)\,ds
-\int_0^\infty  \gvp(\bracket{\xi'}s)\partial_{s}u(s)\,ds, \ 
u\in\cS(\R_+),
\end{eqnarray}
shows that
$$\gg_0 = t_0 + it_1D_n,$$
where $t_0$ and $t_1$ are the operator-valued symbols of order $1/2$ and $-1/2$, 
respectively, with values in $\cL(\cS'(\R_+), \C)$, 
given by  $$t_0 u = -\int_0^\infty
\bracket{\xi'}^{}\gvp(\bracket{\xi'}s) u(s)\,ds,\quad t_1 u = -\int_0^\infty
\gvp(\bracket{\xi'}s) u(s)\,ds.$$ 
Hence $\gg_0$ is of class $1$.}

%%%%%%%%%%%%%%%%%%%%%%%%%%%%%%%%%%%%%%%%%%%%%%%%%%%%%%%%%%%%%%%%%%%%%%%

\dfn{1.6}{%
For $a \in S^\mu (\R^q \times \R^q \times \R^q; E,F)$, the
pseudodifferential operator
$$
\op a : \cS(\R^q,E) \to \cS(\R^q,F)
$$
is defined by
$$
\left(\op a\right) u(y) 
= \iint e^{i(y-\ti y) \eta} a(y,\ti y,\eta) u(\ti y)
\, d \ti y \, \dbar \eta ; \qquad y \in \R^q .
$$
Here 
$\dbar \eta = (2\pi)^{-q} d \eta$. 
If $a$ is independent of $\ti y$,
this reduces to 
\begin{eqnarray}\label{leftsymbol}
\left(\op a\right) u(y) = %(2\pi)^{-q/2}\ 
\int e^{iy\eta} a(y,\eta) \hat u(\eta) \, \dbar\eta;
\end{eqnarray}
 in this case we call $a$ a {\em left} symbol
for $\op a$. If $a$ is independent of $y$, then
\begin{eqnarray}\label{rightsymbol}
\left(\op a\right) u(y) = 
\iint e^{i(y-\ti y)\eta} a(\ti y,\eta) u(\ti y) \, d\ti
y \, \dbar \eta ,
\end{eqnarray}
and $a$ is called a {\em right} symbol.}

%%%%%%%%%%%%%%%%%%%%%%%%%%%%%%%%%%%%%%%%%%%%%%%%%%%%%%%%%%%%%%%%%%%%%%%

\extra{1.7}{Example: Action in the normal direction}{
Let $p \in S^\mu (\R^n \times \R^n)$. 
For fixed $(x',\xi')$, the
function $p(x',\cdot,\xi',\cdot)$ is an element of 
$S^\mu(\R \times \R)$. For $\gs \in \R^2$, 
$p(x',\cdot,\xi',\cdot)$
induces a bounded linear operator
$$
p(x',x_n,\xi',D_n)=  (\op_{x_n} p)(x',\xi'): 
H^{\gs} (\R) \to H^{(\gs_1-\mu,\gs_2)}(\R) ;
$$
by 
$$(\op_{x_n}p)(x',\xi')u(x_n)
=\int e^{ix_n\xi_n}p(x',x_n,\xi',\xi_n)\hat u(\xi_n)\, \dbar\xi_n,$$
see \cite[Theorem 1.7]{Schrohe90} for the boundedness on 
weighted spaces.
We then have
\begin{equation}\label{43.9.1}
\gk_{\bracket{\xi'}^{-1}} 
\left(\op_{x_n} p\right) \gk_{_{\bracket{\xi'}}} = \op_{x_n}p \left(
x', x_n/\bracket{\xi'}, \xi', \bracket{\xi'} \xi_n \right) :
\end{equation}
In fact, for $u \in \cS(\R)$,
\begin{eqnarray*}
\lefteqn{\gk_{\bracket{\xi'}^{-1}} \left(\op_{x_n} p\right) 
(\gk_{\bracket{\xi'}} u)(x_n) }\\
&=&  
\int e^{ix_n \xi_n / \bracket{\xi'}} \bracket{\xi'}^{-1}
p(x',x_n/\bracket{\xi'},\xi',\xi_n) 
\hat u (\xi_n/\bracket{\xi'}) \dbar\xi_n ;
\end{eqnarray*}
and the substitution $\eta_n = \xi_n / \bracket{\xi'}$ 
yields the assertion.
The theorem, below, shows that $\op_{x_n} p$ is 
an operator-valued symbol in the sense of 
\ref{1.4}:
}

%%%%%%%%%%%%%%%%%%%%%%%%%%%%%%%%%%%%%%%%%%%%%%%%%%%%%%%%%%%%%%%%%%%%%%%

\prop{1.8}{In the above situation we have
%For $p \in S^\mu(\R^n \times \R^n)$ we have
$$
\op_{x_n} p \in S^\mu (\R^{n-1} \times \R^{n-1}; H^{\gs}(\R),
H^{(\gs_1-\mu,\gs_2)}(\R)).
$$
}

%%%%%%%%%%%%%%%%%%%%%%%%%%%%%%%%%%%%%%%%%%%%%%%%%%%%%%%%%%%%%%%%%%%%%%%

\Proof.
Given multi-indices $\ga,\gb$, we have to estimate
\begin{eqnarray*}
\lefteqn{\sup_{x',\xi'}\|\bracket{\xi'}^{|\ga|} 
\gk_{\bracket{\xi'}^{-1}} 
(\op_{x_n} (D^\ga_{\xi'}D^\gb_{x'} p))
\gk_{\bracket{\xi'}} \|_{{\cL}(H^{(\gs_1,\gs_2)}(\R),
H^{(\gs_1-\mu,\gs_2)}(\R))}}\\
&=& \sup_{x',\xi'}\|\bracket{\xi'}^{|\ga|} 
\op_{x_n} (D^\ga_{\xi'} D^\gb_{x'} p) \left( x',
    x_n/\bracket{\xi'},\xi',\xi_n \bracket{\xi'} 
\right) \|_{{\cL}(H^{(\gs_1,\gs_2)}(\R),
H^{(\gs_1-\mu,\gs_2)}(\R))}.
\end{eqnarray*}
%where both norms are in ${\mathcal L}(H^\gs(\R), H^{\gs-(\mu,0)}(\R))$. 
Since $D^\ga_{\xi'}D^\gb_{x'} p$ is of order $\mu-|\ga|$ 
we may assume that $\ga = \gb = 0$. 
Now 
$$\Op_{x_n}p(x',x_n/\bracket{\xi'},\xi',\xi_n \bracket{\xi'}):
H^{(\gs_1,\gs_2)}(\R)\to
H^{(\gs_1-\mu,\gs_2)}(\R)) $$
is continuous, and a bound for its norm is given by the suprema 
$$\sup\left\{|D^\ga_{\xi_n}D^\gb_{x_n}\{p(x',
x_n/\bracket{\xi'},\xi',\xi_n \bracket{\xi'})\}|\bracket{\xi_n}^{-\mu}
:x_n,\xi_n\in\R\right\}$$
for a finite number of derivatives. 
Since each of them is  $O(\bracket{\xi'}^\mu)$ 
the proof is complete.
\eproof

%%%%%%%%%%%%%%%%%%%%%%%%%%%%%%%%%%%%%%%%%%%%%%%%%%%%%%%%%%%%%%%%%%%%%%%

\extra{1.10b}{Example: Multiplication operators}{\bli \alph\item
Multiplication $M_{x_n}$ by $x_n$  is an element of 
$S^{-1}(\R^{n-1}\times\R^{n-1}, H^s(\R_+), H^{s+(0,1)}(\R_+))$, $s\in\R^2$.
We may replace the pair of Sobolev spaces by
$(H_0^s(\ol\R_+), H_0^{s+(0,1)}(\ol\R_+))$.
% and $(H^s(\R_+), H^{s+(0,1)}(\R_+))$.
\item
Let $\gvp\in \cC_b^\infty(\ol\R_+)$ vanish to all orders at $0$.
Then multiplication $M_\gvp$ by $\gvp(x_n)$ is an element of 
$S^{-\infty}(\R^{n-1}\times\R^{n-1}, H^s(\R_+), H^{s}(\R_+))$. 
 \eli}

%%%%%%%%%%%%%%%%%%%%%%%%%%%%%%%%%%%%%%%%%%%%%%%%%%%%%%%%%%%%%%%%%%%%%%%

\Proof. (a) follows from the fact that $M_{x_n}$ has the symbol
$a(x',\xi')=1\otimes M_{x_n}$ and that 
$$\gk_{\bracket{\xi'}}^{-1}\,a(x',\xi')\,\gk_{\bracket{\xi'}u}
=\bracket{\xi'}^{-1} M_{x_n}u,\quad u\in \cC^\infty_c(\R).$$

(b) follows from (a).\eproof
%%%%%%%%%%%%%%%%%%%%%%%%%%%%%%%%%%%%%%%%%%%%%%%%%%%%%%%%%%%%%%%%%%%%%%%

%%%%%%%%%%%%%%%%%%%%%%%%%%%%%%%%%%%%%%%%%%%%%%%%%%%%%%%%%%%%%%%%%%%%%%%

\extra{1.11}{Asymptotic summation and classical symbols}{A 
sequence $(a_j)$ of operator-valued symbols of orders $\mu-j$ 
with values in $\cL(E,F)$ can be summed asymptotically to a symbol
$a$ of order $\mu$, and $a$ is unique modulo symbols of order
$-\infty$.

A symbol $a$ of order $\mu$ is said to be {\em classical,}
if it has an asymptotic expansion $a \sim\sum_{j=0}^\infty a_j$ with
$a_j$ of order $\mu-j$ 
satisfying the homogeneity relation
\begin{equation}\label{twistedhom}
a_j(y,\tilde y,\lambda\eta) = \lambda^{\mu-j
}\tilde{\kappa}_\lambda\,
a_j(y,\tilde y, \eta)\,\kappa_{\lambda^{-1}}
\end{equation}
for all $\lambda \ge 1,|\eta| \ge R$ with a suitable constant $R$. 
We write
$a\in S^\mu_{cl}(\R^q\times\R^q\times\R^q;E,F)$. 
For $E=\C^k,~F=\C^l$ we recover the standard 
notion.
}

%%%%%%%%%%%%%%%%%%%%%%%%%%%%%%%%%%%%%%%%%%%%%%%%%%%%%%%%%%%%%%%%%%%%%%%

The key to many results on compositions is the lemma, below, 
which is adapted from 
Kumano-go \cite[Chapter 2, Lemma 2.4]{K}.

%%%%%%%%%%%%%%%%%%%%%%%%%%%%%%%%%%%%%%%%%%%%%%%%%%%%%%%%%%%%%%%%%%%%%%%

\lemma{1.13a}{Let $a\in S^\mu(\R^q\times\R^q\times\R^q;E,F)$. 
For $|\gt|\le 1$ define $a_\gt$ by the oscillatory integral
\begin{eqnarray}
a_\gt(y,\eta)=\iint e^{-iz\zeta} a(y,y+z,\eta+\gt\zeta)dz\dbar\zeta.
\end{eqnarray}
Then the family $\{a_\gt:|\gt|\le 1\} $ is uniformly bounded 
in $S^\mu$; its
seminorms can be estimated by those for $a$. 
%Here we used Kumano-go's notion of oscillatory integrals, cf.\ 
%{\rm \cite[Chapter 1]{K}}, 
%extended to operator-valued symbols.
}

%%%%%%%%%%%%%%%%%%%%%%%%%%%%%%%%%%%%%%%%%%%%%%%%%%%%%%%%%%%%%%%%%%%%%%%

\thm{1.14}{\bli\alph\item 
Let $a \in S^\mu (\R^q \times \R^q \times \R^q; E,F)$. 
Then there is a
{\rm(}unique{\rm)} left symbol $a_L=a_L(y,\eta)$ 
for $\op a$ acting as in \eqref{leftsymbol} and 
a $($unique$)$ right symbol $a_R=a_R(\ti y,\eta)$ 
acting as in \eqref{rightsymbol}.
They are given by the oscillatory integrals
\begin{eqnarray}
a_L(y,\eta)
&=&\iint 
e^{-iy\eta}a(y,y+z,\eta+\zeta)dz\dbar\zeta\text{ and}\\
a_R(\ti y,\eta)
&=&\iint e^{i\ti y\eta}a(\ti y+z,\ti y,\eta+\zeta)dz\dbar\zeta
\end{eqnarray}

Moreover, we have 
\begin{eqnarray}\label{leftsymbol1}
a_L(y,\eta)
&=&\sum_{|\ga|<N}
\left.\frac1{\ga!}\partial_\eta^\ga D_{\ti y}^\ga 
a(y,\ti y,\eta)\right|_{\ti y=y}
+ 
N\sum_{|\gg|=N}\int_0^1 \frac{(1-\gt)^{N-1}}{\gg!}\nonumber\\
&&\times
\iint e^{-iz\zeta}\partial_\eta^\gg D_{\ti y}^\gg 
a(y, y+z,\eta+\gt\zeta)dz\dbar\zeta d\gt
\text{\ and}
\label{al}
\\
a_R(\ti y,\eta)&=&
\sum_{|\ga|<N}
\left.\frac{(-1)^{|\ga|}}{\ga!}\partial_\eta^\ga 
D_{y}^\ga a(y,\ti y,\eta)\right|_{y=\ti y}
+ 
N\sum_{|\gg|=N}\int_0^1  
\frac{(1-\gt)^{N-1}}{\gg!}\nonumber\\
&&\times
\iint e^{+iz\zeta}\partial_\eta^\gg 
D_{\ti y}^\gg a(y, y+z,\eta+\gt\zeta)dz\dbar\zeta d\gt
\label{ar}
\end{eqnarray}
with remainders $N\sum_{|\gg|=N}\ldots$ in 
$S^{\mu-N}(\R^q\times\R^q;E,F)$. 

\item 
Given $a \in S^\mu (\R^q \times \R^q; E,F)$ and $b \in S^{\ti\mu} (\R^q
\times \R^q; F,G)$ there is a left symbol $c \in S^{\mu+\ti\mu} (\R^q
\times \R^q; E,G)$ such that
$$
\op b \circ \op a = \op c .
$$
As usual we write $c=b\#a$. 
We have the asymptotic expansion formula
\begin{eqnarray}\label{1.14.1}
(b\#a)(y,\eta)&=&\sum_{|\ga|<N}\frac1{\ga!}
\partial^\ga_\eta b(y,\eta)
D^\ga_y a(y,\eta)+
N\sum_{|\gg|=N}\int_0^1  
\frac{(1-\gt)^{N-1}}{\gg!}\nonumber\\
&&\times
\iint e^{+iz\zeta}\partial_\eta^\gg 
b(y,\eta+\gt\zeta)D_{y}^\gg a( y+z,\eta)
dz\dbar\zeta d\gt
\end{eqnarray}
with a remainder in $S^{\mu+\tilde\mu-N}$. 
\eli}

%%%%%%%%%%%%%%%%%%%%%%%%%%%%%%%%%%%%%%%%%%%%%%%%%%%%%%%%%%%%%%%%%%%%%%%

\Proof.
The proof for existence and form of the left symbol 
is analogous to that of \cite[Chapter 2, Theorem 2.5]{K}. 
The right symbol is obtained by a simple modification.
The estimates on the remainder follow from 
Lemma \ref{1.13a}.

For the analysis of the composition  let $a_R$ be the right symbol for
$\op a$. Then
$$\op b\circ\op a = \op  b \circ \op a_R =
\op \ti c
$$
with $\ti c(y,\ti y,\eta) = b(y,\eta) a_R (\ti y,\eta)$. 
Choosing $c$ as the left symbol of $\op \ti c$ gives the assertion.
Formula \eqref{1.14.1} follows from \eqref{al} and \eqref{ar}.
For the scalar case see \cite[Chapter 2, Theorems 2.6 and 3.1]{K}.
\eproof
%%%%%%%%%%%%%%%%%%%%%%%%%%%%%%%%%%%%%%%%%%%%%%%%%%%%%%%%%%%%%%%%%%%%%%%

\extra{1.30}{Duality}{%
Let $(E_-, E_0, E_+)$ 
be a triple of Hilbert spaces. We assume that all
are embedded in a common vector space $V$ and that $E_0 \cap E_+ \cap
E_-$ is dense in $E_\pm$ as well as in $E_0$. 
Moreover we assume that
there is a continuous, non-degenerate sesquilinear form 
$\skp{\cdot,\cdot}_E :
E_+ \times E_- \to \C$ which coincides with the inner product of $E_0$
on $(E_+\cap E_0)\times (E_-\cap E_0)$.
We ask that, via
$\skp{\cdot,\cdot}_E$, we may identify $E_+$ with the dual of 
$E_-$ and vice
versa, and that
$$
\| e \|'_{E_-} = \sup_{\| f \|_{E_+} =1} | \skp{f,e}_E | , \quad \| f
\|'_{E_+} = \sup_{\| e \|_{E_-} =1} | \skp{f,e}_E |
$$
furnish equivalent norms on $E_-$ and $E_+$, respectively.
Suppose there is a group action $\gk$ on $V$
which has strongly continuous restrictions to 
$E_0$ and $E_\pm$, unitary on $E_0$,
i.e.,  
$\skp{\gk_\gl e,f}_E = \skp{e,\gk_{\gl^{-1}} f}_E$ for $e,f \in E_0$. 
Then
$$
\skp{\gk_\gl e,f}_E = \skp{e, \gk_{\gl^{-1}} f}_E , 
\quad e \in E_+ , f \in E_- ,
$$
since the identity holds on the dense set 
$(E_+ \cap E_0) \times (E_- \cap
E_0)$. 
In other words, the action
$\gk$ on $E_+$ is  dual to the action $\gk$ on
$E_-$ and vice versa.

Typical examples for the above situation are given by the triples
of weighted Sobolev spaces
$$(H^{-\gs} (\R), L^2 (\R), H^\gs (\R))\ \mbox{~~ and~~}\  
(H^{-\gs}_0 (\ol \R_+),L^2(\R_+), H^\gs (\R_+)), 
\quad \gs \in \R^2.$$

Let $(F_-,F_0,F_+)$ be an
analogous triple of Hilbert spaces with group action $\ti \gk$, 
and let
$a \in S^\mu (\R^q \times \R^q \times \R^q; E_-, F_-)$. 
We define $a^*$ by
$a^* (y,\ti y,\eta) = a(\ti y,y,\eta)^* \in \cL (F_+,E_+)$, where the
last asterisk denotes the adjoint operator with respect to the
sesquilinear forms $\skp{\cdot,\cdot}_E$ and $\skp{\cdot,\cdot}_F$.
%:$$\skp{a(\ti y,y,\eta)^* f,e}_E 
%= \skp{f, a(\ti y,y,\eta) e}_F , \quad e \in E_- ,f \in F_+ .$$
It is not difficult to check that 
$a^* \in S^\mu (\R^q \times \R^q \times\R^q; F_+, E_+)$.

Moreover, we may introduce a continuous non-degenerate 
sesquilinear form
$$
\skp{\cdot,\cdot}_{\cS_E} : 
\cS (\R^q, E_+) \times \cS (\R^q, E_-) \to \C
$$
by $\skp{u,v}_{\cS_E} = \int \skp{u(y), v(y)}_E dy$. 
Analogously we define
$\skp{\cdot,\cdot}_{\cS_F}$.

The symbol $a^*$ induces a continuous mapping 
$\Op a^*: \cS (\R^q, F_+) \to \cS(\R^q, E_+)$. 
This is the unique operator satisfying
$$
\skp{(\Op a^*) u,v}_{\cS_E} = \skp{u, (\Op a) v}_{\cS_F} .
$$
}

%%%%%%%%%%%%%%%%%%%%%%%%%%%%%%%%%%%%%%%%%%%%%%%%%%%%%%%%%%%%%%%%%%%%%%%
%ES superfluous
\forget{\cor{1.30a}{Let $a \in S^\mu (\R^q \times \R^q; E,F)$.
The (left) symbol of the formal adjoint to 
$\op a$ is the left symbol associated to
$a(\tilde y,\eta)^*$. 
Its asymptotic expansion is given by \eqref{leftsymbol1}.}
}%forget
%%%%%%%%%%%%%%%%%%%%%%%%%%%%%%%%%%%%%%%%%%%%%%%%%%%%%%%%%%%%%%%%%%%%%%%
\extra{1.30b}{Change of coordinates}{Let $\chi:\R^q\to\R^q$ be a
smooth diffeomorphism with all derivatives bounded and 
$0<c\le|\det\chi(y)|\le  C$ for all $y$.

Given an operator-valued pseudodifferential operator 
$A=\op a: \cS(\R^q,E)\to \cS(\R^q,F)$ we define its push-forward 
$A^\chi$ under $\chi$ by 
$$(A^\chi(u\circ \chi))(y) = (Au)(\chi(y)).$$
For any lattice in $\R^q$ one finds
functions $\{\gvp_j\}$, $\psi_j$, $j=1,2,\ldots$, 
each of them centered around a lattice point, such that $\sum\gvp_j=1$
and $\gvp_j\psi_j=\gvp_j$. 

Just as in the scalar case, cf.\ \cite[Chapter 2, \S 6]{K}, 
it turns out that $A^\chi$ 
is an operator-valued pseudodifferential operator. 
In fact, modulo regularizing operators,
$$A^\chi=\op a^\chi$$  with the
double symbol 
\begin{eqnarray}\label{chco}
\lefteqn{a^\chi(y,y',\eta) = \sum_{j=1}^\infty \gvp_j(\chi(y))\, 
a(\chi(y), \chi(y'),\nabla_y\chi(y,y')^{-t}\eta)\,\psi_j(\chi(y'))}\\
&&\times|\det\nabla_y\chi(y,y')|^{-1}\ 
|\det\partial_y\chi(y')|.\nonumber
\end{eqnarray}
Here, 
\begin{eqnarray}
\nabla_y\chi(y,y')=\int_0^1 \partial_y\chi(y+s(y-y')\,ds, 
\end{eqnarray}
the superscript $-t$ means the inverse of the 
transpose, and the above formulas only make sense, 
if the lattice is sufficiently 
fine.

The leading term of the left symbol of $a^\chi$ is 
given by 
\begin{eqnarray}\label{leading}
a^\chi(y,y,\eta)= 
\sum \gvp(\chi(y)) a(\chi(y), \chi(y),\partial_y\chi(y)^{-t}\eta),
\end{eqnarray}
noting that $\nabla_y\chi(y,y)=\partial_y\chi(y)$.
}

%%%%%%%%%%%%%%%%%%%%%%%%%%%%%%%%%%%%%%%%%%%%%%%%%%%%%%%%%%%%%%%%%%%%%%%

%%%%%%%%%%%%%%%%%%%%%%%%%%%%%%%%%%%%%%%%%%%%%%%%%%%%%%%%%%%%%%%%%%%%%%%
\extra{1.15}{Semiclassical operators}{We 
shall now consider families $a(\hbar)$ of symbols
with values in $\cL(E,F)$,  
which depend smoothly on the parameter $\hbar\in(0,1]$.
We define an $\hbar$-scaling with the help of the group actions
$\gk$ on $E$ and $\tilde\gk$ on $F$:  
\begin{eqnarray*}
a_\hbar(\hbar;y,\tilde y,\eta)
&=&\tilde\gk^{-1}_\hbar a(\hbar;y,\tilde y,\hbar\eta)\gk_\hbar
\ \ \text{ and}\\
a_{1/\hbar}(\hbar;y,\tilde y,\eta)
&=&\tilde\gk^{}_\hbar a(\hbar;y,\tilde y,\eta/\hbar)\gk^{-1}_\hbar 
.
\end{eqnarray*}
Then:
\bli\alph
\item Let $a=a(\hbar),\ \hbar\in(0, 1]$, be bounded in
$S^\mu(\R^q\times\R^q\times \R^q;E,F)$. Then so is 
$\{a_\hbar:\gve\le\hbar\le 1\}$ for every $\gve>0$.
Moreover, for every $\hbar_0>0$, the mapping $\hbar\mapsto a_\hbar$ 
is continuous in $\hbar_0$ with respect to the topology of $S^\mu$.

\item Let $a$ be as above and write 
$a|_{\rm diag}(\hbar;y,\eta)=a(\hbar;y,y,\eta)$. Then 
\begin{eqnarray}
(a_\hbar)_L-(a|_{\rm diag})_\hbar&=&\hbar\ b_L(\hbar)
\quad\text{ and}\nonumber\\
(a_\hbar)_R-(a|_{\rm diag})_\hbar&=&\hbar\ b_R(\hbar),\nonumber
\end{eqnarray}
with $b_{L,1/\hbar}$ and $b_{R,1/\hbar}$ 
bounded in the topology of $S^{\mu-1}$.
The corresponding seminorms can be  estimated in terms
of the symbol seminorms for $a$.

\item Given bounded families $a(\hbar)$ and $b(\hbar)$
in $S^{\mu}(\R^q\times \R^q;E,F)$ 
and $S^{\mu'}(\R^q\times \R^q;F,G)$,
we find, for each $N$
%let $c(\hbar)=b_\hbar\# a_\hbar$. 
%For each fixed $N$, we then find 
\begin{eqnarray}\label{1.15c.1}
b_\hbar\#a_\hbar -\sum_{|\ga|<N}\frac{\hbar^{|\ga|}}{\ga!}
(\partial^\ga_\eta b)_\hbar
(D^\ga_y a)_\hbar=\hbar^Nr_N(\hbar)
\end{eqnarray}
with $r_{N,1/\hbar} $ bounded in 
$S^{\mu+\mu'-N}(\R^q\times\R^q;E,G).$
The seminorms can be estimated in terms of those 
of $a$ and $b$.
\eli}

%%%%%%%%%%%%%%%%%%%%%%%%%%%%%%%%%%%%%%%%%%%%%%%%%%%%%%%%%%%%%%%%%%%%%%%

\Proof. (a) The first assertion is immediate from 
\eqref{estgk} and the fact that, for $\gve\le\hbar\le 1$, 
the quotient  $\bracket{\hbar\eta}/\bracket{\eta}$ 
is both bounded and bounded away from 
zero. The second follows from the fact that 
\begin{eqnarray*}
\lefteqn{a_\hbar(\hbar;y,\eta)-a_{\hbar_0}(\hbar_0;y,\eta)
=(a_\hbar(\hbar;y,\eta)-a_{\hbar}(\hbar_0;y,\eta))+ 
(a_{\hbar}(\hbar_0;y,\eta)-a_{\hbar_0}(\hbar_0;y,\eta))}
\\
&=&\left(\int_{0}^{1}
\partial_\hbar a(\hbar_0+s(\hbar-\hbar_0);y,\hbar\eta)\,ds\ + 
\int_0^1 \partial_\eta a(\hbar_0;y,\hbar_0\eta+\gt(\hbar-\hbar_0)\eta)\,
 d\gt\right)\ (\hbar-\hbar_0).
\end{eqnarray*}
\\
(b) is immediate from 
\eqref{al} and \eqref{ar}, respectively, together with Lemma \ref{1.13a}.   
\\
(c) In the expansion formula 
\eqref{1.14.1}  let us replace there $a$ by $a_\hbar$ and $b$ by
$b_\hbar$. We have $\partial^\ga_\eta (b_\hbar)=
\hbar^{|\ga|} (\partial^\ga_\eta b)_\hbar$ and $D^\ga_y (a_\hbar)
=(D^\ga_y a)_\hbar$. This leads to the desired 
expansion. For the rescaled remainder 
$r^N_{1/\hbar}$ we obtain the expression in \eqref{1.14.1}
with $\partial^\gg_\eta b(y,\eta+\gt\zeta)$ replaced by
$(\partial^\gg_\eta b)(y,\eta+\hbar\gt\zeta)$.
The boundedness in $S^{\mu+\mu'-N}$ then follows from 
Lemma \ref{1.13a}.  
\forget{\begin{eqnarray*}
r_N(\hbar;y,\eta)&=&\hbar^{-N}
N\sum_{|\gg|=N}\int_0^1  
\frac{(1-\gt)^{N-1}}{\gg!}\nonumber\\
&&\times
\iint e^{+iz\zeta}\partial_\eta^\gg 
(b_\hbar(y,\eta+\gt\zeta))
D_{y}^\gg a_\hbar(y+z,\eta)
dz\dbar\zeta d\gt.
\end{eqnarray*}
We note that %\rand{hbar pr?fen}
$$\partial_\eta^\gg (b_\hbar(y,\eta+\gt\zeta))
=\hbar^N(\partial_\eta^\gg b_\hbar)(y, \eta+\gt\zeta).$$
The corresponding relation for $b$  in connection with 
Lemma \ref{1.13a} then shows the assertion. }
\eproof

%%%%%%%%%%%%%%%%%%%%%%%%%%%%%%%%%%%%%%%%%%%%%%%%%%%%%%%%%%%%%%%%%%%%%%%

\cor{1.15z}{Let $a$ be a pseudodifferential operator and $\chi$ 
a change of coordinates as in {\rm \ref{1.30b}}. 
It follows from \eqref{chco} and {\rm \ref{1.15}(b)} that
$$(a_\hbar)^\chi-(a^\chi)_\hbar=\hbar\op b(\hbar)
$$ 
with $b_{1/\hbar}$ bounded in the topology of $S^{\mu-1}$. 
}
%%%%%%%%%%%%%%%%%%%%%%%%%%%%%%%%%%%%%%%%%%%%%%%%%%%%%%%%%%%%%%%%%%%%%%%

\lemma{1.15a}{Let the supports of 
$\gvp,\psi\in \cC^\infty_b(\R^q)$ have positive distance,
and let $a\in S^\mu(\R^q\times\R^q;
E,F)$. Then for any $N$ we can write 
$$\gvp (\op a_\hbar)\psi=\hbar^N \op r_N(\hbar)$$
with a family $r_N$ of symbols such that 
$r_{N,1/\hbar}$ is bounded in $S^{\mu-N}$.  
 }
%%%%%%%%%%%%%%%%%%%%%%%%%%%%%%%%%%%%%%%%%%%%%%%%%%%%%%%%%%%%%%%%%%%%%%%
\Proof. This is immediate from the expansion formula \eqref{al}.\eproof
%%%%%%%%%%%%%%%%%%%%%%%%%%%%%%%%%%%%%%%%%%%%%%%%%%%%%%%%%%%%%%%%%%%%%%%

\lemma{1.15f}{Let $a=a(\hbar)$, $\hbar\in(0,1]$, be bounded in 
$S^0(\R^q\times\R^q;L^2(\R_+),L^2(\R_+))$.
Then also  $a_\hbar $ is  bounded in that class; its seminorms can be
estimated in terms of those of $a$.
}
%%%%%%%%%%%%%%%%%%%%%%%%%%%%%%%%%%%%%%%%%%%%%%%%%%%%%%%%%%%%%%%%%%%%%%%
\Proof. Since $\gk$ is unitary on $L^2$, we only have to check that 
$$%f\tilde\gk_{\bracket\eta}^{-1}\tilde\gk_\hbar^{-1}
D^\ga_\eta D^\gb_y(a(\hbar;y,\hbar\eta))=O(\bracket{\eta}^{-|\ga|}).$$
This in turn is a consequence of the fact that 
$$D^\ga_\eta D^\gb_y (a(\hbar;y,\hbar\eta))=\hbar^{|\ga|}
(D^\ga_\eta D^\gb_y a)(\hbar;y,\hbar\eta)$$
and that $\hbar\bracket\eta \bracket{\hbar\eta}^{-1}$ is bounded. \eproof
%%%%%%%%%%%%%%%%%%%%%%%%%%%%%%%%%%%%%%%%%%%%%%%%%%%%%%%%%%%%%%%%%%%%%%%

\prop{h.7a}{Let $s_1,s_2, t_1,t_2\ge 0$, $s_1+s_2>0$, $t_1+t_2>0$ and
$m\in\R$. Moreover,  let 
$g\in S^m(\R^q\times\R^q, H^{(-s_1,-t_1)}_0(\ol \R_+), 
H^{(s_2,t_2)}(\R_+))$ with 
$g(y,\eta)=0$ for large $|y|$. 
Then, for each $\gve>0$,  $g$ can be approximated by elements in 
$$S^{-\infty}(\R^q\times\R^q, \cS'(\R_+), \cS(\R_+)),$$ 
which vanish for large $|y|$, in the topology of 
$S^{m+\gve}= S^{m+\gve}(\R^q\times\R^q, 
L^2(\R_+),L^2(\R_+))$.
}
%%%%%%%%%%%%%%%%%%%%%%%%%%%%%%%%%%%%%%%%%%%%%%%%%%%%%%%%%%%%%%%%%%%%%%%
For the proof we need the following well-known result:
%%%%%%%%%%%%%%%%%%%%%%%%%%%%%%%%%%%%%%%%%%%%%%%%%%%%%%%%%%%%%%%%%%%%%%%

\lemma{h.7b}{Let $K\in\cK(L^2(X))$ and $\gve>0$. 
Then there exist $\gvp_1,\ldots,\gvp_N,\psi_1,\ldots,\psi_N\in 
\cC^\infty_c(X^\circ)$ such that   
$$\|K-\sum_{j=1}^N \gvp_j\otimes\psi_j\|_{L^2(X)}<\gve.$$
}
%%%%%%%%%%%%%%%%%%%%%%%%%%%%%%%%%%%%%%%%%%%%%%%%%%%%%%%%%%%%%%%%%%%%%%%
\Proof\ \ of Proposition \ref{h.7a}. 
By composing with the operator 
$\bracket{\eta}^{-m-\gve}\otimes \id$, 
we may assume that $m=-\gve$. 
Choose $\gvp\in\cC^\infty_c([0,\infty))$ 
with $\gvp(t)\equiv 1$ for $t\le 1$, and let 
$$ g_N(y,\eta)=g(y,\eta)\gvp(|\eta|/N),\quad N\in \N.$$ 
Then $g_N$ is a regularizing symbol with values in 
$\cL(H^{(-s_1,-t_1)}_0(\ol\R_+), H^{(s_2,t_2)}(\R_+))$.  
Moreover, $g_N$ tends to $g$ 
in the topology of $S^0$. 
It is therefore sufficient to approximate $g_N$ in $S^0$.

As $g_N$ vanishes for $(y,\eta)$ outside a compact set, 
we have in fact 
$$g_N\in \cS(\R^q\times\R^q)%\pitensor \cS(\R^q)
\pitensor \cL(H^{(-s_1,-t_1)}_0(\ol\R_+),H^{(s_2,t_2)}(\R_+))
\hookrightarrow 
\cS(\R^q\times\R^q)\pitensor\cK(L^2(\R_+)).$$
According to Lemma \ref{h.7b}, each element of $\cK(L^2(\R_+))$
can be approximated in $\cL(L^2(\R_+))$
by an integral operator with a rapidly decreasing kernel. 
%in $\cS(\R_+\times\R_+)$.  
Each of these defines a continuous operator
$\cS'(\R_+)\to \cS(\R_+)$. 
Hence $g_N$ can be approximated, in the topology of 
$$\cS(\R^q\times\R^q)\pitensor \cL(L^2(\R_+))
\hookrightarrow S^{-\infty}(\R^q\times\R^q;L^2(\R_+),L^2(\R_+))$$
by elements
in  the tensor product 
$\cS(\R^q\times\R^q)\pitensor \cS(\R_+\times\R_+)$, 
hence also  by elements  of 
$S^{-\infty}(\R^q\times\R^q;\cS'(\R_+),\cS(\R_+))$ which
vanish for $y$ outside a compact set. 
\eproof
%%%%%%%%%%%%%%%%%%%%%%%%%%%%%%%%%%%%%%%%%%%%%%%%%%%%%%%%%%%%%%%%%%%%%%%
\rem{h.7f}{A symbol $g$ in 
$S^{-\infty}(\R^q\times\R^q;\cS'(\R_+),\cS(\R_+))$ which
vanishes for $y$ outside a compact set induces an operator 
$\op a$ on $L^2(\R^{q+1}_+)$ with an integral kernel which is 
rapidly decreasing and thus can be approximated by a smooth
compactly supported integral kernel.
This slightly improves the statement of 
Proposition \ref{h.7a}.}
%%%%%%%%%%%%%%%%%%%%%%%%%%%%%%%%%%%%%%%%%%%%%%%%%%%%%%%%%%%%%%%%%%%%%%%

\section{Semiclassical Operators in Boutet de Monvel's Calculus}
\setcounter{equation}{0}
%%%%%%%%%%%%%%%%%%%%%%%%%%%%%%%%%%%%%%%%%%%%%%%%%%%%%%%%%%%%%%%%%%%%%%%

%%%%%%%%%%%%%%%%%%%%%%%%%%%%%%%%%%%%%%%%%%%%%%%%%%%%%%%%%%%%%%%%%%%%%%%

\extra{h.1}{Pseudodifferential operators}{Let 
$p\in S^\mu(\R^n\times\R^n)$ and $u\in\cS(\R)$. Then 
\begin{eqnarray}
\op_{x_n}(p_\hbar) u(x_n)
&=&  \int\int_0^\infty e^{i(x_n-y_n)\xi_n} 
p(x',x_n,\hbar\xi',\hbar\xi_n)\ u( y_n)\,dy_n\dbar \xi_n\nonumber\\
&=& \int \int_0^\infty e^{i(x_n/\hbar-y_n)\xi_n} 
p(x',x_n,\hbar\xi',\xi_n)\ 
u(\hbar y_n)\,dy_n \dbar \xi_n\nonumber\\
&=&\kappa_{\hbar^{-1}}\int \int_0^\infty 
e^{i(x_n-y_n)\xi_n} p(x',\hbar x_n,\hbar\xi',\xi_n)
\kappa_{\hbar}u(y_n)\,dy_n\dbar\xi_n\nonumber.
%&&- \hbar\ \int e^{i(x_n/\hbar-y_n)\xi_n} 
%D_{\xi_n}p^1(x,\hbar\xi',\xi_n)\ 
%u(\hbar y_n)\,dy_n d\xi_n\nonumber\\
%&=&\left(\kappa_{\hbar^{-1}}\left(\op_{x_n}p^0\right)_\hbar
%\kappa_{\hbar}\right)u(x_n)+ \hbar\left(\kappa_{\hbar^{-1}}\left(\op_{x_n}r(h)\right)_\hbar
%\kappa_{\hbar}\right)u(x_n),
\end{eqnarray}
In case $p$ is independent of $x_n$, the last term equals 
$\left(\op_{x_n}p\right)_\hbar u(x_n)$
where the subscript indicates that we use the $\hbar$-scaled 
symbol associated to the operator-valued symbol $\op_{x_n}p$.
}

%%%%%%%%%%%%%%%%%%%%%%%%%%%%%%%%%%%%%%%%%%%%%%%%%%%%%%%%%%%%%%%%%%%%%%%

%%%%%%%%%%%%%%%%%%%%%%%%%%%%%%%%%%%%%%%%%%%%%%%%%%%%%%%%%%%%%%%%%%%%%%%
%ES notation simplified: 'symbols' instead of 'boundary symbol operators'
\extra{h.2}{Potential, trace, 
and singular Green boundary symbol operators}{%
We define the $\hbar$-scaled operators as in \ref{1.15},
noting that the group action on $\C$ is the identity.
Specifically, 
$$\begin{array}{rcll}
k_\hbar(x',\xi') &=& 
\gk_\hbar^{-1}k(x',\hbar\xi')&\text{(potential symbols)}\\ 
% boundary symbol operator
t_\hbar(x',\xi') &=& 
t(x',\hbar\xi')\gk_\hbar&\text{(trace symbols)}\\
g_\hbar(x',\xi') &=& \gk_\hbar^{-1} g(x',\hbar\xi')\gk_\hbar
&\text{(singular Green symbols).}
\end{array}
$$
}
%%%%%%%%%%%%%%%%%%%%%%%%%%%%%%%%%%%%%%%%%%%%%%%%%%%%%%%%%%%%%%%%%%%%%%%
%%%%%%%%%%%%%%%%%%%%%%%%%%%%%%%%%%%%%%%%%%%%%%%%%%%%%%%%%%%%%%%%%%%%%%%
\lemma{1.15b}{Let $g\in 
S^\mu(\R^{n-1}\times\R^{n-1},\cS'(\R_+),\cS(\R_+))$, $\mu\in\R$,
and $\gvp\in \cC^\infty_b(\R^n)$, supported in $\R^n_+$. 
%\cap\supp\psi=\emptyset$.
Then, for every $N\in\N$, we have 
$$\gvp\op(g_\hbar)=\hbar^N\op r_N(\hbar)$$
with a family $r_N$ of singular Green symbols 
such that 
$r_{N,1/\hbar}$ is bounded in $S^{\mu-N}$.  
}

%%%%%%%%%%%%%%%%%%%%%%%%%%%%%%%%%%%%%%%%%%%%%%%%%%%%%%%%%%%%%%%%%%%%%%%
\Proof. 
Write $\gvp_N(x) = x_n^{-N}\gvp\in \cC_b^\infty(\R^n)$.
In view of the fact that 
$M_{x_n^N}\gk_\hbar^{-1}= \hbar^N\gk_\hbar^{-1}M_{x_n^N}$ we have  
$$M_\gvp\op g_\hbar 
=\hbar^N M_{\gvp_N}
\op(\gk_{\hbar}^{-1}(x_n^Ng(x',\hbar\xi')\gk_\hbar).$$
This shows the assertion, since ${x_n^N}g$ is a singular Green symbol 
of order $\mu-N$.\eproof  

%%%%%%%%%%%%%%%%%%%%%%%%%%%%%%%%%%%%%%%%%%%%%%%%%%%%%%%%%%%%%%%%%%%%%%%

\lemma{h.4}{Let $g_1$ and $g_2$ be singular Green symbols
of orders $\mu_1$ and $\mu_2$, respectively.
%For $j=1,2$ let
%$g_j\in S^{\mu_j}(\R^{n-1}\times\R^{n-1};\cS'(\R_+),\cS(\R_+))$.
Then 
\begin{eqnarray}
\nonumber
g_{1,\hbar}g_{2,\hbar}-\left(g_1g_2\right)_\hbar=
\hbar\ c(\hbar) 
\end{eqnarray}
for a family $c(\hbar)$, $0<\hbar\le 1$, of 
singular Green symbols 
with $c_{1/\hbar}$ bounded of order $\mu_1+\mu_2-1$. 
%with $c_{1/\hbar}$ bounded in 
%$S^{\mu_1+\mu_2-1} (\R^{n-1}\times\R^{n-1};\cS'(\R_+),\cS(\R_+))$.
}

%%%%%%%%%%%%%%%%%%%%%%%%%%%%%%%%%%%%%%%%%%%%%%%%%%%%%%%%%%%%%%%%%%%%%%%

\Proof. This is immediate from  \ref{1.15}(c) in connection with
Definition \ref{h.2}.\eproof

%%%%%%%%%%%%%%%%%%%%%%%%%%%%%%%%%%%%%%%%%%%%%%%%%%%%%%%%%%%%%%%%%%%%%%%
\rem{h.2a}{In the following Proposition we will show an analog of 
\ref{1.15}(c). 
The statement becomes more involved, since we have to take care of the 
leftover term. We thus fix the notation beforehand.

Let
$p_j\in S^{\mu_j}_{\tr}(\R^n\times\R^n)$, $\mu_j\in \Z,\ j=1,2$.
For fixed $(x',\xi')$, we consider the operators 
$(\op^+_{x_n} p_j)(x',\xi')$.  
Their composition 
gives rise to a $(x',\xi')$-dependent leftover term. 
We denote by 
$l(p_1,p_2)$ its operator-valued singular Green symbol.
In order to study it we follow Grubb \cite[(2.6.18)ff]{Grubb} %\TEXT{CHK} 
and decompose as in \eqref{B.4.1}:
\begin{equation}\label{h.2a.0}
p_j=s_j+q_j %\quad S_j=\sum_{k=0}^{\mu_j}a^{(j)}_k(x,\xi')D_n^k
\end{equation}
with differential symbols $s_j$ and symbols $q_j$ of normal order
$-1$. %zero (this is only necessary if $\mu_j>0$).
We write $s_2(x,\xi)=\sum_{j=0}^{\mu_2}a_j(x,\xi')\xi_n^j$
and $\gg_m=\gg_0\circ D_n^m$ 
Then 
\begin{eqnarray}\label{h.2a.2}
l(p_1,p_2)
=\sum_{m=0}^{\mu_2-1}k_m\gg_m + 
g^+(q_1)g^-(q_2).
\end{eqnarray}
with the potential symbol $k_m$  given by 
\begin{eqnarray}\label{h.2a.3}
k_m(x',\xi')v 
= i r^+\left(p_1(x,\xi',D_n)
\sum_{l=m+1}^{\mu_2} a^{}_l(x,\xi')D^{l-1-m}_n\right)(v\otimes \gd);
\end{eqnarray}
and the singular Green symbols $g^\pm$ defined by
\begin{eqnarray}\label{gminus}
g^+(q_1) = r^+\op_{x_n}(q_1)e^-J\quad\text{ and }\quad
g^-(q_2) = Jr^-\op_{x_n}(q_2)e^+
\end{eqnarray}
with the reflection operator $J: f(x',x_n)\mapsto f(x',-x_n)$. 

We will need a semiclassical version of the
following  well-known statement: 

\thm{h.2z}{Given a pseudodifferential symbol $p$ of order
$\mu$ with the transmission property and
$l\in \N_0$, the prescription
$$v\mapsto r^+\op p(v\otimes D_n^l\delta),\quad v\in\cS(\R^{n-1}),$$
defines a potential symbol $k$ of order $\mu+l+1/2$ whose 
symbol seminorms can be estimated in terms of those of $p$.
Writing  $p=p^0+x_np^1$ with $p^0(x',\xi)=p(x',0,\xi)$,
we have $k=k^0+k^1$, where
\begin{eqnarray}\label{h.2a.1}
k^0(x',\xi')= \op_{x_n}(p^0(x',\xi)\xi_n^l)(\delta_{x_n=0}),
\end{eqnarray}
considered as a multiplication operator on $\cS(\R_+)$, and $k^1$ 
is of order $\mu+l-1/2$. }

For a proof in the spirit of operator-valued symbols see \cite[Lemma 2.11]{SIB}. 
In the semiclassical situation we obtain:
 %, cf.\ \cite[]{Grubb}.
} 
%%%%%%%%%%%%%%%%%%%%%%%%%%%%%%%%%%%%%%%%%%%%%%%%%%%%%%%%%%%%%%%%%%%%%%%

\lemma{h.2b}{Under the above assumptions, 
$$v\mapsto r^+\op_{x_n} p_\hbar (v\otimes \hbar^{l+1/2}D_n^l\delta) ,
\quad v\in \cS(\R^{n-1})$$
defines a  potential boundary symbol operator $k(\hbar)$. We have 
$$k(\hbar) = k^0_\hbar + \hbar  k^1(\hbar)$$
with $k^0$ defined in \eqref{h.2a.1}  and  $k^1_{1/\hbar}$ 
bounded of order $\mu+l-1/2$.}

%%%%%%%%%%%%%%%%%%%%%%%%%%%%%%%%%%%%%%%%%%%%%%%%%%%%%%%%%%%%%%%%%%%%%%%

\Proof. Replacing $p$ by $p\,\xi_n^l$ we can assume $l=0$.
This yields the potential symbol 
\begin{eqnarray}
 \lefteqn{k(\hbar;x',\xi')
=\hbar^{1/2}r^+\int e^{ix_n \xi_n}p(x,\hbar\xi)\, \dbar\xi_n
}\nonumber\\
&=&\hbar^{-1/2} r^+
\int e^{ix_n\xi_n/\hbar}
\left(p^0(x',\hbar\xi',\xi_n)
+x_n p^1(x,\hbar\xi',\xi_n)\right)\, \dbar\xi_n
\nonumber\\
&=& \gk_\hbar^{-1}r^+
\left((\op_{x_n}p^0)(x',\hbar\xi') \delta +
\int e^{ix_n\xi_n} 
\hbar x_np^1(x,\hbar x_n,\hbar\xi',\xi_n)\, \dbar\xi_n
\right)\nonumber\\
&=&k_\hbar^0(x',\xi') +\hbar k^1(\hbar;x',\xi').\nonumber
\end{eqnarray}
Rescaling the potential symbols $k^1(\hbar) $ we obtain
$$k^1_{1/\hbar}(\hbar;x',\xi')= \int e^{ix_n\xi_n} 
 x_np^1(x,\hbar x_n,\xi',\xi_n)\, \dbar\xi_n
= \int e^{ix_n\xi_n} 
 (-D_{\xi_n})p^1(x,\hbar x_n,\xi',\xi_n)\, \dbar\xi_n,$$
 which is uniformly bounded in 
$S^{\mu-1/2}(\R^{n-1}\times\R^{n-1}; \C,\cS(\R_+))$ by
Theorem \ref{h.2z}. % 2.11 in \cite{SIB}.
\eproof

%%%%%%%%%%%%%%%%%%%%%%%%%%%%%%%%%%%%%%%%%%%%%%%%%%%%%%%%%%%%%%%%%%%%%%%

In the following proposition, we shall analyze the relation between
the interior symbols and the leftover term. 

%%%%%%%%%%%%%%%%%%%%%%%%%%%%%%%%%%%%%%%%%%%%%%%%%%%%%%%%%%%%%%%%%%%%%%%

\prop{h.3}{We use the notation of Remark \ref{h.2a}. 
 Then 
\begin{eqnarray}
\nonumber
\op^+_{x_n}p_{1,\hbar}\circ \op^+_{x_n}p_{2,\hbar}-
\op^+_{x_n}(p_{1}p_{2})_\hbar -
l(p_1,p_2)_\hbar
=\hbar\ \left(\op^+_{x_n} c(\hbar) +  d(\hbar)\right) 
\end{eqnarray}
for two  families $c$ and $d$ 
with $c_{1/\hbar}$ bounded in $S^{\mu_1+\mu_2-1}_{\tr}(\R^n\times\R^n)$
and  
$d_{1/\hbar}$ bounded in 
$S^{\mu_1+\mu_2-1} (\R^{n-1}\times\R^{n-1};
\cS'(\R_+), \cS(\R_+))$.
}
%%%%%%%%%%%%%%%%%%%%%%%%%%%%%%%%%%%%%%%%%%%%%%%%%%%%%%%%%%%%%%%%%%%%%%%

\Proof. 
% We have $(P_{1,\hbar}^+)(P_{2,\hbar}^+)=
%(P_{1,\hbar}P_{2,\hbar})^++L(P_{1,\hbar},P_{2,\hbar})$.
We already know from \ref{1.15} (for the standard $\C$-valued case) that
$$\op_{x_n}p_{1,\hbar}\op_{x_n}p_{2,\hbar}
-\op_{x_n}\left(p_{1}p_{2}\right)_\hbar
=\hbar\ \op(c(\hbar)),$$
with  $c_{1/\hbar}$ bounded in $S^{\mu_1+\mu_2-1}$.
In fact, we even obtain boundedness in the corresponding class 
with the transmission property, since we have
an asymptotic expansion for $c$ with terms bounded in the topology 
of symbols with the transmission property.

Let us next consider the  symbol
$l(p_{1,\hbar},p_{2,\hbar})$ of the leftover term 
$L(\op_{x_n}p_{1,\hbar},\op_{x_n}p_{2,\hbar})$, which we 
will compute according to \eqref{h.2a.2}, \eqref{h.2a.3}.
Replacing $\xi$ by $\hbar\xi$, we obtain
$$l(p_{1,\hbar},p_{2,\hbar}) 
=\sum_{m=0}^{\mu_2-1}k_m(\hbar)\gg_m + 
g^+(q_{1,\hbar})g^-(q_{2,\hbar}).$$
Using the notation introduced above, 
the potential symbol $k_m$ is given by 
$$k_m(\hbar;x',\xi')v 
= i r^+\left(\op_{x_n}p_1(x,\hbar\xi)
\sum_{l=m+1}^{\mu_2} 
a^{}_l(x,\hbar\xi')\hbar^l D_n^{l-1-m}\right)(v\otimes \gd).$$
%We know that 
%\begin{eqnarray*}
%\lefteqn{\left(\op_{x_n}p_1(x,\hbar\xi)\right)\ 
%a_l(x,\hbar\xi')\hbar^l D_n^{l-m-1}}\\
%&=&\hbar^{m+1/2}\left(\op_{x_n}\left(p_1(x,\hbar\xi)a_l(x,\hbar\xi')
%\hbar^{l-m-1/2} \xi_n^{l-m-1}\right)
%+\hbar \op_{x_n}r(\hbar)\right)
%\end{eqnarray*}
%for a family $r$ with $r_{1/\hbar}$ bounded in $S^{\mu_1+\mu_2-m-1}_{\tr}$.
By Lemma \ref{h.2b},  
 $\hbar^{-m-1/2}k_m(\hbar)$ is a family of potential symbols of order 
$\mu_1+\mu_2-m-1/2$, 
equal to $k_{m,\hbar}^0$ modulo lower order 
terms. In view of the fact that $\gg_{m,\hbar}=\hbar^{m+1/2}\gg_m$, 
this shows that the composition $k_m(\hbar)\gg_m$  
is a singular Green symbol which equals 
$k_{m,\hbar}^0\gg_{m,\hbar}$ modulo lower order terms of the desired form.
 
The singular Green symbols  $g^\pm(p)=g^\pm(p)(x',\xi')$ 
associated to a pseudodifferential
symbol $p$ of negative normal order are the integral operators
with the  kernels 
%$\tilde g^\pm(p)(x',\xi',x_n,y_n)$ given by
\begin{eqnarray}\label{gplus}
\tilde g^\pm(p)(x',\xi',x_n,y_n)
&=&\left.\int e^{iz\xi_n}p(x',x_n,\xi',\xi_n)\dbar\xi_n
\right|_{z=\pm(x_n+y_n)}
\end{eqnarray}
so that
\begin{eqnarray*}
\tilde g^\pm(p_\hbar)(x',\xi',x_n,y_n)
&=&\hbar^{-1}\, \left.\int 
e^{iz\xi_n}p(x',x_n, \hbar\xi',\xi_n)\dbar\xi_n
\right|_{z=\pm(x_n+y_n)/\hbar}
\end{eqnarray*}
This implies that 
$$g^\pm(p_\hbar)(x',\xi') 
=\gk_\hbar^{-1} g^\pm(p(x',\hbar x_n,\hbar\xi',\xi_n))\gk_\hbar$$
Writing $p=p^0+x_np^1$ as before, we have
%with $p^0(x',\xi)=p(x',0,\xi)$, we have
$$g^\pm(p_\hbar)(x',\xi') 
=\left(g^\pm(p^0)\right)_\hbar+  
\hbar \gk^{-1}_\hbar x_n 
g^\pm\left(p^1(x',\hbar x_n,\hbar\xi',\xi_n)\right) \gk_\hbar$$ 
As %In view of the fact that 
$p^1(x',\hbar x_n,\hbar\xi',\xi_n)$, 
$0<\hbar\le 1,$ is a  bounded 
family of pseudodifferential operators of order $\mu$ 
with the transmission property and
since multiplication by $x_n$ lowers the order by $1$ 
according to \ref{1.10b}(a), 
we obtain the assertion of the proposition by 
first applying this consideration 
to $g^+(q_{1,\hbar})$ and $g^-(q_{2,\hbar})$ 
and then using Lemma \ref{h.4}.
\eproof

%%%%%%%%%%%%%%%%%%%%%%%%%%%%%%%%%%%%%%%%%%%%%%%%%%%%%%%%%%%%%%%%%%%%%%%

\section{The Graph Projection}

%%%%%%%%%%%%%%%%%%%%%%%%%%%%%%%%%%%%%%%%%%%%%%%%%%%%%%%%%%%%%%%%%%%%%%%

\dfn{2.0}{%ES cosmetic
The graph projection of a bounded 
operator $a$ on a Hilbert space $H$ is the operator $\Gr(a)$ on $H\oplus H$
given by  
$$\Gr(a)= \Mat{(1+a^*a)^{-1}& (1+a^*a)^{-1}a^*\\ 
a(1+a^*a)^{-1}& a(1+a^*a)^{-1}a^*}
%-\Mat{0&0\\0&1}
=
\Mat{(1+a^*a)^{-1}& (1+a^*a)^{-1}a^*\\ 
a(1+a^*a)^{-1}& 1-(1+aa^*)^{-1}}.$$
For unbounded operators, it is not a priori clear that the above 
definition makes sense, nor is  the second 
identity obvious. We will now have a closer look.}

%%%%%%%%%%%%%%%%%%%%%%%%%%%%%%%%%%%%%TTT%%%%%%%%%%%%%%%%%%%%%%%%%%%%%%%%%%

\extra{fa.0}{The framework}{%
Let $V\hookrightarrow  H\hookrightarrow V'$ be Hilbert spaces, 
with $V$ dense in $H$, and assume that $V'$ is the dual space of $V$ 
with respect to an extension 
$\skp{\cdot,\cdot}$ of the inner product in $H$. 
Moreover let $a_0:V\to H$ be a bounded operator with adjoint 
$a_0^*:H\to V'$.
We assume that $a_0$ is closable in $H$ and denote by $a$ the closure. 
Explicitly: 
The domain $\cD$ of $a$ consists of all $x$ in $H$ 
for which there exists a sequence $(x_n)$ in $V$ with $x_n\to x$ in $H$ 
and $a_0x_n\to y$ in $H$ for some $y\in H$. 
In that case, we define $ax=y$. }   

\lemma{fa.1}{$\cD$ naturally is a Hilbert space with the inner product
\begin{eqnarray}\label{fa.1.1}
 \skp{\!\skp{x,y}\!}_{\cD}=\skp{x,y}+\skp{ax,ay}.
\end{eqnarray}
$V$ is dense in $\cD$.}

\Proof. 
Clearly, \eqref{fa.1.1} defines an inner product on $\cD$. 
The associated norm $\|x\|_{\cD}=(\|x\|^2+\|ax\|^2)^{1/2}$ 
is the graph norm with respect to which $\cD$ is complete. 
Hence $\cD$ is a Hilbert space. $V$ is dense in $\cD$ by construction. 
\eproof

\lemma{fa.2}{Let $x\in H$. 
Then the element $a_0^*x\in V'$ extends to a continuous 
linear functional on $\cD$, and
\begin{eqnarray}\label{fa.2.1}
\skp{a_0^*x,y}= \skp{x,ay},\quad x\in H, y\in \cD. 
\end{eqnarray}
}

\Proof.  For $v$ in the dense subspace $V$ we have
$$|\skp{a_0^*x,v}|=|\skp{x,av}|\le \|x\|\|av\|
\le \|x\| \|v\|_{\cD}$$ 
so that $a_0^*x$ extends continuously to $\cD$.
The stated identity follows.\eproof 

We now denote by $\cE$ the range of the operator $1+a_0^*a:\cD\to V'.$

%%%%%%%%%%%%%%%%%%%%%%%%%%%%%%%%%%%%%%%%%%%%%%%%%%%%%%%%%%%%%%%%%%%%%%%

\lemma{fa.3}{The elements of $\cE$ define continuous 
linear functionals on $\cD$. The norm of $(1+a_0^*a)x$ on $\cD$ is
$\|x\|_{\cD}$, and the operator $1+a_0^*a$ is injective.}

\Proof. As $\cD\hookrightarrow H$, the first statement follows from
Lemma  \ref{fa.2}. Note that 
\begin{eqnarray}
\|(1+a_0^*a)x\|_{\cL(\cD,\C)} &=&
\sup_{\|y\|_{\cD}\le 1}|\skp{(1+a_0^*a)x,y}|=
\sup_{\|y\|_{\cD}\le1}|\skp{x,y}+\skp{ax,ay}| \nonumber
\\
&=&\sup_{\|y\|_{\cD}\le1}|\skp{\!\skp{x,y}\!}_{\cD}|
= \|x\|_{\cD}\nonumber.
\end{eqnarray} 
%On the other hand $\|x\|_{\cD}$ is attained by choosing $y=x/\|x\|_{\cD}$.
This implies injectivity. \eproof 

%%%%%%%%%%%%%%%%%%%%%%%%%%%%%%%%%%%%%%%%%%%%%%%%%%%%%%%%%%%%%%%%%%%%%%%

\lemma{fa.4}{$\cE$ inherits a Hilbert space structure from $\cD$. 
The associated norm is the norm in $\cD'=\cL(\cD,\C)$. 
This allows us to identify $\cE$ with the dual of $\cD$ with respect 
to the sesquilinear  pairing between $V$ and $V'$.
$\cE$ contains the range of $a_0^*$}  
\Proof. Let $e=x+a_0^*ax$ and $f=y+a^*_0ay$ be two elements of $\cE$. 
Then we let 
\begin{eqnarray}
\skp{\!\skp{e,f}\!}_{\cE}= \skp{\!\skp{x,y}\!}_{\cD},
\end{eqnarray}  
so that the Hilbert space structure of $\cD$ carries over to $\cE$.
The associated norm of $e=(1+a_0^*a)x$ then is equal to 
$\|x\|_{\cD}$, which is the norm of $(1+a_0^*a)x$ in $\cD'$ 
by Lemma \ref{fa.3}.
Moreover, the identity 
$$\skp{x+a_0^*ax,y}= \skp{\!\skp{x,y}\!}_{\cD}$$
shows that the action of $\cE$ on $\cD$ via the sesquilinear pairing 
coincides with the pairing with elements of 
$\cD$, which gives the whole dual space. By Lemma \ref{fa.2},
$a_0^*(H)\subset \cD'=\cE.$  \eproof

\cor{fa.5}{The operator $1+a_0^*a:\cD\to \cE$ is invertible with 
a bounded inverse. 
Moreover,  the operators $(1+a_0^*a)^{-1}$,  
$(1+a_0^*a)^{-1}a_0^*$, $a(1+a_0^*a)^{-1}$, and $a(1+a_0^*a)^{-1}a_0^*$  
define elements of $\cL(H)$.}

%%%%%%%%%%%%%%%%%%%%%%%%%%%%%%%%%%%%%%%%%%%%%%%%%%%%%%%%%%%%%%%%%%%%%%%

\lemma{fa.7}{The restriction of $(1+a_0^*a)^{-1}$ to $H$ is 
self-adjoint, and the operators $a(1+a_0^*a)^{-1}$ and 
$(1+a_0^*a)^{-1}a^*_0$ are adjoints of each other in $\cL(H)$.
}

%%%%%%%%%%%%%%%%%%%%%%%%%%%%%%%%%%%%%%%%%%%%%%%%%%%%%%%%%%%%%%%%%%%%%%%

\Proof. Let $x\in H\subset \cE$ and $z=(1+a_0^*a)^{-1}x\in\cD$.
For $y\in H$,  Lemma \ref{fa.2} implies that
\begin{eqnarray*}
\skp{(1+a_0^*a)^{-1}y,x}= \skp{(1+a_0^*a)^{-1}y,(1+a_0^*a)z}
=\skp{y,z}=\skp{y,(1+a_0^*a)^{-1}x}.
\end{eqnarray*} 
This shows that $(1+a_0^*a)^{-1}$ is selfadjoint. For 
$x,y\in H$ the element $a_0^*x\in\cE$ has a preimage $z$ in $\cD$
under $1+a^*_0a$.  We infer from \eqref{fa.2.1} that
\begin{eqnarray*}
\skp{a(1+a_0^*a)^{-1}y,x}&=&\skp{(1+a_0^*a)^{-1}y,a_0^*x} =
\skp{(1+a_0^*a)^{-1}y,(1+a_0^*a)z}\\
&=& \skp{y,z}
=\skp{y,(1+a^*_0a)^{-1}a^*_0x}.
\end{eqnarray*}
This shows the second statement.
\eproof
%%%%%%%%%%%%%%%%%%%%%%%%%%%%%%%%%%%%%%%%%%%%%%%%%%%%%%%%%%%%%%%%%%%%%%%

\extra{fa.6}{Notation}{In the above we wrote $a_0^*$ 
in order to stress the fact that this 
operator is {\em not} the $H$-valued Hilbert space adjoint of $a$
but an operator with values in $V'$. 
Now that this has been made
clear we shall go back to the simpler notation and write $a^*$. }

%%%%%%%%%%%%%%%%%%%%%%%%%%%%%%%%%%%%%%%%%%%%%%%%%%%%%%%%%%%%%%%%%%%%%%%

%ES bundles inserted
\extra{setup}{The set-up}{In the sequel, 
$$T:\cC^\infty(X, E_1)\to \cC^\infty(X,E_2)$$ 
will be a Fredholm operator of order and class zero in 
Boutet de Monvel's calculus.
Following an idea of Elliott-Natsume-Nest \cite{ENN} we will 
associate to $T$ the operator $A=\gL^{m,+}_-T$ for some  $m>n$ 
and then study the graph projection. 

The operator $T$ induces 
%The restriction $T_m$ of $T$ to $H^m(X,E_1)$ induces 
a Fredholm operator  $T_m:H^m(X,E_1)\to H^m(X,E_2)$. 
Its adjoint $T^*_m:H^{-m}_0(X,E_2)\to H^{-m}_0(X,E_1)$ 
extends the $L^2$-adjoint $T^*$. 
We let
$$A=\Lambda^{m,+}_-T_m:H^m(X,E_1)\to L^2(X,E_2).$$
This is an operator of order $m$ and class zero;
there is no leftover term in the composition.
The adjoint $A^*$  is the operator
$(\Lambda^{m,+}_-T_m)^*=T_m^*\Lambda^m_+:L^2(X,E_2)\to H^{-m}_0(X,E_1).$

We consider the composition $A^*A$ as the bounded operator 
$$T^*_m\gL^m_+\gL^{m,+}_-T_m:H^m(X,E_1)\to H^{-m}_0(X,E_1).$$

On the other hand,  $A=\gL^{m,+}_-T$ 
extends to a bounded operator from
$ L^2(X,E_1)$ to $H^{-m}(X,E_2)$ with adjoint 
$A^*=T^*\gL^{m}_+:H^m_0(X,E_2)\to L^2(X,E_1)$, and 
$AA^*=\gL^{m,+}_-TT^*\gL^{m}_+$ maps $H_0^m(X,E_2)$ to $H^{-m}(X,E_2)$.
}
 
%%%%%%%%%%%%%%%%%%%%%%%%%%%%%%%%%%%%%%%%%%%%%%%%%%%%%%%%%%%%%%%%%%%%%%%
%ES  changed for bundles
\lemma{2.1}{We have natural embeddings 
\begin{eqnarray*}
&&H^m(X,E_j)\hookrightarrow L^2(X,E_j)\hookrightarrow H^{-m}_0(X,E_j)
\quad\text{ and}\\
&&H^m_0(X,E_j)\hookrightarrow L^2(X,E_j)\hookrightarrow H^{-m}(X,E_j),
\end{eqnarray*} 
$j=1,2,$ and topological isomorphisms
$$1 + A^*A:H^m(X,E_1) \to H^{-m}_0(X,E_1)\quad 
\text{and}\quad 1+AA^*:H^m_0(X,E_2)\to H^{-m}(X,E_2).$$
}

%%%%%%%%%%%%%%%%%%%%%%%%%%%%%%%%%%%%%%%%%%%%%%%%%%%%%%%%%%%%%%%%%%%%%%%

\Proof. The embeddings are well-known. The second statement  follows 
from Corollary \ref{fa.5}, applied to the operator
$$ a= \Mat{0&A^*\\A&0}: \begin{array}{ccc}
H^m(X,E_1)&&L^2(X,E_1)\\
\oplus&\to&\oplus\\ 
H^m_0(X,E_2)&& L^2(X,E_2)
\end{array},$$
together with the fact that -- due to elliptic regularity --
the domain of the closure of $A$ is $H^m(X,E_1)$, that of $A^*$ is 
$H^m_0(X,E_2)$. 
\eproof
\medskip

%%%%%%%%%%%%%%%%%%%%%%%%%%%%%%%%%%%%%%%%%%%%%%%%%%%%%%%%%%%%%%%%%%%%%%%

As a consequence, the operators $A(1+A^*A)^{-1}A^*$ and 
$(1+AA^*)^{-1}$ are bounded operators on $L^2(X,E_1)$ 
and $L^2(X,E_2)$, respectively. 

%%%%%%%%%%%%%%%%%%%%%%%%%%%%%%%%%%%%%%%%%%%%%%%%%%%%%%%%%%%%%%%%%%%%%%%

\lemma{2.2}{The restriction of $AA^*$ to $H^{2m}_0(X,E_2)$ 
maps to $L^2(X,E_2)$. }

%%%%%%%%%%%%%%%%%%%%%%%%%%%%%%%%%%%%%%%%%%%%%%%%%%%%%%%%%%%%%%%%%%%%%%%

\Proof. 
We know that $\gL^m_+:H^{2m}_0\stackrel
\cong\longrightarrow H^m_0$.
Now we observe that 
$H^m_0$ naturally embeds into $H^m$, and
thus $TT^*$ defines a bounded map from $H^m_0$ %\hookrightarrow H^m
to $H^m$. 
Finally, $\gL^m_-:H^m\stackrel\cong\longrightarrow L^2$.
\eproof

%%%%%%%%%%%%%%%%%%%%%%%%%%%%%%%%%%%%%%%%%%%%%%%%%%%%%%%%%%%%%%%%%%%%%%%

\lemma{2.3}{$A(1+A^*A)^{-1}A^*=1-(1+AA^*)^{-1}$ on $L^2(X,E_2)$.}

%%%%%%%%%%%%%%%%%%%%%%%%%%%%%%%%%%%%%%%%%%%%%%%%%%%%%%%%%%%%%%%%%%%%%%%

\Proof. 
According to Lemma \ref{2.2}, the composition 
$A(1+A^*A)^{-1}A^*(1+AA^*)$ is defined as a 
bounded operator from $H^{2m}_0$ to $L^2$, and we have 
\begin{eqnarray}\label{2.3.1}
A(1+A^*A)^{-1}A^*(1+AA^*)=A(1+A^*A)^{-1}(1+A^*A)A^*=AA^*.
\end{eqnarray}
We next denote by $\cR$ the range of the restriction of 
$(1+AA^*)^{-1}$ to $L^2$, so that $1+AA^*:\cR\to L^2$ 
is an isomorphism.
According to \ref{2.1} and \ref{2.2}, 
we have $H^{2m}_0\subseteq \cR \subseteq H^m_0$.

As \eqref{2.3.1} extends to $\cR$, 
the compositions, below, are defined on $L^2$, 
and
\begin{eqnarray*}
\lefteqn{A(1+A^*A)^{-1}A^*}\\
&=& A(1+A^*A)^{-1}A^*(1+AA^*)(1+AA^*)^{-1}
=AA^*(1+AA^*)^{-1}=1-(1+AA^*)^{-1}.
\end{eqnarray*}
\eproof
%%%%%%%%%%%%%%%%%%%%%%%%%%%%%%%%%%%%%%%%%%%%%%%%%%%%%%%%%%%%%%%%%%%%%%%

\lemma{2.4}{$\Gr(A)$ and $e=\Mat{0&0\\0&1}$
define idempotents in %the $C^*$-algebra 
$\cK(L^2(X,E_1\oplus E_2))^\sim$, where as usual
the tilde indicates the unitization.

In particular, the difference  $[\Gr(A)]-[e]$
defines a class in 
$K_0(\cK(L^2(X, E_1\oplus E_2))).$}

%%%%%%%%%%%%%%%%%%%%%%%%%%%%%%%%%%%%%%%%%%%%%%%%%%%%%%%%%%%%%%%%%%%%%%%

\Proof. Apart from the one in the lower right corner,
the entries in $\Gr(A)$ are 
%The operators $(1+A^*A)^{-1}$, $(1+A^*A)^{-1}A^*$, 
%$A(1+A^*A)^{-1}$ are 
compact on $L^2$ as a consequence of the 
compact embeddings  $H^m_0\hookrightarrow L^2$ and 
$H^{m}\hookrightarrow L^2$. By 
Lemma \ref{2.3}, the last entry differs from the identity by 
%$A(1+A^*A)^{-1}A^*=1-(1+AA^*)^{-1}$ with 
the compact operator $ (1+AA^*)^{-1}$. 
\eproof

%%%%%%%%%%%%%%%%%%%%%%%%%%%%%%%%%%%%%%%%%%%%%%%%%%%%%%%%%%%%%%%%%%%%%%%

\thm{2.5}{The class $[\Gr(A)]-[e]$
%\eqref{K} 
in  $K_0(\cK(L^2(X,E_1\oplus E_2)))$ equals 
$$
\left[\pi_{\ker A}\right]-\left[\pi_{\ker A^*}\right]=
\left[\pi_{\ker T}\right]-\left[\pi_{\ker T^*}\right].
$$
Here, $\pi_V$ denotes the orthogonal projection onto $V$
with respect to the $L^2$-inner product.}
%%%%%%%%%%%%%%%%%%%%%%%%%%%%%%%%%%%%%%%%%%%%%%%%%%%%%%%%%%%%%%%%%%%%%%%

\Proof. Replacing $A$ by $tA$ for $t\ge 1$, we consider 
$\Gr(tA)$, which is a norm continuous family of idempotents. 
We claim that $\Gr(tA)$ converges to 
$$\Mat{\pi_{\ker T}&0\\0&1-\pi_{\ker T^*}} $$
as $t\to\infty$.
We let $\cH_1=(1+A^*A)^{-1}(L^2)$ %=\{x\in H^m: A^*Ax\in L^2\}$
and $\cH_2=.(1+AA^*)^{-1}(L^2)$. % =\{x\in H_0^m: AA^*x\in L^2\}$
This allows us to consider $A^*A$ and $AA^*$ as an unbounded 
operators on $L^2$ with domains $\cH_1$ and $H_2$, respectively. 
The graph projection is not affected by this 
change. For the unbounded operator, however, the statements are 
well-known. They are a consequence of the fact that $0$ is an 
isolated point in the spectrum.
Alternatively, the statement can be checked by a direct computation. 
\forget{
A first observation is that $A^*Ax$ 
is an element in $ \im T^*$ for $x\in \cH$:
In fact, since $\cH\subset H^{m}$, 
$A^*Ax=T_m^*y$ for $y=\gL_+^m\gL_-^{m,+}T_mx\in H^{-m}_0$.
As $T^*_my\in L^2$, elliptic regularity implies that $y\in L^2$, so 
that $T^*_my=T^*y\in \im T^*$.

Our next observation is that, for $f\in\cH$,
\begin{eqnarray*}
f\perp \ker T\Leftrightarrow (1+t^2A^*A)f\perp\ker T, 
\end{eqnarray*}
where both `$\perp$' refer to the inner product on $L^2$. 
Indeed, $(\ker T)^\perp=\im T^*$, since $T^*$ has closed range.
So if $f\in \im T^*$, then $(1+t^2A^*A)f=f+t^2A^*Af\in\im T^*$,
and vice versa.

On $\ker T \cap\cH$, we have $1+t^2A^*A=1$. 
For $f\perp\ker T$, on the other hand,
the fact that $T$ is Fredholm on $L^2$ implies that
there is a constant $c_0$ such that 
$\|Tf\|_{L^2}\ge c_0\|f\|_{L^2}$, and
\begin{eqnarray}\label{lowerest}
\lefteqn{\|(1+t^2A^*A)f\|^2_{L^2}
\ge
\skp{f,f}+2t^2\skp{Af,Af}\nonumber
= \|f\|^2_{L^2}+2t^2\|Tf\|^2_{H^m}}\\
&\ge& 
\|f\|^2_{L^2}+2t^2\|Tf\|^2_{L^2} 
\ge (1+2c_0^2t^2)\ \|f\|^2_{L^2}. 
\end{eqnarray}
We conclude that $(1+t^2A^*A)^{-1}=1$ on $\ker T$, while 
$\|(1+t^2A^*A)^{-1}\|=O(t^{-2})$ on $(\ker T)^\perp.$
Hence
$$(1+t^2A^*A)^{-1}\to \pi_{\ker T}\quad
\text{ as }t\to\infty.$$

With the same consideration, 
the operator $tA(1+t^2A^*A)^{-1}$ vanishes on $\ker T$ and is 
$O(t^{-1})$ on $(\ker T)^\perp$. As $t\to\infty$, the limit
thus is zero.

It follows from Lemma \ref{fa.7} that its adjoint, 
the operator $(1+t^2A^*A)^{-1}tA^*$, also tends to zero. 

Finally, $t^2A(1+t^2A^*A)^{-1}A^*$ coincides with 
$1-(1+t^2AA^*)^{-1}.$ An argument as above shows that 
$(1+t^2AA^*)^{-1}$ converges to the orthogonal projection onto 
${\ker A^*}={\ker (T_m^*\gL^{m}_+)}$ as $t\to\infty$. 
Since 
%$\gL_+^{m}:L^2\to H^{-m}_0$ is an isomorphism, 
the class of
$\pi_{\ker (T_m^*\gL^{m}_+)} $ in $K_0(\cK)$  coincides with that
of $\pi_{\ker T^*}$ we obtain the assertion.
}%forget
\eproof

%%%%%%%%%%%%%%%%%%%%%%%%%%%%%%%%%%%%%%%%%%%%%%%%%%%%%%%%%%%%%%%%%%%%%%%

\cor{2.6}{For $m>n$, 
$$\ind T=\Tr(1+A^*A)^{-1}-\Tr(1+AA^*)^{-1}.$$}

%%%%%%%%%%%%%%%%%%%%%%%%%%%%%%%%%%%%%%%%%%%%%%%%%%%%%%%%%%%%%%%%%%%%%%%

\Proof. The difference  $\cG(A)- e$ 
%\Mat{0&0\\0&1}$ 
is a trace class 
operator in $\cL(L^2)$ since its four entries are trace class operators. 
This in turn follows from 
the fact that the embeddings $H^m\hookrightarrow L^2$ and 
$L^2\hookrightarrow H^{-m}_0$ 
are trace class. 
According to Theorem \ref{2.5} we then have
$$\ind A =\Tr([\pi_{\ker A}] -[\pi_{\coker A }]) = 
\Tr \left(\cG(A)- e\right)=  
%\Mat{0&0\\0&1}\right) =  
\Tr(1+A^*A)^{-1}-\Tr(1+AA^*)^{-1}.$$\eproof

%%%%%%%%%%%%%%%%%%%%%%%%%%%%%%%%%%%%%%%%%%%%%%%%%%%%%%%%%%%%%%%%%%%%%%%

%ES Johannes' change
\forget{\rem{2.4a}{It is well-known that 
$$K_0(\cK(L^2(X, E_1\oplus E_2)))
\cong K_0(\cK(L^2(X))), $$
so that the difference $[\Gr(A)]-[e]$
can be identified with an element in 
$K_0(\cK(L^2(X))).$ 
We can make this precise: 
We choose bundles $F_1$ and $F_2$ such that 
$E_1\oplus F_1$ and $E_2\oplus F_2$ are trivial with fiber $\C^m$. 
We can then view both the graph projection  and 
the projection $\pi_{F_2}$ onto $F_2$ as elements 
%$\cL(L^2(X,E_1\oplus E_2))$ as a subset 
of $\cL(L^2(X,\C^{2m}))$ and identify, with respect 
to the decomposition  
$$\C^{2m}\cong E_1\oplus F_1\oplus E_2\oplus F_2,$$ 
the  class
$[\Gr(A)]-[e]$ with 
%the difference class 
%$$\Mat{(1+A^*A)^{-1}&0&(1+A^*A)^{-1}A^*&0\\0&0&0&0\\
%A(1+A^*A)^{-1}&0&A(1+A^*A)^{-1}A^*&0\\0&0&0&1}-
%\Mat{0&0&0&0\\0&0&0&0\\0&0&1&0\\0&0&0&1},$$
$$[\Gr(A)+\pi_{F_2}]- [\pi_{E_2\oplus F_2}]$$
in $K_0(\cK(L^2(X,\C^{2m})))$ which is canonically 
isomorphic to $K_0(\cK(L^2(X)))$.
}
}%forget

%%%%%%%%%%%%%%%%%%%%%%%%%%%%%%%%%%%%%%%%%%%%%%%%%%%%%%%%%%%%%%%%%%%%%%%
\section{The Graph Projections of the Symbols}

\extra{2.7}{Notation}{%
We denote by 
\begin{eqnarray*}
p^{m}=\gs_\psi^m(A)\quad\text{and}\quad c^{m}=\gs^m_\partial(A)
\end{eqnarray*} 
the homogeneous principal pseudodifferential symbol 
and the homogeneous principal boundary symbol of the operator  $A$
in \ref{setup}.
Both are invertible.

Locally 
\begin{eqnarray}\label{2.7.1}
c^{m}(x',\xi')=p^{m,+}(x',0,\xi',D_n)+g^{m}(x',\xi')
\end{eqnarray}
with a suitable strictly homogeneous singular Green part $g^{m}$.
We then choose a smooth function $p$ on $T^*X$ which coincides with 
$p^{m}$ for $|\xi|\ge1$ and a smooth singular Green symbol $g$
which coincides with $g^{m}$ for $|\xi'|\ge 1$ and let 
$$c(x',\xi')=p^+(x',0,\xi',D_n)+g(x',\xi').$$

We next  apply Corollary \ref{fa.5} to the operator family  
$$1+c^*(x',\xi')c(x',\xi'):
H^m(\R_+,\widetilde E_1)\to H^{-m}_0(\R_+,\widetilde E_1),
\quad (x',\xi')?\in T^*X.$$
For each $(x',\xi')$, we denote by $\cD_{(x',\xi')}$  
the domain of the closure
and by $\cE_{(x',\xi')}$ the range.
}

%%%%%%%%%%%%%%%%%%%%%%%%%%%%%%%%%%%%%%%%%%%%%%%%%%%%%%%%%%%%%%%%%%%%%%%

\prop{2.8}{For each choice of $(x',\xi')$ we have 
\begin{eqnarray}\label{dom}
\cD_{(x',\xi')} =H^m(\R_+,\widetilde E_1))
\quad\text{and}
\quad \cE_{(x',\xi')}=H^{-m}_0(\ol\R_+,\widetilde E_1)).
\end{eqnarray}
The operator family  
$$1+c^*(x',\xi')c(x',\xi'):H^m(\R_+,\widetilde E_1)
\to H^m_0(\ol\R_+,\widetilde E_1)$$
is pointwise invertible.

The inverse $(1+c^*c)^{-1}$ is an element of 
\begin{eqnarray}\label{2.8.0}
S^{-2m}(\R^{n-1}\times\R^{n-1};H^{-m}_0(\R_+,\widetilde E_1)
,H^m(\R_+,\widetilde E_1)).
\end{eqnarray}

\forget{For large $|\xi'|$,   
$(1+c^*c)^{-1}$ coincides with a boundary symbol of order $-2m$
in Boutet de Monvel's calculus  whose 
pseudodifferential part is induced by the symbol 
$(1+p_0^*p_0)^{-1}$, where $p_0(x',\xi)=p(x',0,\xi)$ for the symbol 
$p$ introduced above. } %forget
}
We shall improve this result in \ref{2.8a}, below.

%%%%%%%%%%%%%%%%%%%%%%%%%%%%%%%%%%%%%%%%%%%%%%%%%%%%%%%%%%%%%%%%%%%%%%%

\Proof. By definition, 
$\cD_{(x',\xi')}\subseteq H^m(\R_+,
\widetilde E_1)$ 
is the domain of the closure
of the operator 
$$c(x',\xi'): H^m(\R_+,\widetilde E_1)\to L^2(\R_+,
\widetilde E_2)$$
in $L^2(\R_+,\widetilde E_1)$.
It consists of all $v\in L^2$ for which there is a sequence
$v_k$ in $H^m$ with $v_k\to v$ in $L^2$ and $c(x',\xi')v_k\to 
w$ in $L^2$ for some $w$ which then is defined to be 
$c(x',\xi')v$.  

As $g(x',\xi'):
\cS'(\R_+,\widetilde E_1)\to \cS(\R_+,\widetilde E_2)$ is
continuous, $g(x',\xi')v_k$ will converge for any convergent
$L^2$-sequence $(v_k)$; hence the domain is independent of $g$.
We therefore have
$$\cD_{(x',\xi')}\subseteq 
\{v\in L^2: p^+(x',0,\xi',D_n)v\in L^2\}.$$ 
In view of the fact that $p$ is elliptic, the last set
is a subset of $H^m(\R_+,\widetilde E_1)$, and we get
the first part of \eqref{dom}. 

According to  Lemma \ref{fa.4}, the space $\cE_{(x',\xi')}$
is the dual space of $\cD_{x',\xi')}$ with respect to the pairing
induced by the $L^2$ inner product. This gives the second statement
in \eqref{dom}.
 
We conclude that 
$$(x',\xi')\mapsto 1+c^*(x',\xi')\,c(x',\xi')$$
is a smooth family of operators in 
$$\cL(H^m(\R_+,\widetilde E_1),
H^{-m}_0(\ol\R_+,\widetilde E_1)).$$ 

As inversion is continuous, $(1+c^*(x',\xi')c(x',\xi'))^{-1}$
also is a smooth family. 
Moreover, for $|\xi|\ge 1$, we have $c=c^{m}$, 
and  both $c(x',\xi')$ and $c^*(x',\xi')$ 
are invertible. We can write
\begin{eqnarray}\label{2.8.1}
(1+c^*c)^{-1}=c^{-1}\ (1+(c^*)^{-1}c^{-1})^{-1}\ (c^{*})^{-1}.
\end{eqnarray}
If $\gvp=\gvp(\xi')$ is an excision function on $\R$ which vanishes for
$|\xi'|\le 1$ and is equal to one for large $|\xi'|$, then the
homogeneity of $c$ implies that, in local coordinates, 
\begin{eqnarray}\label{2.8.2}
\gvp\ c^{-1}\in S^{-m}(\R^{n-1}\times\R^{n-1}; L^2(\R_+), H^m(\R_+)).
\end{eqnarray}
Similarly as in \ref{B.21}, 
$c^*(x',\xi')^{-1}=(c^{-1}(x',\xi'))^*.$ Consequently, 
\begin{eqnarray}\label{2.8.3}
\gvp\ c^{*-1}\in S^{-m}(\R^{n-1}\times\R^{n-1}; H^{-m}_0(\R_+), L^2(\R_+)).
\end{eqnarray}
Next we note that the positivity of $m$ implies that 
$$\gk_{\bracket{\xi'}}^{-1}(c^*)^{-1}(x',\xi')\;c^{-1}(x',\xi') 
\gk_{\bracket{\xi'}}\longrightarrow 0 \text{ in } 
\cL(H^{-m}_0(\R_+),H^m(\R_+)).$$
In particular, $(1+(c^*)^{-1}c^{-1})^{-1}$ is uniformly bounded in 
$\cL(L^2(\R_+))$. 
For $d=1+(c^*)^{-1}c^{-1}$ we then deduce 
%from the fact that the group $\gk$ acts unitarily on $L^2$ 
that 
$$\gk_{\bracket{\xi'}}^{-1}\left( \partial_{\xi_j}d^{-1}\right)
\gk_{\bracket{\xi'}}
=\gk_{\bracket{\xi'}}^{-1} d^{-1}\partial_{\xi_j}d\ d^{-1}
\gk_{\bracket{\xi'}}
=O(\bracket{\xi'}^{-1})$$
in $\cL(L^2(\R_+))$.
Iteration shows that 
\begin{eqnarray}\label{2.8.4}
(1+(c^*)^{-1}c^{-1})^{-1} \in S^{0}(\R^{n-1}\times\R^{n-1}; L^2(\R_+), L^2(\R_+)).
\end{eqnarray}
The smoothness in $(x',\xi')$ and equations \eqref{2.8.1} -- \eqref{2.8.4}
then imply that 
$$(1+c^*c)^{-1}\in S^{-2m}(\R^{n-1}\times\R^{n-1};H^{-m}_0(\R_+), H^m(\R_+)).$$
\eproof

%%%%%%%%%%%%%%%%%%%%%%%%%%%%%%%%%%%%%%%%%%%%%%%%%%%%%%%%%%%%%%%%%%%%%%%

\thm{2.8a}{$(1+c^*c)^{-1}$ is a boundary symbol operator in Boutet de
Monvel's calculus whose pseudodifferential part is 
$r^+_{x_n}(1+p^*p)|_{x_n=0}^{-1}$.} 

%%%%%%%%%%%%%%%%%%%%%%%%%%%%%%%%%%%%%%%%%%%%%%%%%%%%%%%%%%%%%%%%%%%%%%%
\Proof.
We start with a few preliminaries.
We consider the composition $c^*(x',\xi')c(x',\xi')$
with $c$ defined in \eqref{2.7.1}.
For fixed $(x,\xi)$, the adjoint 
$$c^*(x,\xi):H^0_0(\ol\R_+)\cong L^2(\R_+)\to H^{-m}_0(\ol\R_+)$$ 
is given by  
$$c^*(x',\xi')= p^*(x',0,\xi',D_n)+g^*(x',\xi')$$
with the formal adjoints of $p$ and $g$. 

This, however, needs some explanation. In order to keep
the notation light, we will simply write $c^*=p^*(D_n)+g^*$.
The adjoint $g^*$ is of order $m$ and class $0$. 
It naturally maps $L^2(\R_+)$
to $\cS(\R_+)$, which can be viewed as a subspace of 
$H^{-m}_0(\ol\R_+)$ via extension by zero.

In order to apply $p^*(D_n)$ to $u\in L^2(\R_+)$, we first extend $u$
by zero to an element of $L^2(\R)$. Applying $p^*(D_n)$ furnishes a
distribution in $H^{-m}(\R)$, which will  in general
not vanish on $\R_-$: It only induces a functional on 
$H^m(\R_+)$ coming from an element in $H^{-m}_0(\ol\R_+)$.
Indeed, according to \eqref{B.4.1}, we can write 
$p^*=s^*+q^*$ with a polynomial $s^*$ and $q^*$ of normal order $-1$. 
Then $q^*(D_n)$ maps $e^+L^2(\R_+)$ to $L^2(\R)$ so that its
restrictions to the positive and the negative half-line  are
defined. In view of the fact 
that $s^*(D_n)$ preserves the support, 
%we see that, for $u$ in $e^+L^2(\R_+)$ 
the distribution $p^*(D_n)u$ as an element in $H^{-m}_0(\ol\R_+)$ is 
given by 
$$p^*(D_n)u-r^-q^*(D_n)u=p^*(D_n)u-e^-Jg^-(q^*)u$$
with the symbol $g^-(q^*)$  introduced in \eqref{gminus}. 
It is of order $m$ and class $0$.

We therefore have
\begin{eqnarray}
c^*c&=&p^*(D_n)e^+(p^+(D_n)+g)+e^+g^*(p^+(D_n)+g)
\nonumber\\
&&-e^-Jg^-(q^*)(p^+(D_n)+g)\nonumber\\
&=&p^*(D_n)e^+(p^+(D_n)+g)+e^+g_1+e^-Jg_2\label{2.8a.5} 
\end{eqnarray}
with suitable singular Green symbols $g_1$ and $g_2$ of 
orders $2m$ and class $m$. Hence (with $p=p(x',0,\xi)$)
\begin{eqnarray*}
\lefteqn{r^+(1+p^*p)^{-1}(D_n)\, c^*c}\\
 &=& r^+((1+p^*p)^{-1}p^*)(D_n)\,e^+p^+(D_n)
+\ r^+((1+p^*p)^{-1}p^*)(D_n)\,e^+g\\
&&+ \ r^+(1+p^*p)^{-1}(D_n)\,e^+g_1
+\ r^+(1+p^*p)^{-1}(D_n)\,e^-Jg_2.
\end{eqnarray*}
Surprisingly, all the terms on the right hand side can be treated
within Boutet de Monvel's calculus. The first is the 
composition of two truncated pseudodifferential operators. It equals
$$((1+p^*p)^{-1}p^*p)^+(D_n)-l((1+p^*p)^{-1}p^*,p),$$
where the leftover term of the composition is
of order $0$ and class $m$.
The second and the third are compositions of a truncated 
pseudodifferential operator with a singular Green symbol,
thus  singular Green symbols. % operator. 
Both have order zero and class $m$. 
The final term is the composition of the $g^+$-term of the 
pseudodifferential part with $g_2$ and therefore also a 
singular Green symbol of order $0$ and class $m$.

Putting all this together, we find that 
\begin{eqnarray*}
r^+(1+p^*p)^{-1}(D_n)(1+c^*c)= 1+g_3
\end{eqnarray*}
with a singular Green symbol $g_3$ of order $0$ and class $m$.
Hence 
\begin{eqnarray}\label{2.8a.3}
(1+c^*c)^{-1}=r^+(1+p^*p)^{-1}(D_n)-g_3(1+c^*c)^{-1}.
\end{eqnarray}

%%%%%%%%%%%%%%%%%%%%%%%%%%%%%%%%%%%%%%%%%%%%%%%%%%%%%%%%%%%%%%%%%%%%%%%

Let us now have a look at 
$c^*c\,r^+((1+p^*p)^{-1})(D_n)$. 
%We consider this composition on $\cC^\infty(\ol\R_+)$
%which is dense in $H^{-m}_0(\ol\R_+)$. We
We first note that we may consider $e^+$ a trivial action on 
$H^{-m}_0(\ol\R_+)$  and that 
\begin{eqnarray}
%\lefteqn{
p^+(D_n)r^+(1+p^*p)^{-1}(D_n)e^+%}\nonumber\\
&=&
(p(1+p^*p)^{-1})^+(D_n)+ \label{2.8a.4}
l(p,(1+p^*p)^{-1}),
\end{eqnarray}
where the leftover term on the right hand side is of
order $-m$ and class $0$. 
As the first term on the right hand side
% is of order $-m$, 
maps to $L^2(\R)$ we can 
rewrite it as 
\begin{eqnarray}
\lefteqn{(p(1+p^*p)^{-1})(D_n)e^+ - 
e^-r^-(p(1+p^*p)^{-1})(D_n)e^+}\nonumber\\
&=& (p(1+p^*p)^{-1})(D_n)e^+ -e^-Jg^-(p(1+p^*p)^{-1}),
\label{2.8a.6}
\end{eqnarray}
where the $g^-$-term is of order $-m$ and class $0$.

Taking into account \eqref{2.8a.5}, \eqref{2.8a.4} and
\eqref{2.8a.6}
\begin{eqnarray*}
\lefteqn{(1+c^*c)r^+(1+p^*p)^{-1}(D_n)}\\
&=&
1 + p^*(D_n)\left( e^+g_4+e^-Jg_5\right)
+e^+g_6+e^-Jg_7
\end{eqnarray*}
with singular Green symbols $g_4$ and $g_5$ of order
$-m$ and class $0$ and $g_6$ and $g_7$ of order and class $0$. 
Note that the image of the sum of the second and the fourth 
summand necessarily lies in $H^{-m}_0(\ol\R_+)$, 
since this is the case for the others.

Combining this with \eqref{2.8a.3} we conclude that 
\begin{eqnarray*}
\lefteqn{(1+c^*c)^{-1}}\nonumber\\
&=& r^+(1+p^*p)^{-1}(D_n)-
(1+c^*c)^{-1}p^*(D_n)\left( e^+g_4 
+ e^-Jg_5\right)\\
&&-(1+c^*c)^{-1} \left(e^+g_6+e^-Jg_7\right)\nonumber\\
&=&r^+(1+p^*p)^{-1}(D_n)
-r^+(1+p^*p)^{-1}(D_n)p^*(D_n)e^+g_4\\
&&-r^+(1+p^*p)^{-1}(D_n)p^*(D_n)e^-Jg_5
-r^+(1+p^*p)^{-1}(D_n)\left(e^+g_6+e^-Jg_7\right) 
%ES was an error -g_3r^+(1+p^*p)^{-1}(D_n)
\\
&&+g_3(1+c^*c)^{-1}p^*(D_n)\left( e^+g_4+e^-Jg_5\right)+g_3(1+c^*c)^{-1}\left(e^+g_6+e^-Jg_7\right).
\end{eqnarray*}
The first term on the right hand side is the one
we want as the pseudodifferential part. 
The second is the composition of a truncated pseudodifferential
operator of order $-m$ with a singular Green symbol
of order $-m$ and class $0$, thus a singular Green 
symbol of order $-2m$ and class $0$. 
The third is the composition of a $g^+$-type symbol 
of order $-m$ and class $0$ with a 
singular Green symbol of order $-m$ and class $0$,
thus of the same type as the second.
The summands of the fourth term are of the same type as the second 
and the third. 
%ES consequence of omission above
\forget{As $e^+$ is a trivial action on $H^{-m}_0$, the fifth is 
the composition of a singular Green symbol
of order $-m$ and class $m$ with a truncated pseudodifferential
operator of order $-m$, thus a singular Green symbol of order $-2m$ and
class $0$. }
As for the sum of the fifth and sixth, we note that 
\begin{eqnarray*}
e^+g_4+e^-Jg_5&\in& S^{-m}(\R^{n-1}\times\R^{n-1};
\cS'(\R_+),\cS(\R_+)\oplus\cS(\R_-))\text{ and}\\
p^*(D_n)&\in& S^{m} (\R^{n-1}\times\R^{n-1};\cS(\R_+)\oplus\cS(\R_-), H^{-m}(\R)).
\end{eqnarray*}
The composition of both therefore is an element of 
$S^{0} (\R^{n-1}\times\R^{n-1};\cS'(\R_+), H^{-m}(\R))$.
As $\cS(\R_+)\oplus \cS(\R_-)\hookrightarrow H^{-m}(\R)$, we have
$$e^+g_6+e^-Jg_7\in 
S^{0} (\R^{n-1}\times\R^{n-1};\cS'(\R_+), H^{-m}(\R)) 
$$
Moreover, we know that the range of the sum of all these terms is in 
$H^{-m}_0(\ol\R_+)$ so that we
can replace $H^{-m}(\R)$ in the symbol space by $H^{-m}_0(\ol\R_+)$.
We saw in Proposition \ref{2.8} that 
\begin{eqnarray*}
(1+c^*c)^{-1}&\in&S^{-2m} (\R^{n-1}\times\R^{n-1};H^{-m}_0(\R_+), H^m(\R_+)).
\end{eqnarray*}
As $g_3$ is of order $0$ and class $m$ we have 
\begin{eqnarray*}
g_3&\in& S^{0} (\R^{n-1}\times\R^{n-1}; H^m(\R_+),\cS(\R_+)).
\end{eqnarray*}
Hence the total composition is an element of 
$S^{-2m} (\R^{n-1}\times\R^{n-1}; \cS'(\R_+),\cS(\R_+))$ 
thus a singular Green symbol of order $-2m$ and class $0$.   
This shows that 
$(1+c^*c)^{-1}$ is a boundary symbol operator in Boutet de Monvel's
calculus which differs from $r^+(1+p^*p)^{-1}(D_n)$ by a 
singular Green symbol of order $-2m$ and class $0$.  
\eproof

%%%%%%%%%%%%%%%%%%%%%%%%%%%%%%%%%%%%%%%%%%%%%%%%%%%%%%%%%%%%%%%%%%%%%%%
\extra{2.8j}{The inverse of $1+cc^*$}{For fixed $(x',\xi')$, the
operator $c(x',\xi')c^*(x',\xi'):
H^{m}_0(\ol\R_+)\to H^{-m}(\R_+)$ acts on 
$v\in \cC^\infty_c(\R_+)$ by considering the function
$c^*(x',\xi')v$ as an element of $H^0_0(\ol\R_+)\cong L^2(\R_+)$
to which we then apply $c$. We know that $c^*=p^*(D_n)+ g^*$
with the notation introduced above. Hence $c^*v$ is a function 
in $\cS(\R)$; interpreting it as a distribution in $H^0_0(\ol\R_+)$
amounts to restricting it to $\R_+$. As $e^+$ can be considered a
trivial action on $\cC^\infty_c(\R_+)$, the action of $c^*$ 
coincides with that of $p^{*,+}(D_n)+g^*$.
The composition $cc^*$ therefore coincides with the composition of
two boundary symbol operators in Boutet de Monvel's calculus.
As $1+cc^*$ is invertible, the inverse is also given by a boundary
symbol in that calculus, and we obtain the statement below:   }

%%%%%%%%%%%%%%%%%%%%%%%%%%%%%%%%%%%%%%%%%%%%%%%%%%%%%%%%%%%%%%%%%%%%%%%

\cor{2.8k}{The inverse $(1+cc^*)^{-1}$ has the following form:  
\begin{eqnarray}\label{2.8k.1}
(1+cc^*)^{-1}&=& r^+(1+pp^*)^{-1}(D_n)+g_{8}
\end{eqnarray}
with a singular Green symbol $g_{8}$ of order $-2m$ and class $0$. }

%%%%%%%%%%%%%%%%%%%%%%%%%%%%%%%%%%%%%%%%%%%%%%%%%%%%%%%%%%%%%%%%%%%%%%%
\prop{2.10}{By  
$$\Gr(p)= \Mat{(1+p^*p)^{-1}& (1+p^*p)^{-1}p^*\\ 
p(1+p^*p)^{-1}& p(1+p^*p)^{-1}p^*}\in \cC_0(T^*X,\cL(E_1\oplus E_2))^\sim$$
we denote the graph projection of $p$.
The difference of equivalence classes 
$$[\Gr(p)]-\left[\Mat{0&0\\0&1}\right]$$ 
then defines an element in $K_0(\cC_0(T^*X,\cL(E_1\oplus E_2)))$ which 
is independent of the way the smoothing near zero is performed.}

%%%%%%%%%%%%%%%%%%%%%%%%%%%%%%%%%%%%%%%%%%%%%%%%%%%%%%%%%%%%%%%%%%%%%%%

\Proof. Let $p_0$ and $p_1$ be two smooth extensions of $p^{m}$. 
We let 
$p_t=(1-t)p_0+tp_1$, $0\le t\le 1$. For each $t$,  
$p_t$ is a smooth function on
$T^*X$ which coincides with $p^{m}$ on $\{|\xi|\ge 1\}$.  
The associated family  $\Gr(p_t)$ is continuous in $t$;
hence the class $[\Gr(p_t)]$ is constant.
\eproof    
%%%%%%%%%%%%%%%%%%%%%%%%%%%%%%%%%%%%%%%%%%%%%%%%%%%%%%%%%%%%%%%%%%%%%%%

\prop{2.15}{For the graph projection of $c$,   
$$\Gr(c)= \Mat{(1+c^*c)^{-1}& (1+c^*c)^{-1}c^*\\ 
c(1+c^*c)^{-1}& c(1+c^*c)^{-1}c^*}\in
\cC_0(T^*\partial X,\cL(L^2(\R_+,\tilde E_1\oplus \tilde E_2)))^\sim$$
the difference of equivalence classes 
$$[\Gr(c)]-\left[\Mat{0&0\\0&1}\right]$$ 
defines an element in 
$K_0(\cC_0(T^*\partial X,\cL(L^2(\R_+,\tilde E_1\oplus \tilde E_2))))$. 
It does not depend on the way the smoothing near zero in 
\ref{2.7} is performed.
% nor on the choice of the order reduction.
}

%%%%%%%%%%%%%%%%%%%%%%%%%%%%%%%%%%%%%%%%%%%%%%%%%%%%%%%%%%%%%%%%%%%%%%%

\Proof. Let $p_0$ and $p_1$ be two smooth extensions of $p^{m}$ and
$g_0$ and $g_1$ two smooth extensions of $g^{m}$. We let 
$p_t=(1-t)p_0+tp_1$, $g_t=(1-t)g_0+tg_1$ and
$$c_t(x',\xi')=p_t^+(x',\xi',D_n)+g_t(x',\xi'),\quad0\le t\le 1.$$
The associated graph projections $\Gr(c_t)$ depend continuously on $t$ 
in the topology of 
$\cC_0(T^*\partial X,\cL(L^2(\R_+,\tilde E_1\oplus \tilde E_2)))^\sim$.
Hence the class is independent of $t$.
\eproof
%Moreover, \TEXT{missing part}\eproof

%%%%%%%%%%%%%%%%%%%%%%%%%%%%%%%%%%%%%%%%%%%%%%%%%%%%%%%%%%%%%%%%%%%%%%%
%ES Johannes
\forget{\rem{2.15a}{With the same argument and actually with the same addition
of vector bundles as in Remark \ref{2.4a} we may consider $[\Gr(p)]-
[e]$ and $[\Gr(c)]-[e]$ as elements of $K_0(\cC_0(T^*X))$ and 
$K_0(\cC_0(T^*\partial X,\cL(L^2(\R_+))))$.
}
}%forget
%%%%%%%%%%%%%%%%%%%%%%%%%%%%%%%%%%%%%%%%%%%%%%%%%%%%%%%%%%%%%%%%%%%%%%%
\section{The Tangent Semigroupoid}
%%%%%%%%%%%%%%%%%%%%%%%%%%%%%%%%%%%%%%%%%%%%%%%%%%%%%%%%%%%%%%%%%%%%%%%

We recall a few concepts from Aastrup-Nest-Schrohe \cite{ANS}
\dfn{t.1}{By $T^\pm X$ 
we denote the set of all vectors  $(x,v)\in T\widetilde X|_{X}$ 
for which $\exp_x(\pm \gve v)\in X$ for sufficiently small $\gve>0$.  
This is a semi-groupoid with addition of vectors, and
 $T^\pm X=TX^\circ\cup T^\pm X|_{\partial X}$ 

We define $\cT^-X$ as the disjoint union 
$T^-X \cup (X\times X\times ]0,1])$, 
endowed with the fiberwise semi-groupoid structure 
induced by the semi-groupoid structure on $T^-X$ 
and the pair groupoid structure on
$X\times X$.
We glue $T^-X$ to $X\times X \times ]0,1]$ via the charts
$$T^-X\times [0,1]\supseteq U\ni (x,v,\hbar)\mapsto \left\{ 
\begin{array}{lc}
(x,v) & \hbox{for }\hbar =0 \\
(x, \exp_x(-\hbar v),\hbar)& \hbox{for }\hbar \not= 0
\end{array}
\right.$$
and
let $\cT^-X(0)=T^-X$ and $\cT^-X(\hbar)=X\times X\times \{\hbar\}$.

In order to avoid problems with the topology of $\cT^-X$ 
(which is in general not a manifold with corners)
we let $\cC^\infty_c(\cT^-X)=\cC^\infty_c(\cT\widetilde{X})|_{\cT^-X}$.
}

%%%%%%%%%%%%%%%%%%%%%%%%%%%%%%%%%%%%%%%%%%%%%%%%%%%%%%%%%%%%%%%%%%%%%%%

\subsection*{C*-algebras Associated to the Semi-groupoids 
$T^-X$ and $\cT^- X$}
Let  $ \cC^\infty ( T^-X)$ denote the smooth functions on $T^-X$ 
with compact support in $T^-X$.
We introduce  
\begin{eqnarray*}
\pi_0&:&\cC^\infty(T^-X)\to \cL(L^2(T X))\quad\text{and}\\
\pi_0^\partial&:&\cC^\infty(T^-X)\to \cL(L^2(T^+X|_{\partial X}))
\end{eqnarray*} 
acting by 
\begin{eqnarray}
\pi_0(f)\xi(x,v)=\int_{T_mX}f(x,v-w)\xi(x,w)\, dw,\\
\pi_0^\partial(f)\xi (x,v)=\int_{T^+_xX}f(x,v-w)\xi(x,w)\, dw.
\label{kappaf}
\end{eqnarray} 
Note: As $f$ has compact support in $T^-X$, it 
naturally extends (by zero) to $TX$.

%%%%%%%%%%%%%%%%%%%%%%%%%%%%%%%%%%%%%%%%%%%%%%%%%%%%%%%%%%%%%%%%%%%%%%%

\dfn{Cr}{$C^*_r(T^-X)$ is the $C^*$-algebra generated by 
$\pi_0\oplus \pi_0^\partial$, i.e.~by the map
$$\cC^\infty(T^-X)\ni f\mapsto  
(\pi_0(f),\pi_0^\partial(f))\in 
\cL(L^2(TX)\oplus L^2(T^+X|_{\partial X})).$$
}
%%%%%%%%%%%%%%%%%%%%%%%%%%%%%%%%%%%%%%%%%%%%%%%%%%%%%%%%%%%%%%%%%%%%%%%

\rem{Cralt}{According to  Lemmas 2.14 and 2.15 in \cite{ANS}, 
$C^*_r(T^-X)$ has the dense $*$-sub\-al\-gebra
$$\cC^\infty_{tc}(T^-X)= \cC^\infty_c(TX)\oplus
\cC^\infty_c(T\partial X\times\ol\R_+\times \ol\R_+)$$
with the representation $\pi_0\oplus \pi_0^\partial$ of 
$\cC^\infty_c(TX )$  on $L^2(T X)\oplus L^2(T^+X|_{\partial X})$, 
defined as above,  and 
the representation 
$\tilde \pi_0^\partial$ of 
$\cC^\infty_c(T\partial X\times\ol\R_+\times \ol\R_+)$
on $L^2(T^+X|_{\partial X})$ given by 
\begin{eqnarray*}\label{rep2}
\tilde \pi_0^\partial (K)\xi(x,v',v_n)=
\int_{T^+X|_{\partial X}} 
K(x,v'-w',v_n,w_n)\xi(x,w',w_n)\, dw'\,dw_n.
\end{eqnarray*}
}

%%%%%%%%%%%%%%%%%%%%%%%%%%%%%%%%%%%%%%%%%%%%%%%%%%%%%%%%%%%%%%%%%%%%%%%

\extra{t.3}{An ideal in $C^*_r(T^-X)$}{Denote by 
\begin{eqnarray*}
\cF:L^2(TX)&\to&L^2(T^*X)\quad\text{and}\\
\cF':L^2(T\partial X)&\to&L^2(T^*\partial X) 
\end{eqnarray*}
the fiberwise Fourier transforms. 

It was noted in \cite[Lemma 2.15]{ANS} that 
$\cF':L^2(T\partial X)\to L^2(T^*\partial X)$
provides an isomorphism between the ideal of  $C^*_r(T^-X)$
generated by the representation $\tilde\pi^\partial_0$ in  
\ref{Cralt} and $\cC_0(T^*\partial X,\cK)$.  

The Fourier transform allows us two more important identifications:
For $f\in \cC^\infty_c(TX)$, the operator 
%$$\cF \pi_0(f)\cF^{-1}:L^2(T^*X,E_1)\to L^2(T^*X,E_2)$$
$$\cF \pi_0(f)\cF^{-1}:L^2(T^*X)\to L^2(T^*X)$$
is the operator of multiplication by $\widehat f=\cF f$.

At the boundary, the choice of a Riemannian metric allows us to 
identify $T^\pm X|_{\partial X} $ with 
$T\partial X\times \R_{\pm}$ and
the operator 
$$\cF'\pi_0^\partial (f){\cF'}^{-1}:
%L^2(T^*\partial X\times \R_+,\widetilde E_1)\to 
%L^2(T^*\partial X\times \R_+,\widetilde E_2)$$ 
L^2(T^*\partial X\times \R_+)\to 
L^2(T^*\partial X\times \R_+)$$ 
with the boundary 
symbol operator $\cF'(f)(x',0,\xi',D_n)^+.$
}

%%%%%%%%%%%%%%%%%%%%%%%%%%%%%%%%%%%%%%%%%%%%%%%%%%%%%%%%%%%%%%%%%%%%%%%
%ES Johannes cosmetically revised
\rem{bundlevsscalar}{So far we have been working with the graph 
projections of operators and symbols acting in vector bundles. 
Since we want it to give rise to elements in the 
$K$-theory of $C^*_r(T^-X)$ and $C^*_r(\cT^- X)$ we will describe here 
how this can be achieved: 
Choose bundles $F_1$ and $F_2$ such that 
$E_1\oplus F_1$ and $E_2\oplus F_2$ are trivial with fiber $\C^\nu$. 
We can consider the graph projection of an operator $A$ acting between 
sections of $E_1$ and $E_2$ as an element of $\cL(L^2(X,\C^{2\nu}))$. 
Let $\pi_{F_2}$ be the projection onto $F_2$. Note that 
$$\Gr(A)+\pi_{F_2} $$
is a projection in $M_{2\nu}(\cK (L^2(X))^\sim)$
and that 
$$[\Gr(A)+\pi_{F_2}]-[\pi_{E_2\oplus F_2}] $$ 
is the index class of $A$ in $K_0(\cK(L^2 (X)))$.
The same consideration applies to the graph projection of the symbol, 
so that we get a class in $K_0(C^*(T^-X))$. 
Since the construction is stable under the $\hbar$-scaling this gives 
a way of passing from general vector bundles to trivial vector 
bundles.

We will in the rest of the paper use this to identify various graph 
projections acting in vector bundles with projections in trivial bundles.   
}

%%%%%%%%%%%%%%%%%%%%%%%%%%%%%%%%%%%%%%%%%%%%%%%%%%%%%%%%%%%%%%%%%%%%%%%

\prop{t.4}{%
With the identification provided by the Fourier transforms 
described in \ref{t.3}, 
 $\Gr(p)\oplus \Gr(c)$ can be regarded as an element in
$ M_N(C^*_r(T^-X)^\sim)$ for suitable $N$.
}

%%%%%%%%%%%%%%%%%%%%%%%%%%%%%%%%%%%%%%%%%%%%%%%%%%%%%%%%%%%%%%%%%%%%%%%

\Proof. We abbreviate $p(D_n)=p(x',0,\xi',D_n)$ 
and consider the difference 
\begin{eqnarray}\label{t.4.1}
\Gr(c)- r^+\, \Gr(p(D_n))\, e^+
\end{eqnarray}
between the graph projection of $c$ and the truncated operator 
obtained from the graph projection
of the operator $p(D_n)$, i.e.  from
$$\Gr(p(D_n))= 
\Mat{(1+p^*p)^{-1}(D_n)&(1+p^*p)^{-1}(D_n)p^*(D_n)\\
p(D_n)(1+p^*p)^{-1}(D_n)
&p(D_n)(1+p^*p)^{-1}(D_n)p^*(D_n)}$$
acting on $L^2(\R,\widetilde E_1\oplus\widetilde E_2)$.
According to Theorem \ref{2.8a}, % and corollary \ref{2.8k},  
$(1+c^*c)^{-1}-r^+(1-p^*p)^{-1}(D_n)$ is a 
singular Green symbol 
of negative order and thus a compact operator from $H^{-m}_0(\ol\R_+,
\widetilde E_1)$ to $H^m(\R_+,\widetilde E_1)$. Hence  
the difference \eqref{t.4.1} is an element of 
$M_N(\cC_0(T^*\partial X,\cK))$. 
As pointed out in \ref{t.3}, 
conjugation by the boundary Fourier transform maps it  
to an ideal in $M_N(C^*_r(T^-X))$.
We only have to show that conjugation by the 
Fourier transforms maps 
$$\Gr(p)\oplus r^+\,\Gr(p(D_n))e^+$$ 
to an element of $M_N(C^*_r(T^-X)^\sim).$

In order to do this, it actually suffices to find a sequence 
$(f_k)$ in $\cC^\infty_c(TX,\cL(\widetilde E_1\oplus\widetilde E_2))$ 
which converges to $f=\cF^{-1}(\Gr(p)-e)$, where $e$ is the usual projection onto
the second component, with 
respect to the norm $g\mapsto \sup_{x\in X}\|g(x,\cdot)\|_{L^1}$.
Indeed, this will imply that 
$\pi_0(f_k)\to \pi_0(f)$ or, 
equivalently, that as multiplication operators 
%on $L^2(\R,\widetilde E_1\oplus\widetilde E_2)$ 
$\cF f_k\to \Gr(p)-e.$ 
Moreover, as   $\|\pi_0^\partial (\cdot)\|$ is dominated by 
$\|\pi_0(\cdot)\| $, 
also the operators $r^+\cF'( f_k)(D_n)e^+$ will approximate 
$r^+\,\Gr(p(D_n))\, e^+-e$.

It remains to find such a sequence. 
To this end we note that the entries in 
$\Gr(p)-e$
are symbols of orders $\le -m<-n$. 
%As $m>n$ by assumption, they are integrable along the fibers, 
%and their inverse
%Fourier transforms are continuous on $TX$. 
%Moreover, the symbol property implies that for a symbol 
%$q$ of order $\le-n$
%we have in local coordinates 
Hence, for any  $K\ge0$, they satisfy 
$$\sup_{x,v}\left|(1+|v|^2)^K \,(\cF^{-1}q)(x,v)\right|
\le\sup_{x,v} \left|\int e^{iv\xi} (1-\Delta_{\xi})^K
q(x,\xi)\,\dbar\xi\right|<\infty.$$
As $\cC^\infty_c(TX)$  
is dense in the space of functions for which the weighted sup-norm 
on the left hand side is finite and as this norm is larger than 
$\sup_x\|g(x,\cdot)\|_{L^1}$ whenever $K >N/2$, 
we find the desired sequence.
\eproof
\bigskip

%%%%%%%%%%%%%%%%%%%%%%%%%%%%%%%%%%%%%%%%%%%%%%%%%%%%%%%%%%%%%%%%%%%%%%%
%%%%%%%%%%%%%%%%%%%%%%%%%%%%%%%%%%%%%%%%%%%%%%%%%%%%%%%%%%%%%%%%%%%%%%%

\extra{pcoutline}{Parametrix Construction}{%
We fix a collar neighborhood $\partial X\times [0,2)$ of $\partial X$
in $X$. 
Moreover, we choose a covering of $X$ by open sets 
$X_j$ and  coordinate maps $\chi_j:U_j\subset\R^n\to X_j$. 
We assume that a coordinate neighborhood either is 
contained in a collar neighborhood of the boundary (`boundary chart') 
or else
does not intersect a neighborhood of the boundary (`interior chart').
It is also no restriction to suppose that 
the boundary charts all lie in the collar
neighborhood and that the variable normal 
to the boundary, 
$x_n$, is fixed and that changes of coordinates only involve the 
tangential variables.

We fix a partition of unity $\gvp_j$ subordinate to the coordinate
charts and cut-off functions $\psi_j$ supported in the same charts 
with $\psi_j(x)\equiv 1$ in a neighborhood of the support of $\gvp_j$.

}

%%%%%%%%%%%%%%%%%%%%%%%%%%%%%%%%%%%%%%%%%%%%%%%%%%%%%%%%%%%%%%%%%%%%%%%

%%%%%%%%%%%%%%%%%%%%%%%%%%%%%%%%%%%%%%%%%%%%%%%%%%%%%%%%%%%%%%%%%%%%%%%
\lemma{c.2}{Without changing the $K$-classes of $\cG(p)$ and $\cG(c)$
we may assume that the symbol $p$ of \ref{2.7} 
is independent of $x_n$ in a neighborhood of $\partial X$.
In particular, we shall assume that this is the case on the collar 
neighborhood  $\partial X\times[0,1)$ of the boundary. 
}

%%%%%%%%%%%%%%%%%%%%%%%%%%%%%%%%%%%%%%%%%%%%%%%%%%%%%%%%%%%%%%%%%%%%%%%

\Proof. Apply a smooth deformation of the boundary variable.
This will imply continuity of the change of the graph 
projections and thus keep the associated $K$-class constant.
\eproof 

%%%%%%%%%%%%%%%%%%%%%%%%%%%%%%%%%%%%%%%%%%%%%%%%%%%%%%%%%%%%%%%%%%%%%%%

\extra{c.3}{Scaling}{According to \ref{1.15} we define
\begin{eqnarray*}
p_\hbar(x,\xi) &=& p(x,\hbar\xi),\\
g_\hbar(x',\xi')&=& \gk_{\hbar}^{-1} g(x',\hbar\xi')\gk_\hbar,\\
c_\hbar(x',\xi') &=& 
\gk_{\hbar}^{-1} \op^+_{x_n}(p|_{x_n=0})(x',\hbar\xi')\gk_{\hbar}
+g_\hbar(x',\xi')
\end{eqnarray*} 

We denote by $p^j$ the function $p$ in the $\chi_j$-coordinates and 
by $\chi_j^*$ the transport of operators from $U_j$
to $X_j$, i.e. $\chi^*_j\op p^j_\hbar$ is the operator induced on 
$X$ from the operator $\op p^j_\hbar$ on $\R^n$. 
Of course, this only makes sense
when multiplied with suitable cut-off functions from the left and 
the right.  

In the boundary charts we will use the $\hbar$-scaled 
boundary symbol operators $c_\hbar$; we write $c_\hbar^k$ for this
symbol in the $\chi_k$-coordinates  and 
$\op'$ for the quantization map for boundary symbol  operators.

Then we define the operator family 
$A_\hbar:\cC^\infty(X,E_1)\to \cC^\infty(X,E_2)$ by 
\begin{eqnarray*}
A_\hbar=\sum_{\text{boundary charts}}
\gvp_k \chi_k^{*}\op'(c^k_\hbar)\psi_k + 
\sum_{\text{interior charts}}\gvp_j\chi_j^*\op(p^j_\hbar)\psi_j,
=A_{b,\hbar} + A_{i,\hbar}
\end{eqnarray*}
consisting of a boundary part $A_{b,\hbar}$ and an 
interior part $A_{i,\hbar}$.

 }
%%%%%%%%%%%%%%%%%%%%%%%%%%%%%%%%%%%%%%%%%%%%%%%%%%%%%%%%%%%%%%%%%%%%%%%

\lemma{c.4}{Let $\go,\go_1$ be smooth and supported in a single 
boundary  neighborhood. Then 
$$\go \op^+(p_\hbar)\go_1 
= \go\gk_{\hbar}^{-1}\op' ( 
(\op_{x_n}^+p|_{x_n=0})(x',\hbar\xi'))\gk_\hbar\go_1. $$}

%%%%%%%%%%%%%%%%%%%%%%%%%%%%%%%%%%%%%%%%%%%%%%%%%%%%%%%%%%%%%%%%%%%%%%%

\Proof. This follows from the fact that $p$ is independent of $x_n$ 
on the collar neighborhood and the computation in  \ref{h.1}.
\eproof

%%%%%%%%%%%%%%%%%%%%%%%%%%%%%%%%%%%%%%%%%%%%%%%%%%%%%%%%%%%%%%%%%%%%%%%
\lemma{c.3z}{$A_\hbar $ is an operator family in Boutet de Monvel's 
calculus. Its pseudodifferential symbol is of the form  
$p_\hbar+ \hbar q(\hbar)$, where $q_{1/\hbar}$ 
is bounded in $S_{\tr}^{m-1}$. 
Over a boundary chart, its singular Green symbol is of the form 
$g_\hbar+\hbar r(\hbar)$ with $r_{1/\hbar}$ bounded in 
$S^{m-1}(\R^{n-1}\times\R^{n-1}; \cS'(\R_+),\cS(\R_+))$.}
%%%%%%%%%%%%%%%%%%%%%%%%%%%%%%%%%%%%%%%%%%%%%%%%%%%%%%%%%%%%%%%%%%%%%%%
\Proof. This follows from \ref{1.30b}, in particular
\eqref{leading}, in connection with Lemma \ref{c.4}. \eproof 
\medskip 
%%%%%%%%%%%%%%%%%%%%%%%%%%%%%%%%%%%%%%%%%%%%%%%%%%%%%%%%%%%%%%%%%%%%%%%

We shall say that a family $C(\hbar)$ of operators 
given by pseudodifferential operators in the interior and 
operator-valued symbols close to the boundary
%in Boutet de Monvel's calculus 
is semiclassically bounded of order $\mu$, if, for
the family $q(\hbar)$ of pseudodifferential symbols, 
$q_{1/\hbar}$ is bounded in $S^\mu$,  
and, for the family $d(\hbar)$ of operator-valued symbols,  
$d_{1/\hbar}$ is bounded in 
$S^\mu(\R^{n-1}\times\R^{n-1};E,F),$
where $E$ and $F$ have to be specified.

%%%%%%%%%%%%%%%%%%%%%%%%%%%%%%%%%%%%%%%%%%%%%%%%%%%%%%%%%%%%%%%%%%%%%%%
\prop{c.3y}{We have 
\begin{eqnarray*}
\lefteqn{1+A_\hbar^*A_\hbar 
= \sum_{\text{boundary charts}}
\gvp_k \chi_k^{*}\op'(1+c^{k*}_\hbar c^k_\hbar)\psi_k} \\
&&+ 
\sum_{\text{interior charts}}
\gvp_j\chi_j^*\op(1+p^{j*}_\hbar p^j_\hbar)\psi_j
+\hbar R_1(\hbar)\\
&=&\sum_{\text{boundary charts}}
\gvp_k \chi_k^{*}\op'(1+c^{k*} c^k)_\hbar\psi_k \\
&&+ 
\sum_{\text{interior charts}}
\gvp_j\chi_j^*\op(1+p^{j*}p^j)_\hbar\psi_j
+\hbar R_1(\hbar),
\end{eqnarray*}
where the family $R_1$ is semiclassically bounded of order $2m-1$
with operator-valued symbols acting between $E=H^{m}(\R_+)$ and
$F=H^{-m}_0(\ol\R_+)$.

{\rm Here, $p^{j,*}$ is the adjoint of the symbol $p^j$ in 
local coordinates, and  $c^{k,*}$ is the adjoint of the 
operator-valued symbol $c^{k}$.}}

%%%%%%%%%%%%%%%%%%%%%%%%%%%%%%%%%%%%%%%%%%%%%%%%%%%%%%%%%%%%%%%%%%%%%%%
\Proof. For the first identity apply \ref{1.15}(c) using Lemma 
\ref{c.3z}. The second is obvious.\eproof  
%%%%%%%%%%%%%%%%%%%%%%%%%%%%%%%%%%%%%%%%%%%%%%%%%%%%%%%%%%%%%%%%%%%%%%%

\prop{c.3x}{Let $\gvp,\psi\in \cC^\infty_c(\partial X\times (0,1))$
be supported in the intersection of an interior and a boundary 
neighborhood. Then
\begin{eqnarray*}
\gvp\op'(1+c^*c)^{-1}_\hbar\psi - \gvp\op(1+p^*p)^{-1}_\hbar\psi
\end{eqnarray*}
is a semiclassically regularizing pseudodifferential operator; i.e. for 
each $N$, we can write it as $\hbar^N \op r^N(\hbar)$ with
$r^N$ bounded in $S^{-2m-N}$. }
%%%%%%%%%%%%%%%%%%%%%%%%%%%%%%%%%%%%%%%%%%%%%%%%%%%%%%%%%%%%%%%%%%%%%%%

\Proof. We know  that  $(1+c^*c)^{-1}$ is a 
boundary symbol in Boutet de Monvel's calculus, whose 
pseudodifferential part is given by $(1+p|_{x_n=0}^*p|_{x_n=0})^{-1}$.
As the localization of the $\hbar$-scaled singular Green part to the 
interior is semiclassically smoothing by Lemma \ref{1.15b}, 
this implies the assertion.
\eproof

%%%%%%%%%%%%%%%%%%%%%%%%%%%%%%%%%%%%%%%%%%%%%%%%%%%%%%%%%%%%%%%%%%%%%%%

\prop{c.3w}{Define
\begin{eqnarray}\label{c.3x.1}
B(\hbar) &=&\sum_{\text{boundary charts}}
\gvp_k\chi_k^{*}\op'(1+c^{k*} c^k)^{-1}_\hbar\psi_k \\
&&+ 
\sum_{\text{interior charts}}
\gvp_j\chi_j^*\op(1+p^{j*}p^j)^{-1}_\hbar\psi_j\nonumber
\end{eqnarray}
This is an operator family in Boutet de Monvel's calculus 
which is semiclassically bounded of order $-2m$, and 
\begin{eqnarray*}
B(\hbar)(1+A_\hbar^*A_\hbar)&=& 1 +\hbar R_2( \hbar )\quad\text{and}\\
(1+A_\hbar^*A_\hbar)B(\hbar)&=& 1 +\hbar R_3( \hbar )
\end{eqnarray*}
with families $R_2$, $R_3$ which are semiclassically bounded of order 
$-1$ with operator-valued symbols acting on $H^{m}(\R_+)$ for 
$R_2$ and on $H^{-m}_0(\ol\R_+)$ for $R_3$.}
%%%%%%%%%%%%%%%%%%%%%%%%%%%%%%%%%%%%%%%%%%%%%%%%%%%%%%%%%%%%%%%%%%%%%%%
\Proof.  This follows from \ref{1.15}(c) in connection with
Propositions \ref{c.3y} and \ref{c.3x}.\eproof

%%%%%%%%%%%%%%%%%%%%%%%%%%%%%%%%%%%%%%%%%%%%%%%%%%%%%%%%%%%%%%%%%%%%%%%
\cor{c.3v}{We infer from Proposition \ref{c.3w} that 
\begin{eqnarray}\label{c.3v.1}
(1+A_\hbar^*A_\hbar)^{-1} =B(\hbar)- \hbar B(\hbar)R_3(\hbar) 
+ \hbar^2R_2(\hbar)(1+A^*_\hbar A_\hbar)^{-1}R_3(\hbar).
\end{eqnarray}
\forget{It is clear that $(1+A^*_\hbar A_\hbar)^{-1}$ is a bounded family 
on $L^2(X)$. Hence $(1+A^*_\hbar A_\hbar)^{-1}$ differs from 
$B(\hbar)$ by a term which tends to zero in $L^2(X)$ as 
$\hbar\to 0$.}
}

%%%%%%%%%%%%%%%%%%%%%%%%%%%%%%%%%%%%%%%%%%%%%%%%%%%%%%%%%%%%%%%%%%%%%%%
{}From this we want to deduce that $(1+A^*_\hbar A_\hbar)^{-1}$
differs from $B(\hbar)$ by a term which is $O(\hbar)$ in 
$\cL(L^2(X))$. So far, this is not obvious:  The boundary
symbol parts of $R_2$ and $R_3$ act on $H^m$ and $H^{-m}_0$, 
respectively, while we only can guarantee boundedness of the inverse on 
$L^2(X)$. 
To this end we make the following observations:
%%%%%%%%%%%%%%%%%%%%%%%%%%%%%%%%%%%%%%%%%%%%%%%%%%%%%%%%%%%%%%%%%%%%%%%

\lemma{c.9a}{Given $N\in\N$ we find an operator family 
$C(\hbar)$,
$0<\hbar\le1$, in Boutet de Monvel's calculus
such that 
$$C(\hbar)A_\hbar=1+S_N(\hbar)$$
with 
$C$  and $S_N$ 
semiclassically bounded of orders $-m$ and $-N$,
respectively, in Boutet de Monvel's calculus.}

%%%%%%%%%%%%%%%%%%%%%%%%%%%%%%%%%%%%%%%%%%%%%%%%%%%%%%%%%%%%%%%%%%%%%%%
\Proof. Apply a semiclassical parametrix construction in Boutet
de Monvel's calculus, using Proposition
\ref{h.3}.
\eproof

%%%%%%%%%%%%%%%%%%%%%%%%%%%%%%%%%%%%%%%%%%%%%%%%%%%%%%%%%%%%%%%%%%%%%%%
\lemma{c.9b}{The operator  families $(1+A^*_\hbar A_\hbar)^{-1}$, 
$(1+A_\hbar A^*_\hbar)^{-1}$, $A_\hbar(1+A^*_\hbar A_\hbar)^{-1}A_\hbar^*$,
$A_\hbar(1+A^*_\hbar A_\hbar)^{-1}$, and 
$(1+A^*_\hbar A_\hbar)^{-1}A_\hbar^*$
are (after continuous extension) uniformly bounded on the corresponding
$L^2$ spaces.
}
%%%%%%%%%%%%%%%%%%%%%%%%%%%%%%%%%%%%%%%%%%%%%%%%%%%%%%%%%%%%%%%%%%%%%%%
\Proof.
For the first two families the statement is obvious as their operator
norm is bounded by $1$. For the third we use that, by Lemma \ref{2.3}, 
$$A_\hbar(1+A^*_\hbar A_\hbar)^{-1}
A_\hbar^*=1-(1+A_\hbar A^*_\hbar)^{-1}.$$
%ES simpler 
For the fourth we note that on a Hilbert space, the norm of an operator
$T$ equals $\|T^*T\|^{1/2}$. We apply this to 
$T=A_\hbar(1+A^*_\hbar A_\hbar)^{-1}$. By \ref{fa.7}, 
\begin{eqnarray*}
T^*T= (1+A^*_\hbar A_\hbar)^{-1}A_\hbar^*
A_\hbar(1+A^*_\hbar A_\hbar)^{-1}=(1+A_\hbar^*A_\hbar)^{-1}-
(1+A_\hbar^*A_\hbar)^{-2},
\end{eqnarray*}
which is bounded. 
\forget{For the fourth we infer from the identity
$$A_\hbar(1+A_\hbar^*A_\hbar)= (1+A_\hbar A_\hbar^*)A_\hbar$$
that 
\begin{eqnarray}\label{c.9b.1}
(1+A_\hbar A_\hbar^*)^{-1} A_\hbar 
= A_\hbar(1+A_\hbar^*A_\hbar)^{-1}.
\end{eqnarray}
Since $A_\hbar(1+A_\hbar^*A_\hbar)^{-1}$ is the 
adjoint of $(1+A_\hbar^*A_\hbar)^{-1}A^*_\hbar$ by \ref{fa.7} and, 
on a Hilbert space, the norm of an operator
$T$ equals $\|T^*T\|^{1/2}$, equation \eqref{c.9b.1} implies that
\begin{eqnarray*}
((1+A_\hbar^*A_\hbar)^{-1}A_\hbar^*)^*
 ((1+A_\hbar^*A_\hbar)^{-1}A_\hbar^*) =
  A_\hbar(1+A_\hbar^*A_\hbar)^{-1}
  (1+A_\hbar^*A_\hbar)^{-1}A_\hbar^*\\
=(1+A_\hbar A_\hbar^*)^{-1} A_\hbar
(1+A_\hbar^*A_\hbar)^{-1}A_\hbar^*
\end{eqnarray*}
is bounded in view of the boundedness of both factors.}
Duality yields the boundedness of the fifth family.\eproof

%%%%%%%%%%%%%%%%%%%%%%%%%%%%%%%%%%%%%%%%%%%%%%%%%%%%%%%%%%%%%%%%%%%%%%%
\cor{9.c}{$(1+A^*_\hbar A_\hbar)^{-1}-B(\hbar)=O(\hbar)$ in 
$\cL(L^2(X,E_1))$.  
}
\Proof. In the third %ES
term on the right hand side of 
\eqref{c.3v.1} we can write 
\begin{eqnarray}\label{9.c.1}
\lefteqn{(1+A^*_\hbar A_\hbar)^{-1}
=C(\hbar)A_\hbar(1+A^*_\hbar A_\hbar)^{-1}A_\hbar^*C(\hbar)^*}\\
&&-S_N(\hbar)(1+A^*_\hbar A_\hbar)^{-1}A_\hbar^*
C(\hbar)^*-C(\hbar)A_\hbar(1+A^*_\hbar A_\hbar)^{-1}S_N(\hbar)^*
\nonumber
\end{eqnarray}
for some $N>m$.
By Lemma \ref{1.15f}, $R_2C$, $R_2S_N$
$C^*R_3$ and $S_N^*R_3$ are bounded on $L^2$, uniformly in 
$\hbar$. The assertion then follows from \eqref{c.3v.1}.\eproof

%%%%%%%%%%%%%%%%%%%%%%%%%%%%%%%%%%%%%%%%%%%%%%%%%%%%%%%%%%%%%%%%%%%%%%%
\cor{c.10a}{$A_\hbar(1+A^*_\hbar A_\hbar)^{-1}-A_\hbar B(\hbar)=
O(\hbar)$ in $\cL(L^2(X,E_1), L^2(X,E_2))$. %ES  
}

%%%%%%%%%%%%%%%%%%%%%%%%%%%%%%%%%%%%%%%%%%%%%%%%%%%%%%%%%%%%%%%%%%%%%%%

\Proof. We multiply equation \ref{c.3v.1} from the left by $A_\hbar$.
\forget{We deduce from Corollary \ref{c.3v} that 
\begin{eqnarray*}\label{c.10a.2}
\lefteqn{A_\hbar(1+A_\hbar^*A_\hbar)^{-1}}\\ &=&
A_\hbar B(\hbar)- \hbar A_\hbar B(\hbar)R_3(\hbar) 
+ \hbar^2 A_\hbar R_2(\hbar)(1+A^*_\hbar A_\hbar)^{-1}R_3(\hbar).
\end{eqnarray*}
}
The composition $A_\hbar B(\hbar)$ furni\-shes an operator family in 
Boutet de Monvel's calculus which is semiclassically bounded of order
$-m$. According to Lemma \ref{1.15f}, the operator family 
$A_\hbar B(\hbar)R_3(\hbar)$ is therefore uniformly bounded on $L^2$.
In the second term on the right hand side we substitute
according to  equation \eqref{9.c.1}. We note that $A_\hbar R_2C$ and 
$A_\hbar R_2 S_N$ are semiclassically bounded of order $0$, with
the operator-valued symbols acting on $L^2(\R_+)$. Hence
$A_\hbar R_2(\hbar)(1+A^*_\hbar A_\hbar)^{-1}R_3(\hbar)$
is uniformly bounded in $\cL(L^2)$. \eproof

%%%%%%%%%%%%%%%%%%%%%%%%%%%%%%%%%%%%%%%%%%%%%%%%%%%%%%%%%%%%%%%%%%%%%%%
In an analogous way we find:

%%%%%%%%%%%%%%%%%%%%%%%%%%%%%%%%%%%%%%%%%%%%%%%%%%%%%%%%%%%%%%%%%%%%%%%
\cor{c.10d}{$(1+A^*_\hbar A_\hbar)^{-1}A^*_\hbar - B(\hbar)A^*_\hbar=
O(\hbar)$ in $\cL(L^2(X,E_2), L^2(X,E_1))$.  %ES
}
%%%%%%%%%%%%%%%%%%%%%%%%%%%%%%%%%%%%%%%%%%%%%%%%%%%%%%%%%%%%%%%%%%%%%%%

\rem{c.10b}{We know from Lemma \ref{2.3} that, as operators on 
$L^2(X,E_2)$, 
$$
A_\hbar(1+A_\hbar^*A_\hbar)^{-1}A_\hbar^*=1-(1+A_\hbar A_\hbar^*)^{-1}.
$$
In order to determine the structure of the left hand side, we
construct, similarly as before, an approximate inverse 
$\widetilde B(\hbar)$ to $1+A_\hbar A_\hbar^*$, using the structure of
the boundary symbol operator $(1+cc^*)^{-1}$ determined in 
Corollary \ref{2.8k}.
We obtain the following result:
 }

%%%%%%%%%%%%%%%%%%%%%%%%%%%%%%%%%%%%%%%%%%%%%%%%%%%%%%%%%%%%%%%%%%%%%%%
\prop{c.10c}{
\begin{eqnarray}\label{c.10c.0}
\widetilde B(\hbar) &=&\sum_{\text{boundary charts}}
\gvp_k \chi_k^{*}\op'(1+c^{k} c^{k*})^{-1}_\hbar\psi_k \\
&&+ 
\sum_{\text{interior charts}}
\gvp_j\chi_j^*\op(1+p^{j}p^{j*})^{-1}_\hbar\psi_j\nonumber
\end{eqnarray}
defines a semiclassically bounded family in Boutet de Monvel's 
calculus of order $-2m$, and 
$$(1+A_\hbar A_\hbar^*)^{-1}= \widetilde B(\hbar)+\hbar R_4(\hbar);
$$
with a family  $R_4$ which is uniformly bounded on $L^2(X,E_1)$.
}

%%%%%%%%%%%%%%%%%%%%%%%%%%%%%%%%%%%%%%%%%%%%%%%%%%%%%%%%%%%%%%%%%%%%%%%

\dfn{t.6}{The construction of $A_\hbar$ together with 
Propositions \ref{t.4}
allows us to define a section of the continuous field 
$M_N(C^*_r(\cT^- X)^\sim)$ by 
$$s(\hbar)=\left\{\begin{array}{ll}\Gr(A_\hbar),&0<\hbar\le 1\\
\Gr(p)\oplus \Gr(c),&\hbar=0\end{array}\right..$$
 }\medskip
We shall now show that this section is continuous.
We will distinguish the cases $\hbar>0$ and $\hbar =0$.

%%%%%%%%%%%%%%%%%%%%%%%%%%%%%%%%%%%%%%%%%%%%%%%%%%%%%%%%%%%%%%%%%%%%%%%

\prop{t.7}{The section $s$ is continuous on $(0,1]$.}

%%%%%%%%%%%%%%%%%%%%%%%%%%%%%%%%%%%%%%%%%%%%%%%%%%%%%%%%%%%%%%%%%%%%%%%

\Proof. It follows from \ref{1.15}(a) and the fact that the symbol 
topology is stronger than the operator topology  that the mappings 
\begin{eqnarray*}
\hbar\mapsto A_\hbar&\in&\cL(H^m(X,E_1),L^2(X,E_2))\quad\text{and}\\
\hbar\mapsto A_\hbar&\in& \cL(L^2(X,E_1),H^{-m}_0(X,E_2))
\end{eqnarray*}
depend continuously on $\hbar$.
As taking adjoints  and inversion are continuous, we obtain the 
assertion.
\eproof

%%%%%%%%%%%%%%%%%%%%%%%%%%%%%%%%%%%%%%%%%%%%%%%%%%%%%%%%%%%%%%%%%%%%%%%

\thm{t.8}{The section $s$ is continuous in  $\hbar=0$.}

%%%%%%%%%%%%%%%%%%%%%%%%%%%%%%%%%%%%%%%%%%%%%%%%%%%%%%%%%%%%%%%%%%%%%%%
 
\Proof. Consider the four entries of $\Gr(A_\hbar)$. 
By Corollary  \ref{c.3v}, 
$(1+A_\hbar^*A_\hbar)^{-1}$ differs from  $B(\hbar)$
by a term which vanishes in $\hbar=0$. 
It is therefore enough to show the continuity of $B(\hbar)$. 
By construction, the boundary symbol of $B(\hbar)$ is given 
by $(1+c^*c)_\hbar^{-1}$ while the interior symbol is 
$(1+p^* p)_\hbar^{-1}$.

By Proposition \ref{h.7a} %and Remark \ref{h.7f}
we find smoothing symbols $q_k$
converging to $(1+p^*p)^{-1}$  in the topology of $S^0$.
Since we assumed $p$ to be constant near $\partial X$, we may 
assume the same of the $q_k$. 

As for the boundary symbol, we know that $(1+c^*c)^{-1}$ 
is a boundary symbol operator in Boutet de Monvel's 
calculus whose pseudodifferential part is 
$r^+(1+p^*p)|_{x_n=0}^{-1}(D_n)$. 
Denote by $h$ its singular Green part. 
According to Proposition \ref{h.7a} we find a sequence of symbols 
$h_k$ in 
$S^{-\infty}(\R^{n-1}\times\R^{n-1};\cS'(\R_+),\cS(\R_+))$
converging to $h$ in the topology of 
$S^{0}(\R^{n-1}\times\R^{n-1};L^2(\R_+),L^2(\R_+))$.
Replacing, for $\hbar=1$, in the definition of $B(1)$ the 
pseudodifferential symbols over interior charts by $q_k$ and
the boundary symbol operators over boundary charts by 
$r^+q_k|_{x_n=0}(D_n) +h_k $, we obtain a sequence of operators 
$B_k$. 
According to Remark \ref{h.7f}, a further approximation allows us
to assume that $B_k$ are integral operators with 
smooth compactly supported integral kernels.
%Note that the approximation is uniform for $0<\hbar\le 1$
%by \ref{1.15f} and the corresponding result for 
%operators with $\hbar$-scaled kernels. 
Similarly we approximate the other entries
in $\cG(A)-e$. Adding then $e$ again and going over to 
$\hbar$-scaled symbols, we infer from Lemma \ref{1.15f} that 
the approximation is uniform for $0<\hbar\le 1$. 
Hence we obtain a sequence of sections $s_k(\hbar)$ of 
$M_N(C^*_r(\ccT X))$ which approximates $\Gr(A_\hbar)$ uniformly.

By  definition, these sections are continuous for $\hbar>0$
and have a continuous extension to $\hbar=0$ given by the
$N\times N$-matrices of their interior and boundary symbols. 
As this matrix, on the other
hand, tends to $\cG(p)\oplus\cG(c)$, we conclude that $s$ is 
continuous.\eproof

%%%%%%%%%%%%%%%%%%%%%%%%%%%%%%%%%%%%%%%%%%%%%%%%%%%%%%%%%%%%%%%%%%%%%%%

\section{The Fundamental Class for Manifolds with Boundary}
In this section we will describe how the fundamental class
$$\int_{T^*X^\circ}:
H^*_c(T^*X^\circ)=HP^*(\cC^\infty_c(T^*X^\circ)\to \C$$
extends to a fundamental class
$$F:HP^*(\cC^\infty_{tc} (T^-X))\to \C.$$
Let us first assume $T^*X=T^*\partial X \times T^*\R_+$ so that
we can consider the elements of
$\cC^\infty_c(T^-X)$ as elements of 
$\cC^\infty_c(T^*\partial X)\pitensor \cC^\infty_{tc}(T^-\R_+).$
We write an element of $\cC^\infty_{tc}(T^-\R_+)$ as the sum of a 
pseudodifferential symbol $p$ and a singular Green symbol $g$ 
on the boundary.
We then obtain the operator of multiplication by 
$p$ on $T^*\R_+$ and the boundary symbol %(Wiener-Hopf) 
operator 
\begin{eqnarray}\label{6.0.1}
c =p(0,D_n)+g.
\end{eqnarray}
%induced by  $p+g$. 
Following Fedosov we define 
$$\tr'(p+g)=\tr(g),$$
noting that $g$ is an integral operator on $L^2(\R_+)$
with a rapidly decreasing kernel and thus trace class.  
%with the mapping $\Pi':H\to \C$ given by 
%$$\Pi'(f)= $$ 
The functional $\tr'$ is not quite a trace, 
but satisfies 
%-- with the notation introduced in 
%\eqref{6.0.1} -- 
the following fundamental property 
\cite[(2.19)]{Fedosov96}
\begin{eqnarray*} \label{fed} 
{\tr'([p_1+g_1,p_2+g_2])}
&=&
%-i\Pi'
%\frac1{2\pi i}
-i\int \frac{\partial p_1(0,\xi_n)}{\partial \xi_n }p_2(0,\xi_n )d\xi_n=
%i\Pi' 
%-\frac1{2\pi i}
i\int p_1(0,\xi_n )\frac{\partial p_2(0,\xi_n )}{\partial \xi_n }d\xi_n.
\end{eqnarray*}

Given  an element in the cyclic periodic complex, 
i.e. an element in $\cC^\infty_{tc} (T^-X)^{\otimes^{m+1}}$, 
we first introduce the boundary functional
\begin{eqnarray*}
\lefteqn{F_\partial (( f_0\otimes (p_0+g_0))
\otimes (f_1\otimes (p_1+g_1))\otimes \cdots 
\otimes (f_m\otimes (p_m+g_m)))}\\
&=&
\int_{T^*\partial X}f_0\,df_1\cdots df_m
\tr'((p_0+g_0)\cdots (p_m+g_m)).
\end{eqnarray*}
We have used here the splitting of $\cC^\infty_{tc}(T^-X)$ 
into tensor factors. 
For notational convenience we will omit the tensor symbols 
and write $f_j(p_j+g_j)$ 
instead of $f_j\otimes (p_j+g_j)$.

The fundamental class is given by
\begin{eqnarray*}\lefteqn{F(f_0(p_0+g_0)\otimes 
f_1 (p_1+g_1)\otimes \cdots \otimes f_m (p_m+g_m))
=\int_{T^*X}f_0 p_0d(f_1p_1)\cdots d(f_m p_m)}\\
&&+iF_\partial \bigg( \sum_{\sigma \text{ cyclic}} 
\hbox{sgn}(\sigma)f_{\sigma (0)}  (p_{\sigma (0)}+g_{\sigma (0)})
%\otimes  f_{\sigma (1)} (p_{\sigma (1)}+g_{\sigma (1)})
\otimes \cdots 
\otimes f_{\sigma (m)} (p_{\sigma (m)}+g_{\sigma (m)})\bigg).
\end{eqnarray*}

\begin{prop}{cykel}{The fundamental class $F$ is a cocycle 
on the periodic cyclic complex.}
\end{prop}

\Proof. We need to prove that $F((B+b)\underline{a})=0$, 
where $\underline{a}\in CC^{per}_*(\cC^\infty_{tc}(T^-X))$.
By Stokes' theorem the boundary part $F_\partial$ of $F$ vanishes on 
$B\underline{a}$. 
The remaining `nonboundary' part of $F$ clearly vanishes on 
$b\underline{a}$. 
Computing the nonboundary part we get
$$\int_{T^*X}d(f_0 p_0)d(f_1 p_1)\cdots d(f_m p_m)
=\int_{\partial (T^*X)}f_0 p_0\,d(f_1 p_1 )\cdots d(f_m p_m)
$$
%Computing the boundary part we get, before symmetrizing
%\begin{eqnarray*}\lefteqn{F_\partial (b(\underline{a}))}\\
%&=&-i \int_{T^*\partial X} \Big( tr'((p_0+g_0)\cdots (p_n+g_n))\cdot f_0f_1df_2 \cdots df_n\\
%&+&\sum_{j=1}^{n-1}(-1)^jtr'((p_0+g_0)\cdots (p_n+g_n))\cdot f_0df_1\cdots d(f_jf_{j+1})\cdots df_n\\
%&+&(-1)^n tr'((p_n+g_n)(p_0+g_0)\cdots (p_{n-1}+g_{n-1})\cdot f_nf_0df_1\cdots df_{n-1}\Big)\\
%&=&-i\int_{T^*\partial X} \Big( tr'((p_0+g_0)\cdots (p_n+g_n))\\
%&-&tr'((p_n+g_n)(p_0+g_0)\cdots (p_{n-1}+g_{n-1}))\Big)\cdot f_nf_0df_1\cdots df_{n-1}\\
%&=&-i\int_{T^*\partial X} \Pi'(p_0(0,\xi) \cdots p_{n-1}(0,\xi) \frac{\partial p_n(0,\xi)}{\partial \xi})\cdot f_nf_0df_1\cdots df_{n-1}\\
%&=&\int_{T^*\partial X}\int_\R p_0(0,\xi) \cdots p_{n-1}(0,\xi) \frac{\partial p_n(0,\xi)}{\partial \xi}d\xi \cdot f_nf_0df_1\cdots df_{n-1}
%\end{eqnarray*}
We want to compute the boundary part, 
i.e. $F_\partial$, on the cyclic permuted terms appearing in $F$.  
A single cyclic permutation of $b(\underline{a} )$ without $\hbox{sgn}(\sigma )$ is of the form 
\begin{eqnarray*}
&&f_{i+1}(p_{i+1}+g_{i+1})\otimes \cdots \otimes f_0(p_0+g_0) f_1(p_1+g_1)
   \otimes \cdots \otimes f_{i}(p_{i}+g_{i})\\
&&-f_{i+1}(p_{i+1}+g_{i+1})\otimes \cdots \otimes f_1(p_1+g_1) f_2(p_2+g_2)
   \otimes \cdots \otimes f_{i}(p_{i}+g_{i})\\
&&\qquad \qquad \qquad \qquad \vdots\\
&&+(-1)^{i-1}f_{i+1}(p_{i+1}+g_{i+1})\otimes \cdots \otimes 
   f_0(p_0+g_0) \otimes \cdots \otimes f_{i-1}(p_{i-1}+g_{i-1})f_{i}(p_{i}+g_{i})\\
&&+(-1)^if_i(p_i+g_i)f_{i+1}(p_{i+1}+g_{i+1})\otimes \cdots \otimes 
   f_0(p_0+g_0) \otimes \cdots \otimes f_{i-1}(p_{i-1}+g_{i-1})\\
&&\qquad \qquad \qquad \qquad \vdots \\
&&+(-1)^mf_i(p_i+g_i)\otimes  \cdots \otimes 
   f_m(p_m+g_m)f_0(p_0+g_0) \otimes \cdots \otimes f_{i-1}(p_{i-1}+g_{i-1}).
\end{eqnarray*}
We split this expression into the sum of the first $i$ terms and the sum of the 
subsequent $m+1-i$ terms.
The action of $F_\partial$ on the first $i$ terms is
\begin{eqnarray*}
\forget{
&&F_\partial (f_{i+1}(p_{i+1}+g_{i+1})\otimes \cdots \otimes f_0(p_0+g_0) f_1(p_1+g_1)\otimes \cdots \otimes f_{i}(p_{i}+g_{i}))\\
&&-F_\partial (f_{i+1}(p_{i+1}+g_{i+1})\otimes \cdots \otimes f_1(p_1+g_1) f_2(p_2+g_2)\otimes \cdots \otimes f_{i}(p_{i}+g_{i}))\\
&&\vdots\\
&&+(-1)^{i-1}F_\partial (f_{i+1}(p_{i+1}+g_{i+1})\otimes \cdots \otimes f_0(p_0+g_0) \otimes \cdots \otimes f_{i-1}(p_{i-1}+g_{i-1})f_{i}(p_{i}+g_{i}))\\
&=&}%forget
\lefteqn{\tr'\big((p_{i+1}+g_{i+1})(p_{i+2}+g_{i+2})\cdots (p_m+g_m)(p_0+g_0)\cdots (p_i+g_i)\big)} \\
&&\times\int_{T^*\partial X}\Big( f_{i+1}df_{i+2}\cdots d(f_0f_1)\cdots df_i-f_{i+1}df_{i+2}\cdots d(f_1f_2)\cdots df_i+\ldots\\ 
&&\mbox{\ \ }+\ldots(-1)^{i-1}f_{i+1}df_{i+2}\cdots d(f_{i-1}f_i)\Big) .
\end{eqnarray*}
The second factor in this expression can be rewritten as
\begin{eqnarray}
\forget{
\lefteqn{\int_{T^*\partial X}\Big( f_{i+1}df_{i+2}\cdots d(f_0f_1)\cdots 
  df_i-f_{i+1}df_{i+2}\cdots d(f_1f_2)\cdots df_i}\\
&&+\ldots(-1)^{i-1}f_{i+1}df_{i+2}
  \cdots d(f_{i-1}f_i)\Big)\\
&=&
}%forget
\lefteqn{\int_{T^*\partial X}\Big( f_0f_{i+1}df_{i+2}\cdots df_m df_1\cdots 
  df_i+f_1f_{i+1}df_{i+2}\cdots df_0df_2\cdots df_i}\nonumber\\
&&  -\big( f_1f_{i+1}df_{i+2}\cdots df_0df_2\cdots df_i
  +f_2f_{i+1}df_{i+2}\cdots df_1df_3\cdots df_i\big)\nonumber\\
&&+\ldots +(-1)^{i-1}\big( f_{i-1}f_{i+1}df_{i+2}\cdots df_{i-2}df_i
  +f_if_{i+1}df_{i+2}\cdots df_{i-1}\big)\Big)\nonumber\\
&=&\int_{T^*\partial X}\Big( f_0f_{i+1}df_{i+2}\cdots df_mdf_1\cdots df_i
  +(-1)^{i-1}f_if_{i+1}df_{i+2}\cdots df_{i-1}\Big).\label{second}
\end{eqnarray}
A short computation shows that
\begin{eqnarray*}
\lefteqn{d(f_0f_if_{i+1})df_{i+2}\cdots df_m df_1\cdots df_{i-1}
=f_0f_idf_{i+1}\cdots df_m df_1\cdots df_{i-1}}\\
&&+(-1)^{m}f_0f_{i+1}df_{i+2}\cdots df_m df_1\cdots df_i 
  +(-1)^{m-(i-1)}f_if_{i+1}df_{i+2}\cdots df_{i-1}.
\end{eqnarray*}
Stokes' theorem then implies that the term \eqref{second} equals
\begin{eqnarray*}
\forget{
\lefteqn{\int_{T^*\partial X}\Big( f_0f_{i+1}df_{i+2}\cdots df_m df_1\cdots df_i
  +(-1)^{i-1}f_if_{i+1}df_{i+2}\cdots df_{i-1}\Big)}\\
&=&
}%forget
\lefteqn{(-1)^{m+1}\int_{T^*\partial X}f_0f_idf_{i+1}\cdots df_m df_1\cdots df_{i-1}}\\
&=&(-1)^{mi+1}\int_{T^*\partial X}f_0f_idf_1\cdots df_{i-1}df_{i+1}\cdots df_m.
\end{eqnarray*}
Putting everything together, the action of $F_\partial$ on the first sum is given by
\begin{eqnarray*}
&(-1)^{mi+1}&\tr'\big((p_{i+1}+g_{i+1})(p_{i+2}+g_{i+2})\cdots 
(p_m+g_m)(p_0+g_0)\cdots (p_i+g_i)\big) \\
&&\times \ \int_{T^*\partial X}f_0f_idf_1\cdots df_{i-1}df_{i+1}\cdots df_m.
\end{eqnarray*}
The action of $F_\partial$ on the remaining $m+1-i$ terms gives
\begin{eqnarray*}
\forget{
&&(-1)^iF_\partial \big(f_i(p_i+g_i)f_{i+1}(p_{i+1}+g_{i+1})\otimes \cdots \otimes 
f_0(p_0+g_0) \otimes \cdots \otimes f_{i-1}(p_{i-1}+g_{i-1})\big)+\\
&&\ldots +(-1)^mF_\partial \big(f_i(p_i+g_i)
%\otimes f_{i+1}(p_{i+1}+g_{i+1})
\otimes \cdots 
  \otimes f_m(p_m+g_m)f_0(p_0+g_0) \otimes \cdots \otimes f_{i-1}(p_{i-1}+g_{i-1})\big)\\
&=&
}%forget
\lefteqn{\tr'\big( (p_i+g_i)(p_{i+1}+g_{i+1})\cdots (p_m+g_m)(p_0+g_0)\cdots (p_{i-1}+g_{i-1})\big)}  \\
&&\quad\times\int_{T^*\partial X}\Big( (-1)^i f_if_{i+1}df_{i+2}\cdots df_mdf_0\cdots df_{i-1}+\\
&&\ \ \qquad(-1)^{i+1}
  \big( f_if_{i+1}df_{i+2}\cdots df_{i-1}+f_if_{i+2}df_{i+1}df_{i+3}\cdots df_{i-1}\big)+\ldots \\
&&\ \qquad+(-1)^m\big(f_if_mdf_{i+2}\cdots df_{m-1}df_0\cdots df_{i-1}+f_if_0df_{i+2}
  \cdots df_mdf_1 \cdots df_{i-1}   \big) \\
&=&(-1)^{mi}\tr'\Big( (p_i+g_i)(p_{i+1}+g_{i+1})\cdots (p_m+g_m)(p_0+g_0)\cdots (p_{i-1}+g_{i-1})\Big) \\
&&\qquad\times\int_{T^*\partial X} f_if_0df_1\cdots df_{i-1}df_{i+1}\cdots df_m.
\end{eqnarray*} 
All in all we get that the action of $F_\partial$ on this symmetrization is 
\begin{eqnarray*}
\lefteqn{(-1)^{mi}\tr' ([p_i+g_i, (p_{i+1}+g_{i+1})\cdots 
   (p_m+g_m)(p_0+g_0)\cdots (p_{i-1}+g_{i-1})])}\\
&\times& \int_{T^*\partial X} f_if_0df_1\cdots df_{i-1}df_{i+1}\cdots df_m \\\\
&=&(-1)^{mi+1}\ i\int_{T^*\partial X}\int_\R p_0(0,\xi_n) 
   \cdots p_{i-1}(0,\xi_n )p_{i+1}(0,\xi_n )\cdots p_{m}(0,\xi_n )\ \frac{\partial p_i(0,\xi_n )}{\partial \xi_n }d\xi_n \\
&&\qquad \times f_if_0df_1\cdots df_{i-1}df_{i+1}\cdots df_m\\
&=&(-1)^{(m+1)i}\ i\int_{\partial (T^* X)} p_0(0,\xi_n ) \cdots p_{i-1}(0,\xi_n )
   p_{i+1}(0,\xi_n )\cdots p_{m}(0,\xi_n )\ f_if_0df_1\cdots df_{i-1}\\ 
&& \qquad  \times \frac{\partial p_i(0,\xi_n )}{\partial \xi_n }d\xi_n df_{i+1}\cdots df_m\\
&=&{\sgn}(\sigma)\ i\int_{\partial (T^* X)} p_0(0,\xi_n ) \cdots p_{i-1}(0,\xi_n  )p_{i+1}(0,\xi_n )\cdots p_{m}(0,\xi_n ) \ f_if_0df_1\cdots df_{i-1}\\
&&\qquad\times\frac{\partial p_i(0,\xi_n )}{\partial \xi_n }d\xi_n df_{i+1}\cdots df_m,
\end{eqnarray*}
where $\sigma$ is the corresponding permutation. Hence
\begin{eqnarray*}
F(b(\underline{a}))
&=&\sum_{i=1}^m\int_{\partial (T^* X)} p_0(0,\xi_n ) 
\cdots p_{i-1}(0,\xi_n )p_{i+1}(0,\xi_n )\cdots p_{m}(0,\xi_n )  \  
f_if_0df_1\cdots df_{i-1}\\
&&\qquad\times \frac{\partial p_i(0,\xi_n )}{\partial \xi_n }d\xi_n df_{i+1}\cdots df_m,
\end{eqnarray*}
%Symmetrizing we get
%\begin{eqnarray*}
%\sum_{j=1}^n (-1)^{j-1}\int_{T^*\partial X}\int_\R p_0(0,\xi) \cdots p_{j-1}(0,\xi)p_{j+1}(0,\xi)\cdots  p_n(0,\xi)  \cdot f_jf_0df_1\cdots df_{j-1}\frac{\partial p_j(0,\xi)}{\partial \xi}d\xi df_{j+1}\cdots df_n
%\end{eqnarray*}
which is equal to 
$$\int_{\partial (T^*X)} f_0 p_0d(f_1 p_1)\cdots d(f_m p_m).$$
\eproof

\extra{genfund}{The fundamental class in general}%
{For general $X$ the restriction of an element $a$ in 
$\cC^\infty_{tc}(T^-X)$ to $\partial (T^*X)$ can be factorized as 
a sum of elements of the form $f\otimes (p+g)$, i.e. the boundary symbol factorizes. 
We will adopt this notation, i.e. the boundary part of $a$ will be denoted $f\otimes (p+g)$. 
The symbol part in the interior will be denoted by 
$\tilde{a}$, which is a function on $T^*X$. 
It is then straightforward to generalize the fundamental class to nonproduct cases:

Let $\omega$ be a closed differential form on $T^*X$ 
of even degree, which is the pull back of a closed differential 
form on $X$.  We first define $F_{\partial ,\omega}$ by
\begin{eqnarray*}\lefteqn{F_{\partial,\omega} 
( f_0\otimes (p_0+g_0))\otimes (f_1\otimes (p_1+g_1)\otimes \cdots \otimes (f_m\otimes (p_m+g_m)))}\\
&=&\int_{T^*\partial X}f_0df_1\cdots df_m\cdot
\omega \cdot \tr'((p_0+g_0)\cdots (p_m+g_m))
\end{eqnarray*}
and then let 
\begin{eqnarray*}\lefteqn{F(a_0\otimes a_1\otimes \cdots \otimes a_m)}\\
&=&\int_{T^*X}\tilde{a}_0d\tilde{a}_1\cdots \tilde{a}_m \cdot \omega
-iF_{\partial,\omega} \left( \sum_{\sigma \text{ cyclic}} 
\sgn(\sigma)a_0
\otimes  a_{\sigma (1)} \otimes\cdots \otimes a_{\sigma (m)} \right).
\end{eqnarray*}}

\begin{prop}{fundament}{$F_{\omega}$ descends to a map 
$$F_\omega :HP(\cC^\infty_{tc}(T^-X))\to \C.$$}
\end{prop}
\Proof. The same computation as in Proposition \ref{cykel}. \eproof

\section*{The Index Formula}

Like in \cite[Section 2.5]{Connes} the short exact sequence
$$0\to \cC_0(]0,1])\otimes \cK \to C^*_r(\ccT^- X)\to C^*_r(T^-X)\to 0$$
induces an analytic index map
%ES C_r ??
$$\ind_a:K_0(C_r^*(T^-X))\to \Z.$$
In this section we will  give a formula for this index map. 

We have another short exact sequence coming from the interior of the manifold, namely
 $$0\to \cC_0(]0,1])\otimes \cK \to C^*_r( \ccT X^\circ)\to \cC_0(T^*X^\circ)\to 0$$
(noting that $C^*_r(TX^\circ)\cong \cC_0(T^*X^\circ))$ also inducing an analytic index map
$$\ind_a:K_0(\cC_0(T^*X^\circ))\to \Z.$$
According to Connes \cite[Section 2.5]{Connes}  we have in this case
$$\ind_a=\ind_t,$$
where $\ind_t$ denotes the topological index. 
On the other hand we have an isomorphism 
$\Phi :K_0(C^*_r(T^-X))\to K_0(\cC_0(T^*X^\circ))$ and the diagram
$$
\begin{array}{ccc}
K_0(C^*_r(T^-X))&\to& \Z\\
\downarrow  &\nearrow&\\
 K_0(\cC_0(T^*X^\circ))
\end{array}
$$
commutes. 

%ES Johannes + cosmetic
Let $T$ be the Fredholm operator of order and class zero in Boutet de Monvel's
calculus introduced in \ref{setup}, 
and let $\Gr(a)=\Gr(p)\oplus\Gr (c)$ denote the graph projection of a complete 
symbol of $\gL^{m,+}_-T$. From Theorem \ref{t.8} we obtain
\begin{thm}{absind}{
The index of $T$ is given by 
$$\ind T =\ind_a([\cG (a)]-[e]).%=\dim (\ker A)-\dim (\coker A).
$$}
\end{thm}
We define the topological index
$$\hbox{ind}_t:K_0(C^*_r (T^-X))\to \Z$$
 as the composition $\ind_t \circ \Phi$. We thus get an index theorem
\begin{thm}{indtop}{
$\ind_a=\ind_t.$
}\end{thm}
With this notation we can now prove
\begin{thm}{ind}{%ES e subtracted
$$\ind T = 
F_{Td(X)}(\ch ([\cG (a)]-[e])),$$
where $\ch$ denotes the the Chern-Connes character.}
\end{thm}

\Proof. Let $C^n(T^-X)$ denote the subalgebra of $C^*_r(T^-X)$ 
consisting of symbols of order strictly less than $-n$, 
$n$ being the dimension of $X$. 
We note that  $C^n(T^-X)^\sim$ is closed under holomorphic 
functional calculus and hence 
$K_0( C^n(T^-X))=K_0(C^*_r(T^-X))$. 
Also note that $F_\omega$ is defined on $HP(C^n(T^-X))$. 
Using the cohomological form of the topological index we 
get the following commutative diagrams:
$$ 
\xymatrix{K_0(C^n(T^*X^\circ)) \ar@{=}[r] \ar[d]^{\hbox{ch}} & K_0(C^n{T^-(X)})\ar[d]^{\hbox{ch}}\\
HP_{ev}(C^n(T^*X^\circ)) \ar[r] \ar@{=}[d]& HP_{ev} (C^n(T^-X))\ar[d]^{F_{Td(X)}} \cdot \\
H_c^*(T^*X^\circ) \ar[r]^{\int Td(X)\cdot } & \mathbb{Z}
}
$$
and 
$$ \xymatrix{K_0(C^n(T^*X^\circ)) \ar[dd]_{\hbox{ch}}\ar@{=}[r]&K_0 (\mathcal{C}_0(T^*X^\circ ))\ar@{=}[r] \ar[ddr]^{\hbox{ind}_t}& K_0 (\mathcal{C}_0(T^*X^\circ )) \ar[dd]^{\hbox{ind}_a} \\
\\
HP_{ev}(C^n(T^*X^\circ)) \ar@{=}[r] & H_c^*(T^*X^\circ )\ar[r]_{\int Td(X) \cdot }& \mathbb{Z} }$$
and
$$ \xymatrix{ K_0(\mathcal{C}_0(T^*X^\circ)) \ar[ddr]_{\hbox{ind}_a} \ar@{=}[r]&K_0(C^*_r(T^-X)) \ar[dd]_{\hbox{ind}_a} \ar@{=}[r]&K_0(C^n(T^-X)) \ar[ddl]^{\hbox{ind}_a}\\
\\
 &\mathbb{Z} }$$
from which follows that $F_{Td(X)} \circ \ch=\ind_a$ on $K_0(C^n(T^-X))$.
\eproof

\end{document}